\documentclass[11pt,twoside]{article}
\usepackage{fancyhdr}
\usepackage[colorlinks,citecolor=blue,urlcolor=black,linkcolor=blue,bookmarks=false]{hyperref}
\usepackage{amsfonts,epsfig,graphicx}
\usepackage{afterpage}
\usepackage{comment}
\usepackage{amsmath,amssymb,amsthm} 
\usepackage{array}
\usepackage{dsfont}
\usepackage{fullpage}
\usepackage[T1]{fontenc} 
\usepackage{epsf} 
\usepackage{graphics} 
\usepackage{amsfonts,amsmath}
\usepackage[sort]{natbib} 
\usepackage{psfrag,xspace}
\usepackage{color,etoolbox}
\usepackage{subcaption}

\def\ind{\perp\!\!\!\perp}

\newcommand{\Pb}{\mathbb{P}}
\newcommand{\Pn}{\mathbb{P}_n}

\newcommand{\E}{\mathbb{E}}
\newcommand{\R}{\mathbb{R}}

\newcommand{\norm}[1]{\lVert#1\rVert}
\newcommand{\anorm}[1]{\left\lVert#1\right\rVert}
\newcommand{\vertline}{\hspace{0.1cm}\middle|\hspace{0.1cm}}

\def\logit{\text{logit}}
\def\expit{\text{expit}}

\DeclareMathOperator*{\argmin}{arg\,min}

\DeclareSymbolFont{bbold}{U}{bbold}{m}{n}
\DeclareSymbolFontAlphabet{\mathbbold}{bbold}
\newcommand{\one}{\mathbbold{1}}
\newtheorem{theorem}{Theorem}
\newtheorem{lemma}{Lemma}

\newtheorem{proposition}{Proposition}

\theoremstyle{definition}

\theoremstyle{remark}

\newtheorem{remark}{Remark}

\begin{document}

\def\spacingset#1{\renewcommand{\baselinestretch}%
{#1}\small\normalsize} \spacingset{1}

\raggedbottom
\allowdisplaybreaks[1]


  \title{\vspace*{-.4in} {Causal Inference with High-Dimensional Treatments}}
  \author{\\ Patrick Kramer$^1$, Edward H. Kennedy$^1$, Isaac M. Opper$^2$ \\ \\ \\
    $^1$Department of Statistics \& Data Science, \\
    Carnegie Mellon University \\ 
	\\
   $^2$Luskin School of Public Affairs,  \\
   University of California, Los Angeles \\ \\ \\
    \texttt{pkramer@andrew.cmu.edu, edward@stat.cmu.edu, iopper@ucla.edu} \\
\date{}
    }

  \maketitle
  \thispagestyle{empty}

\begin{abstract}
    In this work, we consider causal inference in various high-dimensional treatment settings, including for single multi-valued treatments and vector treatments with binary or continuous components, when the number of treatments can be comparable to or even larger than the number of observations. These settings bring unique challenges: first, the treatment effects of interest are represented by a high-dimensional vector rather than a scalar; second, positivity violations are often unavoidable; and third, estimation can be based on a smaller effective sample size. We first discuss fundamental limits of estimating effects here, showing that consistent estimation is impossible without further assumptions. We go on to propose novel doubly robust estimators for mean potential outcomes of a high-dimensional single multi-valued treatment. We analyze the proposed estimators under sparsity assumptions, giving finite-sample risk bounds and showing that consistent estimation is possible under these conditions. Moreover, we derive minimax lower bounds in a sparse and structure-agnostic model to characterize optimal rates of convergence and show our risk bounds are unimprovable. We then generalize our proposed estimators as a sparse pseudo-outcome regression framework with constrained regression estimators and error guarantees under sparsity, allowing estimation of generic functionals and different types of high-dimensional treatments. We apply the framework to derive estimators of the mean potential outcomes for high-dimensional vector treatments. Finally, we illustrate the proposed methods through a simulation and an empirical application.
\end{abstract}


\section{Introduction}

Many causal inference problems involve a large number of treatments relative to the number of observations, which we refer to as a \textit{high-dimensional} treatment. For example, high-dimensional treatments naturally occur in medicine, e.g., for cancer radiation therapy, where the treatment consists of different possible locations and intensities (\citet{nabi2022semiparametric}) or when comparing different healthcare providers (\citet{susmann2024doublyrobust}). Other examples of high-dimensional treatments include genetic variants and environmental exposures (\citet{mitra2022future}), text and image embeddings (\citet{feder2022causal}), marketing campaigns (\citet{sharma2020hi}), or schools in a dataset containing information about student performance. We can study the effects of such high-dimensional treatments by estimating the mean potential outcomes for each possible treatment level, i.e., the outcome that we would observe if we intervened and applied this particular treatment level to all units. Comparing these mean potential outcomes across different treatment levels enables a comparison of the effectiveness of these treatments, for example, allowing for the identification of the best possible treatment option.\\

There are unique challenges when estimating the effects of high-dimensional treatments, where the number of treatments $k$ is large relative to the sample size $n$. These include an increasingly complex target parameter, inevitable positivity violations when the treatment levels are discrete, and potentially smaller effective sample sizes due to fewer observations at each treatment level. For these reasons, consistent estimation and fast convergence rates are impossible in this setting without further structure. This naturally suggests a sparsity-based approach to enable efficient estimation, in particular, using estimators that have the flavor of a Lasso or best subset selection estimator \citep{tibshirani1996regression, chatterjee2014assumptionless, hastie2020best, wainwright2019highdimensional, rigollet2023high}.

\subsection{Our Contributions}

Our work addresses the fundamental challenges of high-dimensional treatments by proposing efficient estimators for mean potential outcomes in various high-dimensional treatment settings, in combination with provable theoretical guarantees and optimality results. In particular, this paper presents the following main contributions:
\begin{itemize}
    \item In Section \ref{section:single-treatment}, we \textit{propose efficient estimators} of mean potential outcomes for high-dimensional single multivalued treatments.  We \textit{prove double robustness} and \textit{provide error guarantees} under sparsity in Theorem~\ref{thm:single-treatment}. 
    
    \item In Section \ref{section:high-dim-trt}, we \textit{illustrate fundamental limits of estimating mean potential outcomes} in the high-dimensional treatment regime when no further structure is assumed, which motivates incorporating sparsity. We go on to \textit{prove minimax optimality} of our proposed estimators by deriving minimax lower bounds under a sparse and structure-agnostic model in Theorems~\ref{thm:minimax-lb-sparse-regime} and \ref{thm:minimax-lb-structure-agnostic}, which is, to the best of our knowledge, the first lower bound in the structure-agnostic model for a high-dimensional target.
    \item In Section \ref{section:extensions}, we first generalize the results of Section \ref{section:single-treatment} to general target parameters and other types of high-dimensional treatments. In Section \ref{section:pseudo-outcome-regression}, we \textit{propose a novel sparse pseudo-outcome regression framework} for estimation of general high-dimensional statistical functionals. We \textit{propose efficient estimators} and \textit{derive fast convergence rates} under sparsity assumptions in Theorem \ref{thm:master}. Then, in Section \ref{section:vector-treatment}, we apply the proposed general framework to \textit{estimate mean potential outcomes for high-dimensional vector treatments}. We \textit{present efficient doubly robust estimators} and \textit{provide error guarantees} for binary vector treatments in Theorem~\ref{thm:binary-vector-trt} and for continuous vector treatments in Theorem~\ref{thm:continuous-vector-trt}.
    \item In Section \ref{section:empirical-data-analysis}, we demonstrate our proposed estimators of Section \ref{section:single-treatment} by applying them to simulated and real data. In Section \ref{section:simulation}, we \textit{conduct a simulation} to illustrate the theoretical guarantees and performance of our proposed estimators on generated data with a high-dimensional treatment in both a sparse and non-sparse setting. In Section \ref{sec:data-analysis}, we \textit{apply our proposed method} to estimate the effect of schools on standardized test scores using data from the North Carolina Education Research Data Center.
\end{itemize}
Section \ref{section:discussion} then concludes with a discussion of the presented results and provides possible directions for future work.

\subsection{Related work}

Efficient estimation of mean potential outcomes has been extensively studied in the literature for high-dimensional covariates (\citet{belloni2014inference}, \citet{farrell2015robust}, \citet{d2017overlap}, \citet{athey2018approximate},  \citet{bradic2019sparsity},  \citet{smucler2019unifying}, \citet{wang2020debiased}, \citet{jiang2022new}, \citet{chernozhukov2022debiased}, \citet{liu2023root}, \citet{zeng2024causal}), as well as for high-dimensional outcomes (\citet{du2025causal}), and continuous treatments (see \citet{diaz2013targeted}, \citet{kennedy2017nonparametric}, \citet{bonvini2022fast}, \citet{schindl2024incremental} and references therein). In comparison, high-dimensional treatments have received less attention.\\

Some recent work on high-dimensional treatments has focused on dimension reduction using parametric models \citep{nabi2022semiparametric, goplerud2022estimating, andreu2024contrastive, lin2025learning}. The idea is to find a lower-dimensional representation of the high-dimensional treatment that preserves the causal relationship. Another approach in the literature is a plug-in regression estimator based on a deep neural network, as proposed by \citet{sharma2020hi}. While these approaches address high-dimensionality and are validated with experiments, we aim to contribute an alternative sparsity-based approach that maintains the full treatment structure, achieves provable theoretical guarantees, and is even minimax optimal.\\

\citet{susmann2024doublyrobust} study treatment-group-specific effects, including both direct and indirect standardization parameters, motivated by the application of provider profiling. They allow for multiple treatment levels and derive doubly robust efficient estimators for the target parameters in this setting. Moreover, \citet{mcclean2025comparing} study the comparison of causal parameters with many treatments under potential positivity violations. They present a criterion for comparability of causal parameters across different treatment levels, provide examples of parameters that satisfy this criterion, and propose doubly robust estimators for those parameters. In another recent paper, \citet{xiang2025double} consider the estimation of mean potential outcomes in various treatment regimes that allow for multiple treatment levels, including single multi-valued and vector treatments. They propose doubly robust estimators and show asymptotic normality. In all the aforementioned papers, the estimators are analyzed under the regime in which the number of treatments $k$ is treated as fixed and $n\to\infty$. However, such error guarantees may be unrealistic for high-dimensional treatments where $k$ is large relative to $n$. In our work, we first demonstrate that accounting for the dependence of error bounds on $k$ is crucial in this setting, by showing that estimating mean potential outcomes exhibits a slow $k/n$ rate in mean squared error. Ignoring this $k$-dependence obscures actual error rates and can severely overstate accuracy. In our paper, the number of treatments $k$ and the sample size $n$ can be any two numbers; in particular, we require no specific asymptotics, and $k$ need not be proportional to or scale with $n$. We present finite-sample error guarantees for our proposed estimators, which depend on $n$ and $k$ and are valid for any such values, revealing how a large $k$ worsens error guarantees.\\

Finally, we note that our proposed pseudo-outcome regression framework is similar in spirit to the procedure suggested by \cite{foster2019orthogonal}, who present an empirical risk minimization framework for general loss functions that involve nuisances. Their results are general and work for a wide range of loss functions. However, in our high-dimensional treatment setting, it is helpful to exploit the specific treatment structure more explicitly to obtain a doubly robust nuisance dependence and to avoid issues due to violations of strong positivity. For a detailed discussion of this, we refer to Section \ref{sec:two-stage-procedures}.

\section{Setup \& Notation}

We consider an independent and identically distributed (iid) sample $(Z_1,...,Z_n)$ where $Z=(X,A,Y) \sim \Pb$ for covariates $X \in \R^d$, treatment $A \in \mathcal{A}$, and outcome $Y \in \R$. The set $\mathcal{A}$ is a generic set for now and is defined explicitly once we consider specific treatment settings introduced in Section~\ref{section:high-dim-trt}, such as single multi-valued treatments or vector treatments. Let $\mu_a(x)=\E(Y\mid X=x, A=a)$ denote the outcome regression. Let $\pi_a(x)$ denote the treatment probabilities $\Pb(A=a\mid X=x)$ if $A$ is discrete, or the conditional density of $A$ given $X$ if $A$ is continuous. Throughout we assume the positivity condition that $\pi_a(x) > 0$ for all $a\in\mathcal{A}$ (see Section \ref{section:positivity-violation} for discussion). Further, let $\widehat\varpi_a=\frac{1}{n}\sum_{i=1}^n \mathds{1}(A_i=a)$ be the empirical proportions of the treatment, and, for discrete $A$, let $\varpi_a=\Pb(A=a)$ be the marginal treatment probabilities. \\

Our goal is to study estimation of the treatment-specific mean
\begin{equation}
\int \mu_a(x) \ d\Pb(x) \label{eq:parameter}
\end{equation}
across many treatments $a\in\mathcal{A}$, in particular, when the treatment is high-dimensional (which will be made rigorous in Section~\ref{section:high-dim-trt}). Letting $Y^a$ denote the potential outcome under treatment $A=a$, the quantity (\ref{eq:parameter}) also equals the counterfactual mean $\E(Y^a)$, if we make the additional causal assumptions of consistency ($Y=Y^A$) and no unmeasured confounding ($A \ind Y^a \mid X$). However, all our results apply to the statistical quantity defined in \eqref{eq:parameter}, regardless of whether the causal assumptions hold. \\

We denote the Lebesgue measure by $\lambda$. We use  $\| v \|_1 = \sum_j |v_j|$ and $\| v \|_0 = \sum_j \one(v_j \neq 0)$ for the $\ell_1$- and $\ell_0$-norms for vectors $v$. For a function $f$ we let $\| f \|^2_{P,2}= \int f(z)^2 \ dP(z)$ denote the squared $L_2(P)$ norm. We write sample averages with the empirical distribution shorthand $\Pn(f) = \frac{1}{n} \sum_{i=1}^n f(Z_i)$. When we use sample splitting with separate folds $D_s$, we denote the empirical measure over $D_s$ by $\Pb^s_{n}$.

\section{High-Dimensional Treatments}\label{section:high-dim-trt}

In this section, we present various types of high-dimensional treatments considered in this paper. We characterize each of the different treatment settings, highlight their unique challenges, compare them with one another, and finally display the fundamental limits of how well one can possibly estimate mean potential outcomes in a high-dimensional regime.

\subsection{Types of High-Dimensional Treatments}

In this subsection, we distinguish between two structurally different types of high-dimensional treatments. We rely on a dichotomy distinguishing single multi-valued and vector treatments, which was also discussed, for example, in \citet{xiang2025double}.

\paragraph{Single Multi-Valued Treatments} are treatments that can take values from a set of unordered discrete values. More specifically, $A\in\{1,\dots,k\}$, i.e., $A$ can be one of $k$ possible treatment levels. Note that it does not matter how we denote the elements of the set of possible treatment levels since the levels are unordered. For example, the set of possible values can be $k$ different medical treatments denoted by $1,\dots,k$, and if patient $i$ received treatment number $j$, we observe $A_i=j$. \\

A \textit{vector treatment}, also called \textit{multiple treatment} in \citet{xiang2025double}, is a vector of different treatments, i.e., a combination of multiple individual treatments. More specifically, $A=(A_1,\dots,A_k)$ for individual treatments $A_j$. Depending on which values each individual treatment can take, we want to distinguish vector treatments into two categories further: binary and continuous vector treatments.

\paragraph{Binary Vector Treatments}
We call $A$ a \textit{binary vector treatment}, or also \textit{binary multiple treatment}, if each individual treatment of the treatment combination is binary, i.e., $A=(A_1,\dots,A_k)$ and $A_1,\dots,A_k\in\{0,1\}$. In other words, $A\in\{0,1\}^k$. Intuitively, this means that for each unit, we observe a combination that indicates which of the treatments $A_1,\dots,A_k$ was received by the unit. Note that this can be seen as a special case of single multi-valued treatments, which allows for modeling additional structure compared to the single multi-valued treatment case.

\paragraph{Continuous Vector Treatments}
We call $A$ a \textit{continuous vector treatment}, or also \textit{continuous multiple treatment} or just \textit{continuous treatment}, if $A\in\mathbb{R}^k$ is a continuous random variable. Intuitively, this means that for each unit we observe a combination of treatment doses.\\

\begin{remark}
We note a subtlety in our notation. While $A_1,\dots,A_k$ refer to the treatment components of an observed treatment vector $A$, the term $A_i$ can also refer to the $i$-th of the $n$ observations of the treatment variable. The meaning of the subindex will be clear from the context. Moreover, the $j$-th treatment component of observation $i$ is denoted by $A_{i,j}$.
Throughout, we also refer to \textit{discrete} versus \textit{continuous} treatments at times. By \textit{discrete treatments} we refer to either the single multi-valued or binary vector treatment setting, whereas \textit{continuous treatments} refer to continuous vector treatments.\\
\end{remark}

In this paper, for both single multi-valued and vector treatments, our results allow the number of treatments $k$ and sample size $n$ to be any two values, enabling the analysis of arbitrary treatment dimensions. In particular, we require no specific asymptotics, and, although covered by our theory, $k$ does not have to scale with $n$. In this work, by \textit{high-dimensional treatment}, we refer to the setting where the number of treatment levels $k$ can be large relative to $n$ (without relying on any particular asymptotic scaling).

\subsection{Positivity Violations}\label{section:positivity-violation}

In standard settings with $A$ binary, a common assumption for identification and estimation of mean potential outcomes is strong positivity, i.e., $\Pb(A=1\mid X)\geq\varepsilon>0$ with probability one for some $\varepsilon$ that is independent of the sample size $n$. This assumption says that each unit has a non-trivial chance of receiving the treatment. This assumption may be untenable for high-dimensional treatments, particularly those that are discrete, where strong positivity must necessarily be violated for most treatment levels. Therefore, throughout this paper, we allow for so-called \textit{near-violations of positivity}.\\

By \textit{near-violation of positivity}, we refer to the situation when weak positivity $\pi_a(X)>0$ is satisfied (so that the mean potential outcomes can still be identified), but $\pi_a(X)$ cannot be lower bounded by a strictly positive constant that is independent of $n$ and $k$, i.e., strong positivity is violated. Intuitively, this describes the situation where propensity scores are near zero due to a large number of treatments; however, they are never exactly zero. We note that it is possible to extend our results to address even violations of weak positivity, by handling different target parameters, such as incremental effects \citep{kennedy2019nonparametric,schindl2024incremental}.\\

Note here a fundamental difference between discrete and continuous high-dimensional treatments. For discrete high-dimensional treatments, where $k=k(n)$ is an unbounded sequence in $n$, near-violation of positivity is unavoidable. However, near-positivity violations need not occur in continuous treatments. Even with large $k$ (the number of individual treatments of the continuous treatment vector), it may be reasonable to assume that the density $\pi_a(x)$ can be bounded away from zero by a constant independent of $n$ and $k$.\\

Furthermore, there is an important distinction between single multi-valued and vector treatments to be made. For a single multi-valued treatment, even if it is uniformly distributed, we only have roughly $n/k$ observations at each level, whereas for vector treatments, we have an effective sample size of $n$. In particular, for single multi-valued treatments, we inevitably have very few observations at some treatment levels, as illustrated in Figure \ref{fig:treatment-levels-with-few-obs} in the Appendix. This implies that, when estimating mean potential outcomes, we can tolerate smaller sample sizes in comparison to the number of treatments for vector treatments than for single multi-valued treatments (for details refer to Sections \ref{section:single-treatment} and \ref{section:vector-treatment}).\\

\begin{figure}[ht]
    \centering
    \renewcommand{\arraystretch}{1.5}
    \begin{tabular}{ p{0.23\textwidth} |>{\centering\arraybackslash}p{0.22\textwidth} >{\centering\arraybackslash}p{0.22\textwidth} >{\centering\arraybackslash}p{0.22\textwidth} }
        \hline
        \textbf{Treatment type} & \textbf{High-dimensional estimation target} & \textbf{Inevitable near-violation of positivity} & \textbf{Smaller effective sample size} \\
        \hline
        Single multivalued & $\checkmark$ & $\checkmark$ & $\checkmark$ \\
        Binary vector      & $\checkmark$ & $\checkmark$ & $\times$ \\
        Continuous vector  & $\checkmark$ & $\times$     & $\times$ \\
        \hline
    \end{tabular}
    \caption{Comparison of the challenges that each high-dimensional treatment type brings.}
    \label{tab:treatment-challenges-overview}
\end{figure}

In summary, these near-violations of positivity, combined with a high-dimensional target parameter, make the estimation of mean potential outcomes challenging. It is possible to show that without any further structural assumptions, estimation accuracy cannot be better than $k/n$ (in a minimax sense), which is problematic especially when the number of treatments $k$ is large in comparison to the sample size $n$. For details on this and rigorous minimax lower bounds, we refer to Appendix \ref{sec:no-free-lunch}. Moreover, for an overview of the different treatment types and their inherent challenges, we refer to Figure \ref{tab:treatment-challenges-overview}.

\section{Estimation and Minimax Optimality}\label{section:single-treatment}

In this section, we estimate mean potential outcomes for high-dimensional single multi-valued treatments. We propose Lasso and best subset selection estimators, show their double robustness, and provide error guarantees under sparsity. This reveals that faster convergence rates, in particular, rates faster than the $k/n$ minimax rate in Proposition \ref{prop:classiclowerbd}, are achievable with additional structure. Finally, we derive minimax lower bounds in a sparse and structure-agnostic minimax framework, showing that our proposed estimators are optimal and their error rates are unimprovable without adding more assumptions.

\subsection{Setup \& Model Assumption}

Throughout this section, we assume that the treatment variable $A$ is a single multi-valued treatment, i.e., it takes values in the set $\mathcal{A}=\{1,\dots,k\}$, where the number of possible treatment levels $k$ is potentially large. Our goal is to estimate the mean potential outcomes $\E(Y^a), a=1,\dots,k$. While this section focuses on single multivalued treatments and mean potential outcomes as estimation targets, we also propose estimators for other target parameters and high-dimensional vector treatments in Section \ref{section:extensions} by generalizing the results in this section to a general sparse pseudo-outcome regression framework. \\

To achieve faster convergence rates beyond the slow $k/n$-rate when estimating mean potential outcomes (see Proposition \ref{prop:classiclowerbd} in Appendix \ref{sec:no-free-lunch}), we need to assume additional structure. A natural option is to introduce a sparse model for the mean potential outcomes. Sparsity is a very common assumption in high-dimensional regression when the number of features $d$ is large relative to $n$, e.g., as discussed by \citet{wainwright2019highdimensional},  \citet{rigollet2023high}, and many others. For high-dimensional linear regression, without any further assumptions, the minimax risk is $d/n$, where $d$ is the number of features, which is problematic when $d$ is large relative to $n$. Then, to overcome this slow rate, sparsity of the coefficient vector is assumed, and faster rates that depend on $d$ only logarithmically can be derived. The $k/n$ minimax rate of Proposition \ref{prop:classiclowerbd} in the high-dimensional treatment regime can be viewed as an analog to the high-dimensional regression minimax rate $d/n$. Therefore, inspired by high-dimensional regression, we hope for a logarithmic dependence on $k$ once we assume the mean potential outcomes to be sparse. \\

More formally, we assume the sparse model
\begin{equation}\label{eq:saturated-model}
    \E(Y^a) = \int_{\mathbb{R}^d} \mu_a(x)d\Pb(x) = \psi_0 + \psi_a = \psi_0 + \sum_{j=1}^k \psi_j \mathds{1}(a=j)\quad\text{where}\quad\norm{\psi}\leq s
\end{equation}
for some intercept $\psi_0$ that is assumed to be known for simplicity (refer to Remark \ref{remark:intercept} for a discussion when the intercept needs to be estimated) and $\norm{\cdot}$ either the $L_0$- or $L_1$-norm. We note that the linearity in this model is not actually an assumption for multi-valued treatments, since the model is saturated, whereas for vector treatments it rules out interactions (which could be included as separate new treatments, albeit at the expense of increasing dimension). 
Also note that we require causal identifying assumptions for the first equality to hold, but without these assumptions our results still hold for estimating $\int \mu_a(x)d\Pb(x)$. Intuitively, the sparsity in (\ref{eq:saturated-model}) says that the mean potential outcomes do not deviate too much from the intercept $\psi_0$. Our goal is to estimate the functional $\psi:\mathcal{P}\to\mathbb{R}^k$ by leveraging the introduced sparsity to achieve fast convergence rates. In the next subsection, we propose our estimators.

\subsection{Proposed Estimators}

Our proposed estimator of $\psi$ in (\ref{eq:saturated-model}) can be viewed as a two-stage procedure using sample splitting. In the first stage, on the first fold of the data, we estimate the uncentered efficient influence function of the mean potential outcomes (which we also refer to as \textit{pseudo-outcome}), which depends on the unknown regression function and the propensity score. In the second stage of our procedure, we use these estimated pseudo-outcomes in a sparse linear least squares regression on the second fold of the data. For a comparison with other two-stage procedures, such as \cite{foster2019orthogonal}, we refer to the discussion in Section \ref{section:pseudo-outcome-regression} when we generalize the proposed estimators of this section. In the following, we describe both stages of our proposed method.\\

\noindent\textbf{First Stage.} For nuisance parameters $\eta=(\mu,\pi)$, the pseudo-outcome, i.e., the uncentered influence function of the mean potential outcome shifted by the intercept, is given by
\begin{equation*}
    \varphi_a(Z;\eta) = \frac{\mathds{1}(A=a)}{\pi_a(X)}\left\{Y-\mu_a(X)\right\} + \mu_a(X) - \psi_0
\end{equation*}
for all $a\in\mathcal{A}$. Contruct estimators $\widehat\eta=(\widehat\mu, \widehat\pi)$ of the nuisance parameters on the fold $D_1$ of size $n$ that is independent of our sample $D_2=(Z_1,\dots,Z_n)$.
Then, set the estimated pseudo-outcomes
\begin{equation}\label{def:estimated-pseudo-outcome-single-treatment}
\widehat\varphi_{a}(Z)=\varphi_{a}(Z,\widehat\eta)=\varphi_a(Z,\widehat\mu, \widehat\pi)
\end{equation}
to be the plug-in estimate of the uncentered influence function. It is important to highlight two aspects. First, using the efficient influence function as pseudo-outcome comes with the benefit that our estimated pseudo-outcome is close to the unobserved mean potential outcomes up to the estimation error from the nuisance parameters, as $\E\{\varphi_a(Z;\eta)\}=\psi_a$. Second, the error from nuisance estimation will turn out to be doubly robust due to the definition of the efficient influence function via the von Mises expansion, yielding a second-order remainder term. More specifically, we show that the nuisance penalty has a second-order dependence on
\begin{equation*}
    \delta_n(a) = \anorm{\frac{\pi_a}{\widehat\pi_a} - 1}_{\Pb_X,2} \quad\text{and}\quad \epsilon_n(a) = \anorm{\mu_a - \widehat\mu_a}_{\Pb_X,2}.
\end{equation*}\\

\noindent\textbf{Second Stage.} On the second fold $D_2$, we regress those estimated pseudo-outcomes obtained in the first stage on treatment indicators. Define the Lasso estimator
\begin{equation}\label{def:lasso-single-treatment}
    \widehat\psi_{\text{lasso}} = \argmin_{\norm{\beta}_1\leq s} \widehat R(\beta)
\end{equation}
and the best subset selection estimator
\begin{equation}\label{def:bestsubset-single-treatment}
    \widehat\psi_{\text{subset}} = \argmin_{\norm{\beta}_0\leq s} \widehat R(\beta)
\end{equation}
as the minimizer of the empirical risk $\widehat R$ subject to an $L_1$- and $L_0$-constraint, respectively, where the empirical risk is given by
\begin{equation}\label{eq:emp-risk-single-trt}
    \widehat R(\beta) = \Pb^2_{n}\left\{-2 \sum_{a=1}^k\left[ \frac{\mathds{1}(A=a)}{\widehat\pi_a(X)}\{Y-\widehat\mu_a(X)\}  + \widehat\mu_a(X) - \psi_0 \right]\beta_a \widehat\varpi_a + \beta_A^2\right\}
\end{equation}
where $\widehat\varpi_a$ is the empirical proportion on $D_1$. By imposing an $L_0$- or $L_1$-constraint on the coefficient vector in this regression, we aim to leverage sparsity of the vector $\psi$ and achieve fast convergence rates even when $k$ is large in proportion to $n$. Note that we could switch to a penalized (instead of constrained) regression in the second stage while maintaining the same theoretical guarantees. For a discussion of constrained versus penalized Lasso and best subset selection approaches, we refer to \citet{wainwright2019highdimensional} or \citet{hastie2020best}.\\

\begin{remark}[Motivation for this risk function]\label{remark:motivation-risk-function-single-trt}
    With access to the mean potential outcomes, our goal would be to find
    $$
        \psi = \argmin_{\norm{\beta}\leq s} \sum_{a=1}^k \varpi_a\{\E(Y^a) -\psi_0 - \beta_a\}^2 = \argmin_{\norm{\beta}\leq s} \sum_{a=1}^k \varpi_a\left[ -2\{\E(Y^a)-\psi_0\}\beta_a + \beta_a^2 \right]
    $$
    where $\varpi_a=\Pb(A=a)$ is the marginal treatment probability. The uncentered influence function of the oracle risk $R(\beta)=\sum_{a=1}^k \varpi_a\left[ -2\{\E(Y^a)-\psi_0\}\beta_a + \beta_a^2 \right]$  equals
    $$
        -2\sum_{a=1}^k\E(Y^a)\beta_a\mathds{1}(A=a) + 2\sum_{a=1}^k \varpi_a\beta_a\E(Y^a) - 2\sum_{a=1}^k \varpi_a\beta_a \varphi_a(Z) + \beta_A^2,
    $$
    where $\varphi_a$ is the uncentered influence function of $\E(Y^a)$, assuming $\psi_0=0$ for simplicity. We now use parts of the efficient influence function to obtain the expression
    $$
        - 2\sum_{a=1}^k \varpi_a\beta_a \varphi_a(Z) + \beta_A^2 = -2 \sum_{a=1}^k\left[ \frac{\mathds{1}(A=a)}{\pi_a(X)}\{Y-\mu_a(X)\}  + \mu_a(X) - \psi_0 \right]\beta_a \varpi_a + \beta_A^2,
    $$
    which equals the oracle risk $R(\beta)$ in expectation. We can now plug in an estimate $\widehat\pi$ of the propensity score and $\widehat\mu$ of the regression function, $\widehat\varpi_a$ for $\varpi_a$, and estimate the mean with the sample average to obtain the estimated risk $\widehat R(\beta)$ in \eqref{eq:emp-risk-single-trt}. This motivates the choice of our empirical risk for the Lasso and best subset selection estimator.
    \\
\end{remark}

\begin{remark}[Connection to Gaussian sequence model]
    Since the Lasso and best subset selection estimator are obtained through a constrained regression on indicators via the above risk function, they are heuristically very similar to soft- and hard-thresholding estimators in the Gaussian sequence model, where $\overline Y^a\sim N\left(\psi_0 + \psi_a, \frac{1}{N_a}\right)$ with $N_a\sim\mathrm{Bin}(n, \Pb(A=a))$. However, this analogy is limited: First, the sample size $N_a$ is random; specifically, we only have around $n/k$ observations at each level instead of $n$. Second, we do not observe $\overline Y^a$ directly; instead, it has to be estimated using the covariates. Lastly, we do not assume Gaussianity of the outcome, nor are our (estimated) pseudo-outcomes even sub-Gaussian. In the Appendix \ref{sec:gaussian-sequence-model}, we make this connection rigorous.\\
\end{remark}

\begin{remark}[Alternative risk function]
    As an alternative, one might wish to weight the treatments according to some fixed weight function (i.e., that is not chosen data-dependent) in the above risk function. While this alternative risk also appears natural, we point out some important caveats of using fixed weights in Appendix \ref{sec:alternative-risk}, to which we refer for details.\\
\end{remark}

\begin{remark}[Generalization]
    We note that this two-stage procedure can be generalized to generic high-dimensional target parameters and extended to vector treatments. In Section~\ref{section:pseudo-outcome-regression}, we present a general sparse pseudo-outcome regression framework, which can be used to construct estimators for arbitrary high-dimensional parameters and treatment types. The risk function and estimators proposed above are special cases, as discussed in Remark \ref{remark:recovering-single-multivalued-treatment-setting}.
\end{remark}

\subsection{Error Guarantee}

In this section, we present the error rate of our proposed Lasso and best subset selection estimator under exact sparsity, i.e., when $\norm{\psi}_0\leq s$. We demonstrate that a fast Lasso rate can be achieved and that fast convergence rates are possible, even when the number of treatments is large relative to the sample size. For an analysis under approximate sparsity, i.e., when $\norm{\psi}_1\leq s$, we refer to Appendix \ref{section:single-treatment-approximate-sparsity}.

\begin{theorem}\label{thm:single-treatment}
    Let $\widehat\psi\in\{\widehat\psi_\text{lasso}, \widehat\psi_\text{subset}\}$ be either the Lasso estimator defined in (\ref{def:lasso-single-treatment}) or the best subset selection estimator defined in (\ref{def:bestsubset-single-treatment}) minimizing the empirical risk in (\ref{eq:emp-risk-single-trt}). Assume the model in (\ref{eq:saturated-model}) with exact sparsity, i.e., $\norm{\psi}_0\leq s$ for some sparsity constraint $s$. Moreover, assume that for all $a\in\{1,\dots,k\}$
    \begin{itemize}
        \item[(i)] (Nearly uniform treatment probabilities) $\frac{c}{k}\leq\varpi_a=\Pb(A=a)\leq \frac{C}{k}$,
        \item[(ii)] (Estimated propensities close to empirical weights) $|\widehat\varpi_a/\widehat\pi_a(X)| \leq B$ almost surely, and
        \item[(iii)] (Boundedness)
            $\left|\widehat\mu_a(X)\right|, |Y|, |\psi_0| \leq B$ almost surely
    \end{itemize}
    for constants $c,C,B>0$ that are independent of $k$ and $n$. Then,
    \begin{equation*}
        \E\left(\sum_{a=1}^k \varpi_a(\widehat\psi_a - \psi_a)^2\right) \lesssim s\frac{\log k}{n} + \frac{s}{k}(\delta_n\epsilon_n)^2
    \end{equation*}
    for $\varpi_a=\Pb(A=a)$, $k\geq 2$, $n \geq \max\{\gamma k, s^4/k^2\}\log k$ for a large enough constant $\gamma$, when assuming that $\delta(a)\leq \delta_n$ and $\epsilon_n(a)\leq \epsilon_n$ for all $a\in\{1,\dots,k\}$.
\end{theorem}

The error guarantee given in the theorem consists of two parts: the oracle rate $s\frac{\log k}{n}$ and the nuisance estimation error $s(\delta_n\epsilon_n)^2$. The oracle error rate can be viewed as the estimation error that we would suffer even if we used the unobserved potential outcomes $Y^a$ in the regression. The nuisance estimation error arises from the need to estimate both the regression function and the propensity score, i.e., the price paid for not observing the potential outcomes.\\

Notably, the oracle rate $s\frac{\log k}{n}$ reveals that we can significantly improve on the slow $k/n$ convergence rate stated in Proposition \ref{prop:classiclowerbd} when assuming sparsity of the vector $\psi$. Note that our error now only depends on $k$ logarithmically instead of linearly. This shows that consistent estimation is possible even in the high-dimensional treatment regime.\\

The nuisance estimation error term contains a second-order product $(\delta_n\epsilon_n)^2$ of the regression function and propensity score estimation error. For this reason, we refer to our estimator as \textit{doubly robust}. When assuming parametric models for the nuisance estimation, it is enough to correctly specify one of the two nuisance models for consistent estimation. Even when flexible nonparametric machine learning tools are used to estimate the nuisance parameters, is it sufficient for the product of the error rates to be small enough. In particular, we achieve the oracle rate $s\frac{\log k}{n}$ whenever $\delta_n\epsilon_n\lesssim \sqrt{\frac{k\log k}{n}}$. Hence, the individual nuisance error rates can be slower as long as the product of the two satisfies the rate requirement. For instance, this is the case when $\delta_n,\epsilon_n\asymp \left(\frac{k\log k}{n}\right)^{1/4}$. Note that the requirement $\delta_n\epsilon_n\lesssim \sqrt{\frac{k\log k}{n}}$ to achieve the oracle error rate corresponds to $\delta_n\epsilon_n\lesssim 1/\sqrt{n}$ in the classical doubly robust estimation setting up to a $\log k$ inflation factor. The difference arises from the fact that our effective sample size is $n/k$ instead of $n$. In the following remarks, we comment on the assumptions of Theorem \ref{thm:single-treatment} and provide an alternative upper bound for the error that is slightly more informative.\\

\begin{remark}[Exact versus approximate sparsity]\label{remark:exact-versus-approximate-sparsity}
    The above theorem assumes exact sparsity of $\psi$. Depending on the application, this might appear unreasonable. Instead, one might resort to approximate sparsity, where we assume that $\norm{\psi}_1\leq s$. Even under approximate sparsity, we are able to derive an error guarantee of the above flavor; however, with some distinctions. First, in this setting, we can only provide an error guarantee for the Lasso estimator $\widehat\psi_\text{lasso}$ based on our analysis technique, not for the best subset selection estimator. Second, by assuming only approximate sparsity, we have to pay the price of obtaining a slower convergence rate $s\sqrt{\frac{\log k}{kn}}$ instead of the fast rate $s\frac{\log k}{n}$. Third, obtaining a fast rate under exact sparsity typically requires a restricted eigenvalue assumption on the design matrix of the regression. In the single multivalued treatment case, this assumption is satisfied by design. For details on the results under approximate sparsity, we refer to Appendix \ref{section:single-treatment-approximate-sparsity}.\\
\end{remark}

\begin{remark}[Sparsity constraint]\label{remark:sparsity-constraint}
    The theorem uses the fact that the constraint tuning parameter $s$ in the Lasso estimator (\ref{def:lasso-single-treatment}) and the best subset selection estimator (\ref{def:bestsubset-single-treatment}) is chosen as the true sparsity, that is, $s=\norm{\psi}_0$ for best subset selection and $s=\norm{\psi}_1$ for the Lasso. However, the sparsity of the vector $\psi$ is often unknown, which makes this choice infeasible in practice. Instead, the constraint tuning parameter in the Lasso and best subset selection estimator can be chosen by cross-validation. We note without further details that it is possible to state a cross-validation procedure based on minimizing an estimated risk and provide theoretical guarantees in the form of oracle inequalities, demonstrating the validity of the procedure.\\
\end{remark}

\begin{remark}[Assumptions on treatment probabilities]
    Assuming that the treatment probabilities $\varpi_a$ are nearly uniform implies that each treatment has a roughly equal chance of being applied. It essentially rules out the possibility that some treatments are too rare. We further assume that the estimated propensity $\widehat\pi_a(X)$ is close to the empirical weights $\widehat\varpi_a$ on $D_1$. This is a mild assumption since $\widehat\pi_a(X)$ is also constructed on $D_1$ and can therefore be thresholded from below to satisfy Assumption~(ii) (potentially at the expense of a larger nuisance penalty $\delta_n(a)$).\\
\end{remark}


\begin{remark}[Ultra-high-dimensional regime]
    Our theorem states the rate for $k\log k\lesssim n$. When $k\gtrsim n$, the same rate can be derived, but then  $(\delta_n\epsilon_n)^2\lesssim 1\lesssim \frac{k\log k}{n}$ and the oracle rate dominates. In that regime, the oracle rate is not optimal and can be beaten by a trivial estimator of zero. This trivial estimator has a rate of order $s/k$, which  turns out to be minimax optimal  in that regime. Thus, meaningful estimation is not possible when $k\gtrsim n$, even under sparsity. This phenomenon also occurs in standard high-dimensional regression when $\log d>N$, i.e., the ultra-high-dimensional setting (\citet{verzelen2012minimax}), for $d$  dimension and $N$ sample size. The case $k\gtrsim n$ is analogous to the ultra-high-dimensional regime for our problem, noting that our effective sample size is $n/k$ instead of $n$, so $\log d > N$ reduces to $k\log k > n$ for $N=n/k$ and $d=k$.\\
\end{remark}

\begin{remark}[Intercept]\label{remark:intercept}
    In (\ref{eq:saturated-model}), we assumed the intercept $\psi_0$ to be known. For many practical applications, this may not hold, and the intercept must be estimated. However, for many choices of the intercept, this estimation task is much easier in comparison to estimating the high-dimensional vector $\psi$, and the error rate corresponding to the estimation of $\psi_0$ is negligible. For instance, a reasonable choice of the intercept could be the average of all mean potential outcomes, i.e., $\psi_0=\frac{1}{k}\sum_{a=1}^k \E(Y^a)$. An estimator $\widehat\psi_0$ could be the average of the estimated uncentered efficient influence functions of the mean potential outcomes. The estimated pseudo-outcomes in (\ref{def:estimated-pseudo-outcome-single-treatment}) can then be defined using the estimate $\widehat\psi_0$ instead of the true $\psi_0$. In the proof of Theorem \ref{thm:single-treatment}, we then additionally need to consider the estimation error arising from $\widehat\psi_0$, which, however, turns out to yield a negligible doubly robust rate with a straightforward calculation. Consequently, the same error guarantee in Theorem \ref{thm:single-treatment} can be achieved when the intercept, chosen as the average mean potential outcome, has to be estimated. A similar behavior is expected for other choices of the intercept, such as the median, the mean potential outcome at one specific level, or general linear combinations of the target parameter.\\
\end{remark}

\begin{remark}[Alternative Upper Bound]\label{remark:weighted-nuisance-error}
    Instead of the upper bound stated in the theorem, it is possible to derive the following upper bound with a slightly more informative nuisance error term:
    \begin{equation*}
        \E\left(\sum_{a=1}^k \varpi_a(\widehat\psi_a - \psi_a)^2\right) \lesssim s\frac{\log k}{n} + sk\E\left[ \max_{a\in\{1,\dots,k\}} \left\{\widehat\varpi_a \delta_n(a)\epsilon_n(a) \right\}^2 \right].
    \end{equation*}
    This result is shown along the way in the proof of Theorem \ref{thm:single-treatment}. This upper bound is more informative since it states the nuisance error in terms of the estimation error of the regression function and propensity score at each specific treatment level $a\in\{1,\dots,k\}$. Hence, we do not rely on upper bounding $\delta_n(a),\epsilon_n(a)$ by $\delta_n,\epsilon_n$ across all treatment levels. Most importantly, this nuisance error term reveils that $\delta_n(a)$ and $\epsilon_n(a)$ are weighted by the empirical proportions $\widehat\varpi_a$ of the treatment variable on $D_1$. Hence, if we only have very few observations at a treatment level available for nuisance estimation, then the term allows for greater nuisance estimation errors at this level. In particular, if a treatment level is entirely unobserved, then $\widehat\varpi_a=0$ and nuisance estimation at this level is allowed to be arbitrarily inaccurate without increasing the overall error.\\
\end{remark}

\begin{remark}[Propensity score estimation]
    We note that propensity score estimation can be increasingly challenging for a larger number of treatments. More specifically, suppose we have access to an estimator $\widehat\pi_a$ which estimates $\pi_a$ with error $R$ in the sense that $\norm{\widehat\pi_a - \pi_a}=R$. Now assume that $\widehat\pi_a\approx\frac{1}{k}$. Then,
    $$
        \anorm{\frac{\pi_a}{\widehat\pi_a} - 1} = \anorm{\frac{\pi_a - \widehat\pi_a}{\widehat\pi_a}} \approx kR,
    $$
    i.e., the additive error is multiplied by the number of treatments. This suggests that targeting the inverse propensity score itself (e.g., as in \cite{vanderlaan2026automatic}) could be beneficial in high-dimensional regimes, which specifically controls the error measured in terms of $\anorm{\frac{\pi_a}{\widehat\pi_a} - 1}$, i.e., the quotient of $\pi_a$ and $\widehat\pi_a$. We leave the investigation of robust methods for estimating propensity scores in the high-dimensional treatment regime to future work.
\end{remark}

\subsection{Minimax Lower Bounds}

In this section, we derive minimax lower bounds for estimating mean potential outcomes for single multi-valued treatments under sparsity. More specifically, we state an oracle minimax lower bound as well as a minimax lower bound for the nuisance estimation error in the structure-agnostic minimax framework, both under exact sparsity. Finally, we claim optimality of our proposed Lasso and best subset selection estimator.

\subsubsection{Motivation \& Setup}

Minimax rates are generally defined as $R_n=\inf_{\widehat\psi} \sup_{P\in\mathcal{P}} \E_P(\ell(\widehat\psi-\psi))$ for some statistical model $\mathcal{P}$ and loss function $\ell$. Intuitively, the minimax rate describes the best possible (worst-case) error of any estimator.\\

Studying minimax rates is of crucial interest for several reasons: First, it helps answer the question of whether a particular estimator can be improved or is already optimal. If a derived minimax lower bound and the error rate of a particular estimator match, then this estimator is optimal; i.e., every other estimator would perform worse for at least one distribution in the model. In case a lower bound on the minimax risk is smaller than the error rate of a particular estimator, two scenarios are possible: Either we can improve the minimax lower bound and the estimator is optimal, or the lower bound is tight, and a better estimator can be found. Secondly, minimax rates reveal the fundamental limits of an estimation problem and determine what estimation error one can hope for. Finally, minimax rates allow for a comparison between different estimation problems in terms of their difficulty.\\

It is worth noting that when the goal is to prove the optimality of a given estimator, it is appropriate to make stronger assumptions in the model for the minimax result than in the model used to derive the upper bound on the error of the given estimator. Intuitively, this can be explained by the fact that additional assumptions only make estimation easier, so every minimax lower bound for a smaller model $\mathcal{P'}$ is also a minimax lower bound for every larger model $\mathcal{P''}\supseteq\mathcal{P}'$.\\

In this paper, our primary motivation for studying minimax rates is to verify the optimality of our proposed Lasso and best subset selection estimators for estimating mean potential outcomes in the exactly-sparse single multi-valued treatment regime.

\subsubsection{Minimax Rate under Exact Sparsity}

In the following, we state the minimax rate for estimating mean potential outcomes for high-dimensional single multi-valued treatments under exact sparsity. Subsequently, we conclude optimality of the oracle rate achieved in Theorem \ref{thm:single-treatment}.

\begin{theorem}\label{thm:minimax-lb-sparse-regime}
    Let $s\geq 9$ and $k\geq 8s$. Further, let $\mathcal{P}$ denote the set of all distributions for which:
    \begin{itemize}
        \item[(i)] $\frac{C'}{k}\leq \pi_a(X)\leq \frac{C''}{k}$,
        \item[(ii)] $Y$ is binary,
        \item[(iii)] $(A,Y)\ind X$, and
        \item[(iv)] the vector $\psi=(\psi_1,\dots,\psi_k)$ of the functionals $\psi_a=\int \mu_a(x)d\Pb(x)$ is sparse in the sense that at most $s$ entries are different from $1/2$.
    \end{itemize}
    Then,
    \begin{equation*}
        \inf_{\widehat\psi} \sup_{P\in\mathcal{P}} \E\left\{\sum_{a=1}^k \varpi_a(P)(\widehat\psi_a - \psi_a(P))^2\right\} \geq C\cdot \min\left(s\frac{\log(k/s)}{n}, \frac{s}{k}\right)
    \end{equation*}
    for $\varpi_a(P)=P(A=a)$ and $C$ a constant depending only on $C'$ and $C''$.
\end{theorem}

Most importantly, note that the model assumptions made in the above theorem are strictly stronger than the ones used in the upper bound results for our Lasso and best subset selection estimator in Theorem \ref{thm:single-treatment}: Our minimax model assumes nearly uniform propensity scores in (i), which implies nearly uniform treatment probabilities. Moreover, Assumption (ii) implies the required boundedness of the outcome in Theorem \ref{thm:single-treatment}. Assumption (iii) further shrinks the model. Assumption (iv) implies that the deviations from the intercept are exactly sparse, i.e., $\norm{\psi-\psi_0}_0\leq s$ for $\psi_0=1/2$.\\

Given that the above model in the minimax result is contained in the model in the upper bound result of Theorem \ref{thm:single-treatment}, we can conclude optimality of the oracle rate achieved by our estimators $\widehat\psi\in\{\widehat\psi_{\mathrm{lasso}}, \widehat\psi_{\mathrm{subset}}\}$. It is worth mentioning that there is a slight mismatch between the terms $\log (k/s)$ in the minimax rate and $\log k$ in the upper bound. However, in the regime where $k/s\asymp k^{\gamma}$ for some $\gamma>0$, the rate $\log(k/s)$ is equivalent to $\log k$ (up to constants), in which case we can also write the minimax lower bound as
\begin{equation*}
    \inf_{\widehat\psi} \sup_{P\in\mathcal{P}} \E_P\left\{\sum_{a=1}^k \varpi_a(P)(\widehat\psi_a - \psi_a(P))^2\right\} \gtrsim \min\left(s\frac{\log k}{n}, \frac{s}{k}\right),
\end{equation*}
which now exactly matches the upper bound result in Theorem \ref{thm:single-treatment}.\\

The main idea behind the proof of the theorem is to construct distributions that are close, but for which the target parameter $\psi$ is separated as much as possible. Then, intuitively, no estimator can perform better than this separation, since the distributions are statistically indistinguishable. The classical Le Cam method constructs two such distributions, which suffice to derive minimax lower bounds for nonparametric regression at a point (\cite{tsybakov2009introduction}). However, for our purposes, this method would give a lower bound that is not tight enough, so we use multiple hypotheses. More specifically, we employ Fano's method over a pruned and sparsified hypercube, as described in \cite{tsybakov2009introduction}.\\

The derived minimax risk can be viewed as an analogue of the minimax risk for high-dimensional linear regression, which is of the order $s\log d/n$ (e.g., refer to \cite{raskutti2011minimax}), where our derived rate additionally reflects the effective sample size of $n/k$ in the high-dimensional single multi-valued treatment regime (since the error is weighted by the treatment probabilities $\varpi_a\approx1/k$).\\

The above results prove the optimality of the oracle rate of our proposed estimators, and therefore the optimality of the estimators if they achieve this oracle rate, i.e., the nuisance estimation error is of a smaller or equal order. However, it does not reveal whether the nuisance penalty in Theorem \ref{thm:single-treatment} is minimax optimal, which is covered in the following subsection.

\subsubsection{Minimax Nuisance Error Dependence under Exact Sparsity}

In the following, we show the optimality of the nuisance error dependence of our proposed Lasso and best subset selection estimator. This is done by stating the minimax rate in the structure-agnostic minimax framework.\\

In the classical minimax framework, one typically assumes the relevant distribution is smooth or otherwise structured to derive minimax optimal rates. For instance, standard one-step estimators are known to be sub-optimal in many smoothness classes \citep{bickel1988estimating, birge1995estimation}. However, this approach does not allow for a structure-agnostic estimation of the nuisance components, as optimal estimators carefully leverage the assumed structure of the model.\\

In the \textit{structure-agnostic} minimax framework, introduced by \cite{balakrishnan2023fundamental}, we take nuisance estimators as given (in our case $\widehat\pi$ and $\widehat\mu$), and they are treated as black-box estimators. More specifically, these estimators can be obtained using any estimation procedure applied to a separate fold of the data, which is independent of the data used to run the regression. We further assume that those estimators achieve given error rates $\delta_n$ and $\epsilon_n$, i.e.,
\begin{equation*}
    \max_{a\in\{1,\dots,k\}} \anorm{\frac{\pi_a}{\widehat\pi_a}-1}\lesssim \delta_n \quad \text{and} \quad \max_{a\in\{1,\dots,k\}} \anorm{\widehat\mu_a - \mu_a}\lesssim \epsilon_n.
\end{equation*}
Importantly, these error rates are treated as unknown.\\

When taking the structure-agnostic viewpoint, we aim to determine the optimal error for estimating our statistical functional $\psi$ when we have access to nuisance parameter estimates that achieve the above error guarantees. Therefore, we define the model of all distributions such that this error guarantee of the nuisance parameter estimators is satisfied, and derive the minimax risk under this model. Intuitively, we investigate optimality locally around the given nuisance parameter estimators. The following theorem states the minimax risk for this structure-agnostic model.

\begin{theorem}\label{thm:minimax-lb-structure-agnostic}
    Suppose that we have two pilot estimators $\widehat\pi$ and $\widehat\mu$ that achieve the error rates
    \begin{equation}\label{error-rate}
        \max_{a\in\{1,\dots,k\}} \anorm{\frac{\pi_a}{\widehat\pi_a}-1}\leq C_\pi\delta_n \quad \text{and} \quad \max_{a\in\{1,\dots,k\}} \anorm{\widehat\mu_a - \mu_a}\leq C_\mu\epsilon_n
    \end{equation}
    for some $\delta_n, \epsilon_n=o(1)$ and constants $C_\pi, C_\mu>0$ such that
    \begin{equation*}
        \varepsilon\leq\widehat\mu_a\leq 1-\varepsilon,\quad\frac{C_1}{k}\leq\widehat\pi_{a'}(X)\leq\frac{C_2}{k},\quad\text{and}\quad \widehat\pi_k(X)\geq \varepsilon
    \end{equation*}
    almost surely for all $a\in\{1,\dots,k\}, a'\in\{1,\dots,k-1\}$ for some $\varepsilon\in(0,1/2), C_1,C_2>0$ independent of $n$ and $k$, and $\widehat\pi_k(x)=1-\sum_{j=1}^{k-1} \widehat\pi_j(x)$.
    Let $\mathcal{P}=\mathcal{P}(\delta_n, \epsilon_n)$ denote the model where $Y\in\{0,1\}$ is binary, $X$ is uniformly distributed on $[0,1]^d$,
    \begin{equation*}
        \frac{C_1'}{k}\leq\pi_{a'}(X)\leq \frac{C_2'}{k},\quad \pi_k(X)\geq \varepsilon'
    \end{equation*}
    almost surely for all $a'\in\{1,\dots,k-1\}$ for constants $C_1',C_2',\varepsilon'>0$ independent of $n$ and $k$, the pilot estimators achieve the error rate in (\ref{error-rate}), and the vector $\psi$ with entries $\psi_a=\int\mu_a(x)d\Pb(x), a=1,\dots,k$ is $s$-sparse, i.e., only $s$ entries are different from the value $\theta$, where $\psi_k=\theta$.
    Then, the minimax rate is lower bounded as
    \begin{equation*}
        \inf_{\widehat\psi} \sup_{P\in\mathcal{P}} \E_P\left\{\sum_{a=1}^k \varpi_a(P)\left(\widehat\psi_a - \psi_a(P)\right)^2\right\} \gtrsim \frac{s}{k}(\delta_n\epsilon_n)^2
    \end{equation*}
    where $\varpi_a(P)=P(A=a)$.
\end{theorem}

This theorem shows that, among all distributions for which the pilot estimators achieve the rates $\delta_n$ and $\epsilon_n$, no estimator can do better than $\frac{s}{k}(\delta_n\epsilon_n)^2$. Consequently, the nuisance error dependence of our Lasso and best subset selection estimator in Theorem \ref{thm:single-treatment} is optimal and cannot be improved. In combination with the minimax risk in the sparse regime shown in Theorem \ref{thm:minimax-lb-sparse-regime}, this allows us to claim optimality of the Lasso and best subset selection estimator.\\

Note that the assumptions in the above theorem are stronger than the ones used in the upper bound result in Theorem \ref{thm:single-treatment}. Hence, the model in the above minimax result is strictly smaller, allowing us to conclude optimality of our estimators. The near uniform treatment probabilities assumption of Theorem \ref{thm:single-treatment} is satisfied when the propensity scores are nearly uniform, where assuming $\pi_k(X)\geq\epsilon$ on the last treatment level can be viewed as a stronger assumption since, intuitively, it should make estimation for treatment level $k$ easier. Since $Y$ is chosen to be binary, the outcome boundedness assumption is also satisfied. Moreover, the above theorem assumes a slightly stronger version of exact sparsity, as it additionally requires $\psi_k=\theta$. Finally, we can assume that the data used to construct the pilot estimators $\widehat\pi$ and $\widehat\mu$ satisfies $\Pb_n^1\{\mathds{1}(A=a)\}\asymp 1/k$ for all treatments $a\in\{1,\dots,k\}$ (since the above results treats $D_1$ as given and fixed); then, the assumption of nearly uniform estimated propensity scores implies Assumption (ii) of Theorem \ref{thm:single-treatment}.\\

The overall strategy to prove the theorem is to construct distributions that are close, but such that the target parameter $\psi$ is maximally separated under these distributions. 
Our construction is similar to that of \cite{jin2025structureagnostic} but repeated across many treatment levels and ensuring the assumed sparsity is respected. These distributions are constructed using the method of fuzzy hypothesis \citep{birge1995estimation, ibragimov1987some, ingster2003nonparametric, nemirovski2000topics, robins2009semiparametric, tsybakov2009introduction, jin2025structureagnostic, kennedy2024minimax}. More specifically, this entails constructing two sets of distributions, along with a prior. This prior can be viewed as a mixture, providing us with a pair of such mixtures. Then, we show that these mixtures are close in Hellinger distance, while the functional is separated by the distance of the desired lower bound. This approach is formalized by Lemma~\ref{lem:minimax} (which is adapted from Theorem 2.15 in \cite{tsybakov2009introduction}).\\

\begin{figure}[ht!]
    \centering
    \begin{tabular}{p{0.15\textwidth} cc}
        
        & \textbf{Null distribution $P_\lambda$} & \textbf{Alternative distribution $Q_\lambda$} \\[1em]
        
        \textbf{Active \mbox{levels} $a\in S$} & 
        \raisebox{-0.5\height}{\includegraphics[width=0.35\textwidth]{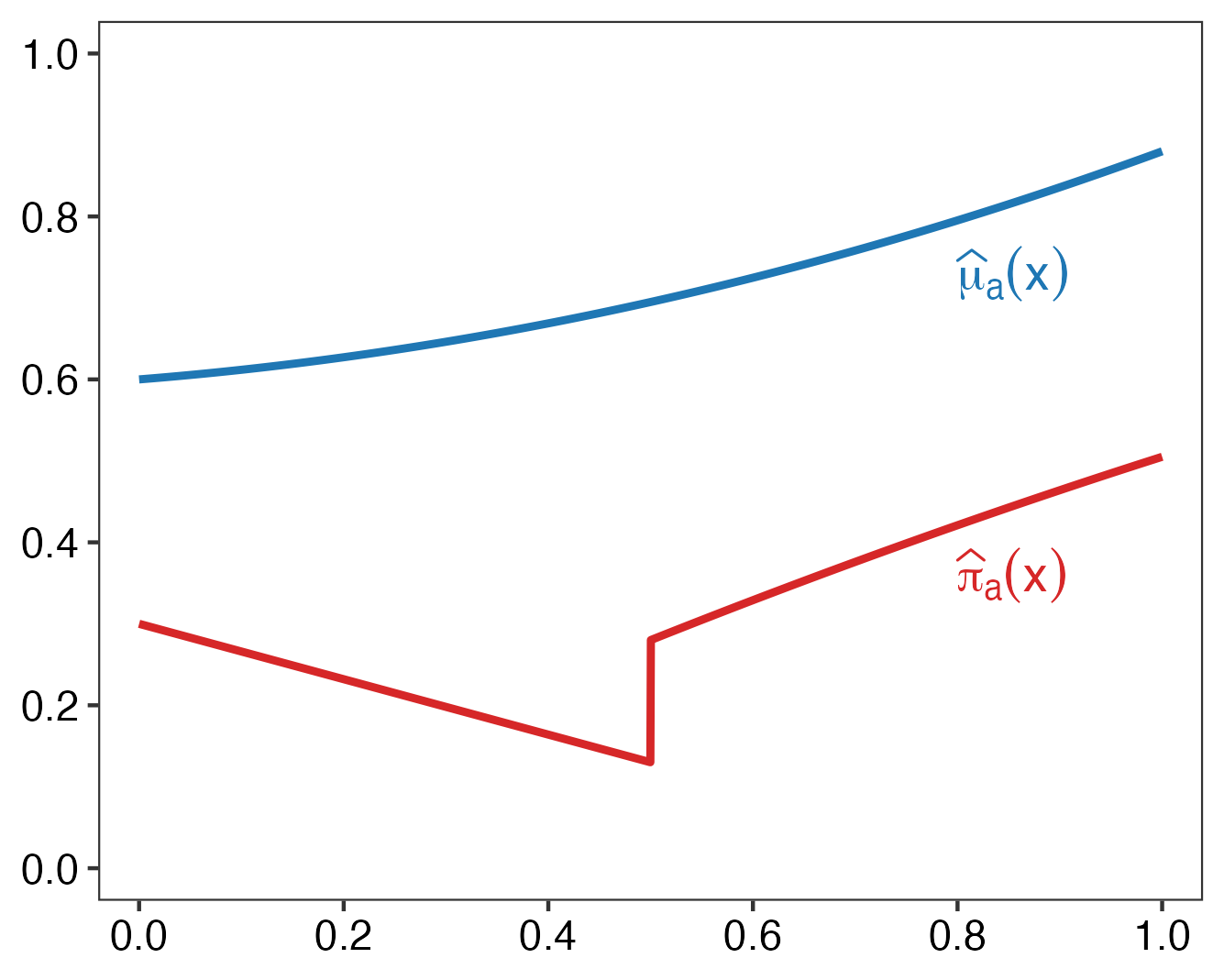}} & 
        \raisebox{-0.5\height}{\includegraphics[width=0.35\textwidth]{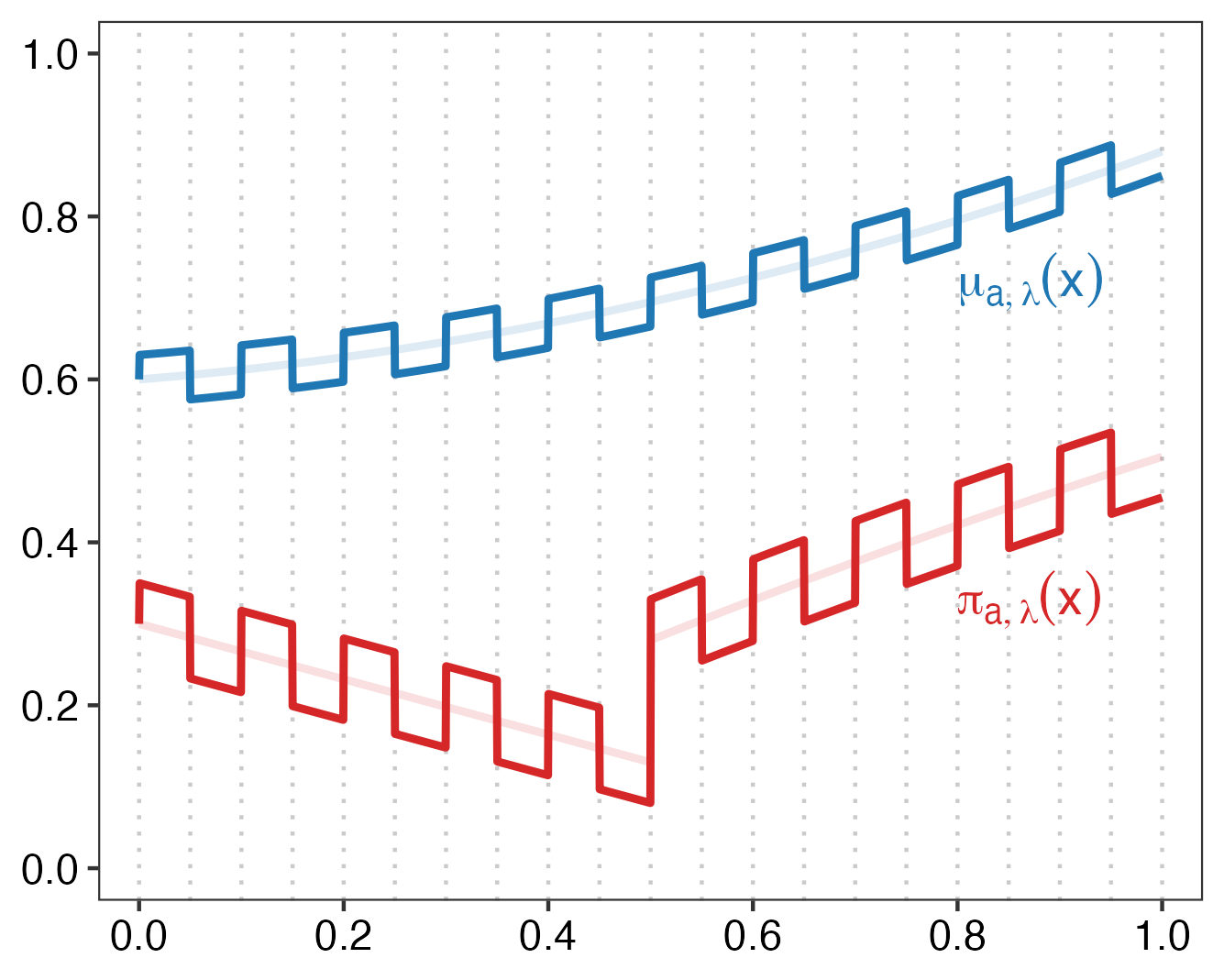}} \\[2em]

        \textbf{Non-active levels $a\notin S, a\neq k$} & 
        \raisebox{-0.5\height}{\includegraphics[width=0.35\textwidth]{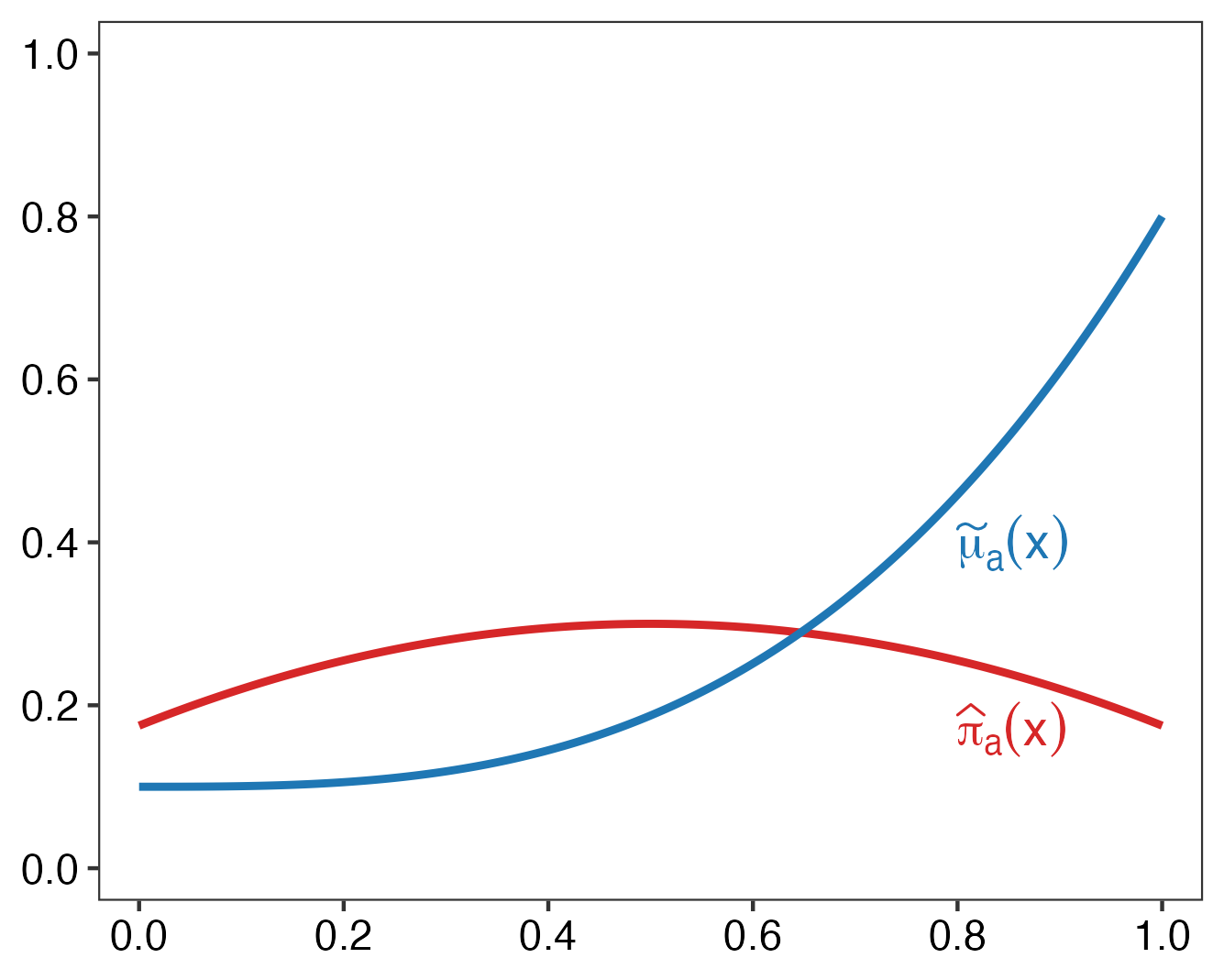}} & 
        \raisebox{-0.5\height}{\includegraphics[width=0.35\textwidth]{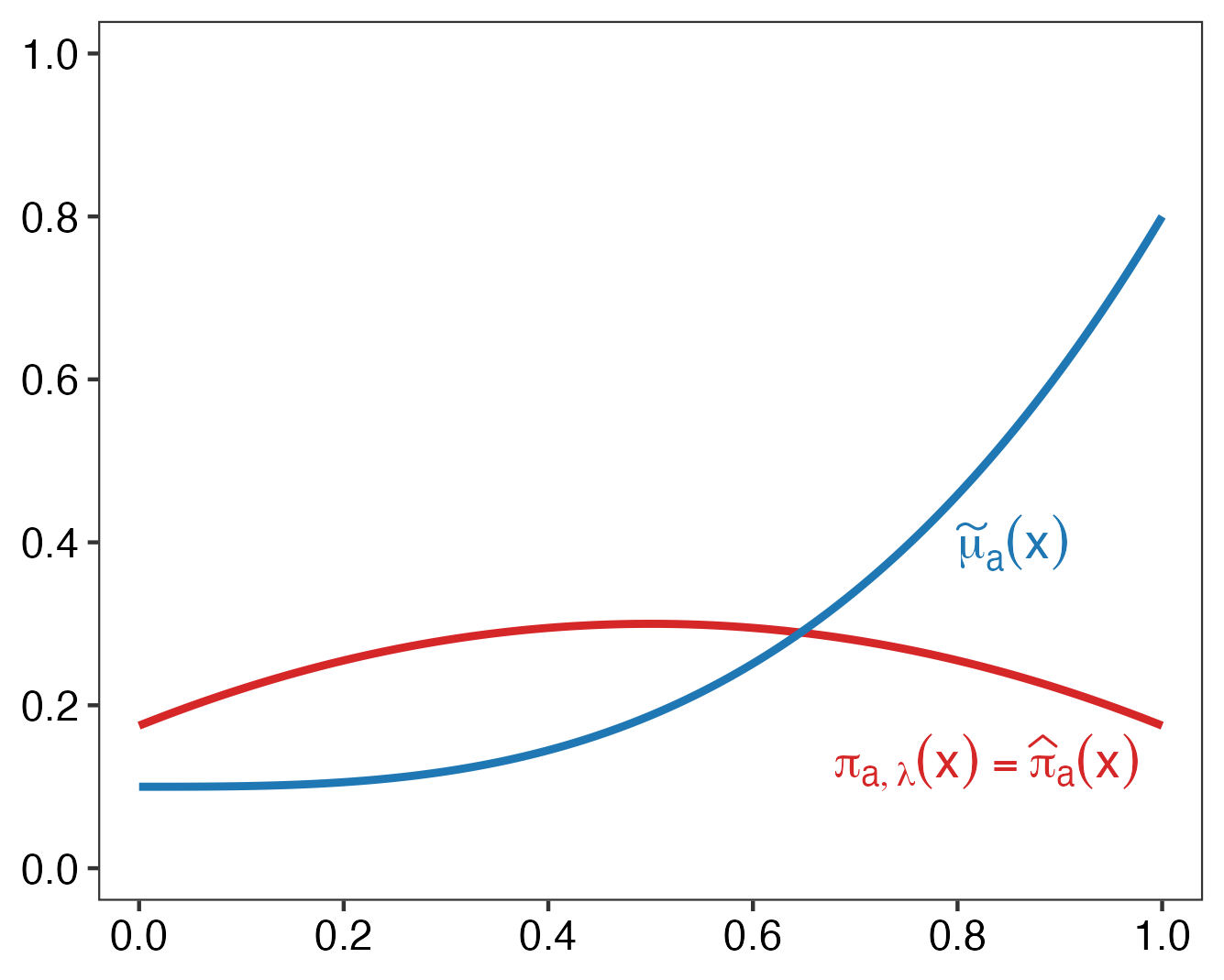}} \\[2em]

        \textbf{Non-active level $a=k$} & 
        \raisebox{-0.5\height}{\includegraphics[width=0.35\textwidth]{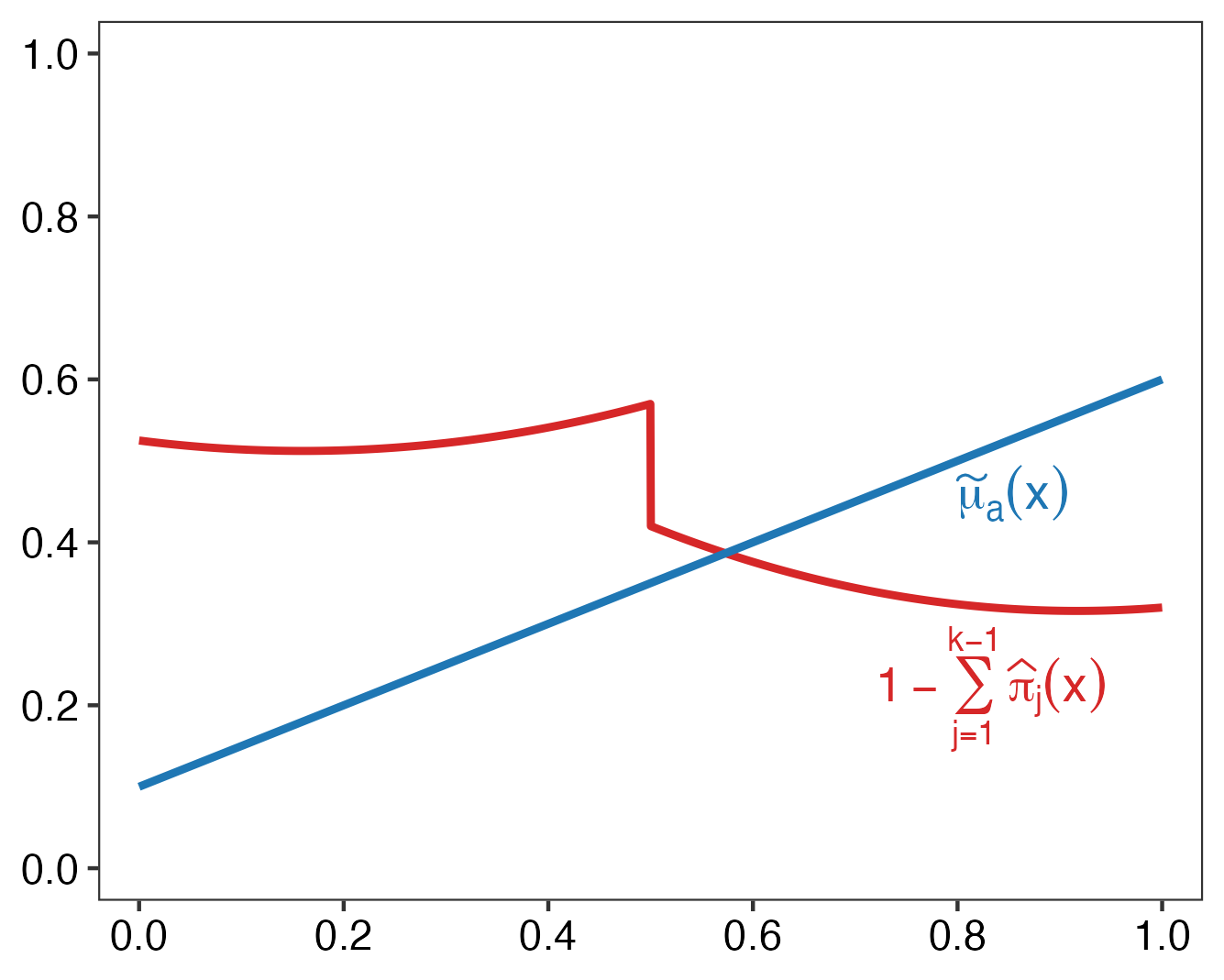}} & 
        \raisebox{-0.5\height}{\includegraphics[width=0.35\textwidth]{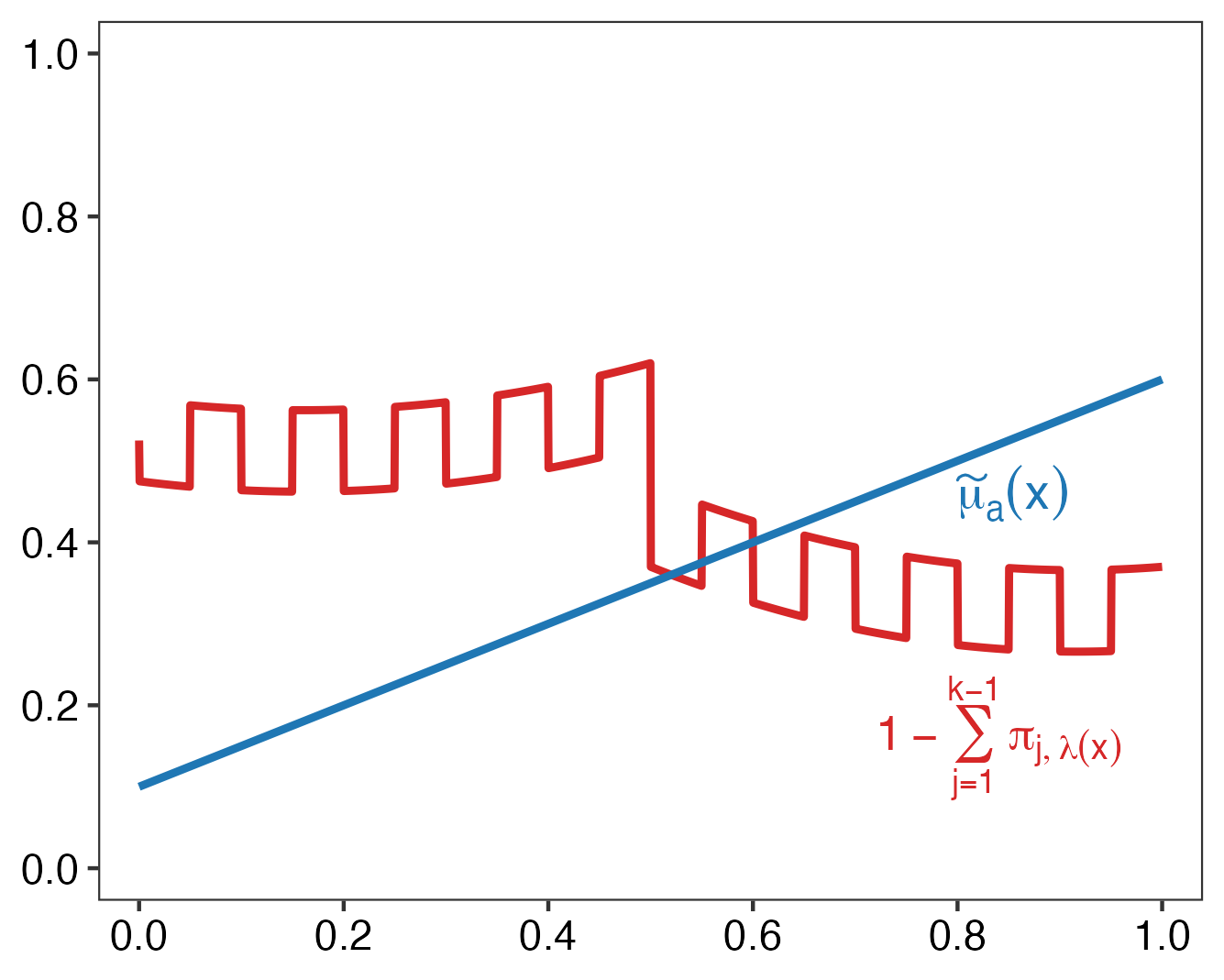}} \\
    \end{tabular}
    \caption{Visualization of the construction of the distributions in the method of fuzzy hypothesis used to prove Theorem \ref{thm:minimax-lb-structure-agnostic}. For active treatment levels, bumps are added to the nuisance functions to obtain the alternative. For non-active levels $a\neq k$, we do not alter the nuisance functions, so the null and alternative distribution coincide. For the remaining treatment level $k$, we leave the regression function unaltered and choose the propensity score so that it sums to 1.}
    \label{fig:minimax-construction}
\end{figure}

Next, we provide more details on how these mixtures are constructed. \cite{jin2025structureagnostic} construct the null distribution through the given nuisance estimates, and the alternative distribution by partitioning the covariate space and adding random asymmetric bumps (indexed by $\lambda$) to the estimates of the propensity score and regression function with height $\delta_n$ and $\epsilon_n$, respectively. In constrast to the setting of \cite{jin2025structureagnostic} with a binary treatment, our construction must additionally consider that we have $k$ different treatment levels and sparsity must be respected. If we simply constructed the alternative distribution by adding bumps to the nuisance estimates at all treatment levels, exact sparsity would be violated, as the functional under this distribution would differ from the intercept at all levels, hence not respecting the model assumptions. Therefore, we construct the alternative distribution by only adding bumps to the nuisance function for active treatment levels, i.e., all levels in an active set $S=\{a\mid \psi_a\neq\theta\}$ of cardinality at most $s$ (the sparsity constraint). At non-active levels, we do not alter the nuisance functions (except for the propensity score at the last level $k$, which we choose to normalize the sum across all levels to one). Then, the model assumptions are respected by the alternative distribution, and the functional differs from the intercept $\theta$ on at most $s$ treatment levels. For a visualization of this construction, we refer to Figure \ref{fig:minimax-construction}. The functional seperation of magnitude $(\delta_n\epsilon_n)^2$ at active treatment levels can now be derived as in \cite{jin2025structureagnostic}. Since there are $s$ active levels, we overall obtain a functional separation of $s(\delta_n\epsilon_n)^2$, which gives the desired lower bound after scaling by $k$.

\section{Extensions}\label{section:extensions}

In this section, we first generalize the theory developed in Section \ref{section:single-treatment} to cover both single multivalued and vector treatments, as well as arbitrary target parameters beyond mean potential outcomes. Then we apply this general framework to propose estimators for mean potential outcomes under high-dimensional vector treatments.

\subsection{Sparse Pseudo-Outcome Regression Framework}\label{section:pseudo-outcome-regression}

In this section, we propose a general sparse pseudo-outcome regression framework which can be used to obtain fast convergence rates for the estimation of an arbitrary high-dimensional statistical functional $\psi\in\mathbb{R}^k$ associated with level-specific outcomes $f(Z;a)$ for some (fixed) function $f$, where the levels $a$ are given by the possible values of the observed random variable $A$, i.e., $a\in\mathcal{A}=\mathrm{supp}(A)$. We propose constrained Lasso and best subset selection estimators and provide error guarantees under sparsity, showing that efficient estimation is achievable even in high-dimensional regimes.

\subsubsection{Motivation}\label{sec:two-stage-procedures}

The presented framework is motivated by the estimators proposed in Section \ref{section:single-treatment} for estimating mean potential outcomes of a single multivalued treatment. In this case, we relate $\psi$ to these mean potential outcomes for each treatment level, i.e., $\E(Y^a)=\psi_0+\psi_a$ for $a\in\{1,\dots, k\}$. As one might be interested in estimating other target parameters or in the high-dimensional vector treatment setting, we propose a sparse pseudo-outcome regression framework that applies to arbitrary statistical functionals $\psi$ associated with treatment-level-specific outcomes. In particular, it recovers the estimation of mean potential outcomes in the high-dimensional single treatment setting (refer to Section \ref{section:single-treatment}) and vector treatment setting, both binary and continuous (refer to Section~\ref{section:vector-treatment}), as special cases. Although these are our primary applications of the proposed framework, it is essential to note that it can be applied to many more high-dimensional target parameters, such as incremental effects.\\

Our proposed framework can be viewed as a two-stage procedure using sample splitting. In the first stage, on the first fold of the data, we estimate so-called pseudo-outcomes (which represent our estimation target). With the term \textit{pseudo-outcome}, we refer to a function of the observed data that acts as a replacement of the actual (unobserved) outcome and, in particular, is equal to it in expectation. Resorting to a pseudo-outcome is necessary since the actual outcome is often unobserved, e.g., potential outcomes. In the example of estimating $\E(Y^a)$, a suitable pseudo-outcome would be the uncentered efficient influence function of this parameter. Often, such a pseudo-outcome depends on unknown nuisance parameters, such as regression functions or the propensity score. Consequently, we must use an estimated version of the pseudo-outcome, such as the estimated uncentered influence function obtained by plugging in estimates of the nuisance parameters. Throughout this section, as the choice of the pseudo-outcome depends on the estimation target, we assume the estimated pseudo-outcomes as given and we are agnostic about how they are obtained (for an explicit construction refer to Sections \ref{section:single-treatment} and \ref{section:vector-treatment}). In the second stage of our procedure, we now use these estimated pseudo-outcomes in a sparse linear least squares regression on the second fold of the data.\\

Two-stage procedures with estimated pseudo-outcomes using sample splitting, which allow for structure-agnostic estimation, have been studied in the literature before. Our proposed framework is of a similar spirit to the two-stage procedures proposed by \citet{kennedy2023towards} for CATE estimation, where estimated pseudo-outcomes are regressed on covariates, or \citet{foster2019orthogonal}, who present a general empirical risk minimization framework where the loss function can depend on nuisance parameters; both these frameworks allow for a structure-agnostic estimation in both stages. Furthermore, in previous lines of work, including \citet{ai2003efficient} and \citet{rubin2005general}, regression with general pseudo-outcomes has been studied, however, without the use of sample splitting, requiring restrictions on the estimation procedures in the first or second stage. In our framework, we are agnostic about the first-stage procedure, however, explicitly set the second-stage procedure to be a constrained linear least squares regression.\\

Our pseudo-outcome regression framework has several advantages over the procedure of \citet{foster2019orthogonal} when handling high-dimensional treatments. First, strong positivity is necessarily violated in the high-dimensional treatment case, which would be required to obtain sufficiently fast rates from \cite{foster2019orthogonal}. Specifically, the fast rates for plug-in empirical risk minimization in Theorem 3 therein include constants $B_1$ and $B_2$, which depend on the smoothness and convexity of the loss function. However, in our high-dimensional treatment setting, those constants scale with $k$ (e.g., for $B_1$, either the Lipschitz constant $L$ or the inverse convexity parameter $1/\lambda$ scales with $k$, depending on the weighting in the risk function). Therefore, the obtained rate from their theorem would be slower due to this additional $k$-dependence. Second, by exploiting additional structure, we obtain a doubly robust nuisance error term, rather than an $L_4$-error that is not doubly robust, an essential feature of our proposed estimators. Third, the error rates presented here are derived using a different proof technique, which may be of independent interest. In particular, although we analyze a constrained version of our regression estimators, our proofs can be straightforwardly adapted to penalized versions (whereas the proof of Example 3 in \cite{foster2019orthogonal} relies on bounding the critical radius of the function class of high-dimensional linear predictors with sparsity constraint; instead, to switch to the penalized version, Theorems 1 and 2 could be invoked, however, requiring an additional derivation of the error rate for the penalized Lasso as second-stage procedure).\\

In the following subsection, we introduce the sparse pseudo-outcome regression in its most general form, allowing for arbitrary estimated pseudo-outcomes and their related high-dimensional statistical functional $\psi$.

\subsubsection{Proposed Estimators}

Let $\psi\in\mathbb{R}^k$ denote any high-dimensional statistical functional throughout this section, which is our estimation target. Further, let $\mathcal{A}=\mathrm{supp}(A)$.\\

Let $\widehat\varphi_1(Z)$ and $\widehat\varphi_{2,a}(Z)$ be generic estimated pseudo-outcomes for each $a\in\mathcal{A}$. These could, for instance, be obtained by constructing pseudo-outcomes $\varphi_1(Z;\eta)$ and $\varphi_{2,a}(Z;\eta)$ (e.g., based on the uncentered efficient influence function of the target parameter) that depend on nuisance parameters $\eta$, and defining the estimated pseudo-outcome as plug-in estimates $\widehat\varphi_1(Z)=\varphi_1(Z;\widehat\eta)$ and $\widehat\varphi_{2,a}(Z)=\varphi_{2,a}(Z;\widehat\eta)$ for a given estimate $\widehat\eta$ of the nuisance parameters. We use sample splitting and assume that the estimation of the functions $\widehat\varphi_1$ and $\widehat\varphi_{2,a}$ is done on a separate independent fold $D_1$ of size $n$; in particular, independent of our sample $(Z_1,\dots,Z_n)$ which we denote by $D_2$. \\

Furthermore, let $V_a\in\mathbb{R}^k$ denote deterministic vectors for all $a\in\mathcal{A}$. We regress the estimated pseudo-outcomes on these predictors. The explicit form of $V_a$ is determined by the model assumption that connects $\psi$ to the level-specific outcome.\\

Finally, let $\widehat\Pb$ and $\widetilde\Pb$ be measures on $\mathcal{A}$, the support of the treatment variable. $\widehat\Pb$ is allowed to depend on the sample $D_1$, and $\widetilde\Pb$ is allowed to depend on both $D_1$ and $D_2$.\\

Given the previous notation and setup, we now present our proposed estimators. Define the empirical risk
\begin{equation} \label{eq:emp-risk}
    \widehat R(\beta) =  -2\Pb^2_{n}\left\{ \widehat\varphi_1(Z)V_A^T\beta + \int_\mathcal{A}\widehat\varphi_{2,a}(Z)V_a^T\beta d\widehat\Pb(a)\right\} + \int_\mathcal{A}(V_a^T\beta)^2d\widetilde\Pb(a).
\end{equation}
Then, define the Lasso estimator
\begin{equation}\label{def:lasso}
    \widehat\psi_{\text{lasso}} = \argmin_{\norm{\beta}_1\leq s} \widehat R(\beta)
\end{equation}
and the best subset selection estimator
\begin{equation}\label{def:bestsubset}
    \widehat\psi_{\text{subset}} = \argmin_{\norm{\beta}_0\leq s} \widehat R(\beta)
\end{equation}
as the minimizer of this empirical risk subject to an $L_1$- and $L_0$-constraint, respectively. By imposing an $L_0$- or $L_1$-constraint on the coefficient vector in this regression, we aim to leverage sparsity of the vector $\psi$ and achieve fast convergence rates even when $k$ is large in proportion to $n$.\\

\begin{remark}[Recovering the estimators of Section \ref{section:single-treatment}]\label{remark:recovering-single-multivalued-treatment-setting}
    We note that the proposed estimators in Section \ref{section:single-treatment} are special cases of the above best subset selection and Lasso estimators. This can be seen by setting the estimated pseudo-outcomes $\widehat\varphi_{2,a}(Z)$ according to \eqref{def:estimated-pseudo-outcome-single-treatment}, $\widehat\varphi_1(Z)=0$, $V_a=(\mathds{1}(a=1),\dots,\mathds{1}(a=k))^T$ to be the vector of treatment indicators, and $\widehat\Pb=\Pb_n^1$ and $\widetilde\Pb=\Pb_n^2$ to be the empirical measures on $D_1$ and $D_2$. Then the general empirical risk in \eqref{eq:emp-risk} simplifies to the empirical risk in \eqref{eq:emp-risk-single-trt} when estimating mean potential outcomes of a single multi-valued treatment.\\
\end{remark}

In the following subsection, we analyze the proposed estimators under sparsity, providing error guarantees that validate the proposed pseudo-outcome regression framework.

\subsubsection{Master Theorem}

In this section, we provide an error guarantee for our proposed estimators under sparsity. Note that we put no restrictions on how the given pseudo-outcomes are obtained when formulating our estimators. However, for $\widehat\psi_\text{lasso}$ and $\widehat\psi_\text{subset}$ to be accurate estimators, we require that the pseudo-outcomes are appropriately related to the estimation target $\psi$. Therefore, we state an error guarantee for our estimators that depends on individual error rates characterizing the quality of our chosen pseudo-outcomes. Such an error guarantee can be derived under both exact sparsity, i.e., $\norm{\psi}_0\leq s$, or approximate sparsity, i.e., $\norm{\psi}_1\leq s$, for some sparsity constraint $s$. In the following, we provide this error guarantee under exact sparsity. For details on the analysis under approximate sparsity, we refer to Appendix \ref{section:master-theorem-under-approximate-sparsity}.

\begin{theorem}[Error rate under exact sparsity] \label{thm:master}
Let $\widehat\psi\in\{\widehat\psi_\text{lasso}, \widehat\psi_\text{subset}\}$ be either the Lasso estimator defined in (\ref{def:lasso}) or the best subset selection estimator defined in (\ref{def:bestsubset}) minimizing the empirical risk in (\ref{eq:emp-risk}). Suppose exact sparsity, i.e., $\norm{\psi}_0\leq s$ for some sparsity constraint $s$. Moreover, assume for all $j\in\{1,\dots,k\}$ and $a\in\mathcal{A}$ that
\begin{enumerate}
    \item[(i)] (Almost sure boundedness) $|\widehat\varphi_1(Z)V_{A,j}|$, $\left| \int_\mathcal{A} \widehat\varphi_{2,a}(Z)V_{a,j} d\widehat\Pb(a) \right|, |V_A^T\psi V_{A,j}|\leq C_1r_1(k)$ almost surely,
    \item[(ii)] (Second moment boundedness) $\Pb\left\{ (\widehat\varphi_1(Z)V_{A,j})^2 \right\}$, $\Pb\left\{\left(\int_\mathcal{A} \widehat\varphi_{2,a}(Z)V_{a,j} d\widehat\Pb(a)\right)^2\right\}$,\\ $(\int_\mathcal{A}V_a^T\psi V_{a,j}d\widehat\Pb(a))^2 \leq C_1\widehat r_1(k)$ almost surely where $\widehat r_1(k)$ can depend on $D_1$ and\\ $\E_{D_1}\{\widehat r_1(k)^p\} \leq Rr_1(k)^p$ for $p=1,2$,
    \item[(iii)] (Nuisance penalty) $\left| \Pb\left\{ \E(\widehat\varphi_1(Z)\mid X, A=a)p_{A\mid X}(a\mid X) + \widehat\varphi_{2,a}(Z) - V_a^T\psi \right\} \right|\leq C_2r_{2,a}(n,k)$ where $p_{A\mid X}$ is the conditional density of $A$ given $X$ with respect to $\widehat\Pb$,
    \item[(iv)] (Restricted eigenvalues)
    $C_3\norm{v}_2^2 \leq v^T\left(\int_\mathcal{A} V_aV_a^T d\widetilde\Pb(a)\right) v$ with probability at least $1-\alpha$ for all $v$ such that $\norm{v}_1 \leq 2\sqrt{s}\norm{v}_2$ where $\alpha\leq C_4\frac{\log k}{n}$,
    \item[(v)] (Difference in measures) $\E\left[\max_{j\in\{1,\dots,k\}}(\widetilde\Pb - \widehat\Pb)\left\{V_A^T\psi V_{A,j}\right\}^2\right]\leq C_5\frac{r_1(k)\log k}{n}$, and
    \item[(vi)] (Integrability) $(\mathcal{A}, \widehat\Pb)$ is $\sigma$-finite and the map
    $$
        (a, Z) \mapsto \left( \E(\widehat\varphi_1(Z)\mid X, A=a)p_{A\mid X}(a\mid X) + \widehat\varphi_{2,a}(Z) - V_a^T\psi \right) V_{a,j}
    $$
    is $\widehat\Pb\otimes\Pb$-integrable
\end{enumerate}
for rates $r_1 = r_1(k)\gtrsim 1, r_2=r_2(n,k)$, and constants $R,C_1,C_2,C_3,C_4,C_5>0$ that are independent of $n$ and $k$.
Then,
\begin{align*}
    \E\left( \int_\mathcal{A}\left\{ V_a^T(\widehat\psi - \psi) \right\}^2 d\widetilde\Pb(a) \right) &\lesssim s\frac{r_1(k)\log k}{n} + s\E\left( \max_{j\in\{1,\dots,k\}} \left(\int_\mathcal{A} r_{2,a}(n,k)|V_{a,j}| d\widehat\Pb(a) \right)^2\right)
\end{align*}
for $r_1(k)\log k \lesssim n$.
\end{theorem}

The preceding theorem states an error rate that consists of two parts: the oracle rate $s\frac{r_1(k)\log k}{n}$ and the error from estimating the pseudo-outcomes given by the second summand. Intuitively, the oracle rate is the error that we would have if we had access to the oracle outcomes and could use those in the regression. Most importantly, this part of the rate depends only on $\log k$ instead of $k$ by leveraging the sparsity constraint on the vector $\psi$, which is desirable especially when $k$ is large in proportion to $n$. The pseudo-outcome estimation error rate is the price we pay for not observing the actual pseudo-outcomes and having to estimate them. Specifically, the rate corresponds to the additional error that arises from the distance of $\widehat\varphi_1(Z)$ and $\widehat\varphi_{2,a}(Z)$ to the outcome $V_a^T\psi$. If these estimated pseudo-outcomes are defined as plug-in estimates $\varphi_1(Z;\widehat\eta)$ and $\varphi_{2,a}(Z;\widehat\eta)$, then this error term corresponds to the error arising from nuisance estimation.\\

In this paper, the main application of the previous theorem is the estimation of the functionals $\int \mu_a(x)d\Pb(x)$ (which equal $\E(Y^a)$ under the usual identification assumptions) for different high-dimensional treatment settings. For instance, we use it to prove Theorems~\ref{thm:single-treatment}, \ref{thm:binary-vector-trt}, and \ref{thm:continuous-vector-trt} as special cases by verifying the conditions of Theorem \ref{thm:master}.\\

We proceed with remarks that discuss the assumptions of the previous theorem in more detail and provide guidance on the choice of a suitable pseudo-outcome.\\

\begin{remark}[Discussion of the assumptions]
    Assumptions (i) and (ii) are mild since the rate $r_1(k)$ can be of any order. Since $r_1(k)$ contributes to the oracle part of the overall rate, one should aim to verify these assumptions with $r_1(k)$ as small as possible. In the same way, Assumption (iii) is mild. When applied in the single or vector treatment setting, the rate $r_{2,a}(n,k)$ corresponds to a second-order nuisance penalty. More specifically, it is the product of the errors from estimating the regression function and the propensity score. The restricted eigenvalue condition in Assumption (iv) is a standard assumption in the high-dimensional regression literature to obtain so-called fast Lasso rates in the random design, i.e., a dependence on $\frac{\log k}{n}$ instead of $\sqrt{\frac{\log k}{n}}$. When we apply the above theorem to single multi-valued and vector treatments, respectively, we show that this assumption can be satisfied. We also show that Assumptions (v) and (vi) are satisfiable for appropriate choices of $\widehat\Pb$ and $\widetilde\Pb$. Finally, for comments on the choice of the sparsity constraint $s$, we refer to Remark \ref{remark:sparsity-constraint}, and for a discussion of exact versus approximate sparsity to Remark \ref{remark:exact-versus-approximate-sparsity}.\\
\end{remark}

\begin{remark}[Discrete version]\label{remark:discrete-case}
    When the treatment variable $A$ is discrete, we can simplify the general empirical risk in \eqref{eq:emp-risk} and state a special case of the above theorem. For generic estimated pseudo-outcomes $\widehat\varphi_a(Z)$, the discrete version of the general empirical risk can be written as
    \begin{equation} \label{eq:emp-risk-discrete}
        \widehat R(\beta) =  \Pb^2_{n}\left\{ -2\sum_{a\in\mathcal{A}}\widehat\varphi_a(Z)V_a^T\beta \widehat\varpi_a+ (V_A^T\beta)^2\right\}
    \end{equation}
    where $\widehat\varpi_j=\Pb^1_{n}\{\mathds{1}(A=j)\}$. For further details, we refer to Appendix \ref{sec:discrete-master-thm}. \\
\end{remark}

\begin{remark}[Choosing a suitable pseudo-outcome]
    The pseudo-outcomes $\widehat\varphi_1(Z)$ and $\widehat\varphi_{2,a}(Z)$ can often be chosen based on the (uncentered) efficient influence function of the target parameter. For instance, suppose that the estimation target (oracle outcome) is $\E\{\mu_a(X)\}$, then the uncentered efficient influence function $\varphi_a(Z; \pi, \mu)$ (depending on the propensity score and regression function) equals the outcome in expectation. Using estimators $\widehat\pi$ and $\widehat\mu$ of the nuisance parameters, we obtain an estimated efficient influence function $\varphi_a(Z;\widehat\pi,\widehat\mu)$. For instance, in the case when the treatment is discrete, we can now set $\widehat\varphi_{2,a}(Z)=\varphi_a(Z;\widehat\mu,\widehat\pi)$ and $\widehat\varphi_1(Z)=0$. Using efficient influence function-based pseudo-outcomes in our proposed framework leads to a nuisance estimation error rate $r_{2,a}(n,k)$ in Assumption (iii) that is a second-order term consisting of squares or products of the nuisance errors.
\end{remark}

\subsection{Vector Treatments}\label{section:vector-treatment}

In this section, we apply the sparse pseudo-outcome regression framework of Section \ref{section:pseudo-outcome-regression} to estimate mean potential outcomes for both binary and continuous vector treatments. We present the proposed Lasso and best subset selection estimator of Section \ref{section:pseudo-outcome-regression} for this specific scenario and provide error guarantees under sparsity. These results demonstrate that the proposed estimators are doubly robust and can achieve fast convergence rates in the high-dimensional treatment regime, provided that additional sparsity is assumed.

\subsubsection{Setup \& Model Assumption}

Throughout this section, we assume that the treatment variable $A$ is a vector treatment, i.e., is a vector $A=(A_1,\dots,A_k)$ consisting of $k$ individual treatments, where the number of treatments $k$ is potentially large. Our goal is to efficiently estimate the mean potential outcomes $\E(Y^a)$, i.e., the outcome if a unit received treatment combination $a=(a_1,\dots,a_k)$.\\

To leverage the pseudo-outcome regression framework, we want to connect our oracle outcomes $\E(Y^a)$ with a sparse statistical functional $\psi$. In this vector treatment setting, we do this by assuming a linear marginal structural model of the form
\begin{equation}\label{eq:marginal-structural-model}
    \E\left(Y^{(a_1,\dots, a_k)}\right) = \int_{\mathbb{R}^d} \E(Y \mid X=x, A_1=a_1,\dots, A_k=a_k)d\Pb(x) = \psi_0 + \sum_{j=1}^k \psi_j a_j\;\text{ where }\;\norm{\psi}\leq s
\end{equation}
for some intercept $\psi_0$ that is assumed to be known for simplicity (refer to Remark \ref{remark:intercept} for a discussion on unknown intercepts) and $\norm{\cdot}$ either the $L_0$- or $L_1$-norm. Note that we require typical identification assumptions for the first equality to hold. Even without identification assumptions, the following theory still applies to estimating $\int \mu_a(x)d\Pb(x)$. Intuitively, the model assumption in (\ref{eq:marginal-structural-model}) says that only a few treatments of every treatment combination have a nontrivial effect on the outcome, and also that the mean potential outcomes do not deviate too much from some intercept value $\psi_0$. Our goal is to estimate the statistical functional $\psi:\mathcal{P}\to\mathbb{R}^k$ using our pseudo-outcome regression framework, leveraging the introduced sparsity to achieve fast convergence rates. \\

In the following, we cover binary and continuous vector treatments separately, as both settings require the general pseudo-outcome regression framework of Section \ref{section:pseudo-outcome-regression} to be applied slightly differently. This is because, in the discrete case, the treatment probabilities can be estimated with their empirical proportions, while, in the continuous case, we have to estimate the treatment density. We start with binary vector treatments in the following section.

\subsubsection{Binary Vector Treatments}\label{sec:binary-vector-treatment}

In this section, we propose estimators for mean potential outcomes for binary vector treatments and analyze their error under sparsity.\\

To estimate $\psi$ in (\ref{eq:marginal-structural-model}), we leverage the pseudo-outcome regression framework, more specifically the Lasso estimator defined in (\ref{def:lasso}) and best subset selection estimator defined in (\ref{def:bestsubset}) for a specific choice of $\widehat\varphi_a$ in the discrete version \eqref{eq:emp-risk-discrete} of the general empirical risk (\ref{eq:emp-risk}) that is minimized. For nuisance parameters $\eta=(\mu,\pi)$, define the pseudo-outcomes
\begin{equation*}
    \varphi_a(Z;\eta) = \frac{\mathds{1}(A=a)}{\pi_a(X)}\left\{Y-\mu_a(X)\right\} + \mu_a(X) - \psi_0.
\end{equation*}
for all $a\in\mathcal{A}$. Note that this is the difference of the uncentered efficient influence function for $\E\{\mu_a(X)\}$ and the intercept, so $\E\{\varphi_a(Z;\eta)\}=\E\{\mu_a(X)\} - \psi_0=a^T\psi$. Then, set
\begin{equation}\label{def:estimated-pseudo-outcome-binary-vector-treatment}
\widehat\varphi_{a}(Z)=\varphi_{a}(Z,\widehat\eta)=\varphi_{a}(Z,\widehat\mu, \widehat\pi)
\end{equation}
to be the plug-in estimator of the pseudo-outcome based on estimated $\widehat\mu, \widehat\pi$ that are constructed on a separate sample $D_1$ that is independent of our sample $D_2=(Z_1,\dots,Z_n)$. Next, we set the features $V_a$ to be the treatment combinations themselves, i.e., $V_a=a$. Plugging these expressions into (\ref{eq:emp-risk-discrete}), our proposed Lasso and best subset selection estimators in (\ref{def:lasso}) and (\ref{def:bestsubset}) then minimize the empirical risk
\begin{equation}\label{eq:emp-risk-binary-vector-trt}
    \widehat R(\beta) = \Pb^2_{n}\left\{-2 \sum_{a\in\{0,1\}^k}\left[ \frac{\mathds{1}(A=a)}{\widehat\pi_a(X)}\{Y-\widehat\mu_a(X)\}  + \widehat\mu_a(X) - \psi_0 \right]a^T\beta \widehat\varpi_a + (A^T\beta)^2\right\}
\end{equation}
subject to an $L_1$- or $L_0$-constraint, respectively, where $\widehat\varpi_a=\Pb^1_{n}\{\mathds{1}(A=a)\}$ denotes the empirical proportion of treatment $a$ on $D_1$. This risk function can be motivated analogously to the risk for single multivalued treatments in Remark \ref{remark:motivation-risk-function-single-trt}.\\

Next, we present the error guarantee for these proposed estimators under exact sparsity, i.e., when $\norm{\psi}_0\leq s$. We show a fast Lasso rate can be achieved and that fast rates are possible even when the vector treatment is high-dimensional. For an analysis under approximate sparsity, we refer to Appendix \ref{section:vector-treatment-approximate-sparsity}.

\begin{theorem}\label{thm:binary-vector-trt}
    Let $\widehat\psi\in\{\widehat\psi_\text{lasso}, \widehat\psi_\text{subset}\}$ be either the Lasso estimator defined in (\ref{def:lasso}) or the best subset selection estimator defined in (\ref{def:bestsubset}) minimizing the empirical risk in (\ref{eq:emp-risk-binary-vector-trt}) for $\widehat\varphi_a$ as in (\ref{def:estimated-pseudo-outcome-binary-vector-treatment}). Assume the model defined in (\ref{eq:marginal-structural-model}) and suppose exact sparsity, i.e., $\norm{\psi}_0\leq s$ for some sparsity constraint $s$. Moreover, assume that for all $a\in\{0,1\}^k$
    \begin{itemize}
        \item[(i)] (Estimated propensity score close to empirical weights) $|\widehat\varpi_a/\widehat\pi_a(X)| \leq B$,
        \item[(ii)] (Restricted eigenvalues) $v^T\E(AA^T)v \geq D\norm{v}_2^2$ for all $v\neq 0$,
        \item[(iii)] (Sample covariance condition) $\sup_{\norm{v}_2=1, \norm{v}_1\leq 2\sqrt{s}} \left|v^T\left(\E(AA^T) - \frac{1}{n}\mathbb{A}^T\mathbb{A}\right)v \right|\leq \frac{D}{2}$ with probability at least $1-\alpha$ for $\alpha\leq C_\alpha\frac{\log k}{n}$ where $\mathbb{A} = \begin{pmatrix}A_1^T \\ \vdots \\ A_n^T\end{pmatrix}$, and
        \item[(iv)] (Boundedness) $|Y|, |\widehat\mu|, |\psi_0| \leq B$
    \end{itemize}
    almost surely for constants $B,C_\alpha, D>0$ that are independent of $k$ and $n$.
    Then,
    \begin{equation*}
        \E\left( \left\{ A^T(\widehat\psi - \psi) \right\}^2 \right) \lesssim s\frac{\log k}{n} + s(\delta_n\epsilon_n)^2
    \end{equation*}
    for $s^4\log k \lesssim n$, when $\delta_n(a)\leq\delta_n$ and $\epsilon_n(a)\leq\epsilon_n$ for all $a\in\{0,1\}^k$ where
    $$
        \delta_n(a) = \anorm{\frac{\pi_a}{\widehat\pi_a} - 1} \quad\text{and}\quad \epsilon_n(a) = \anorm{\mu_a - \widehat\mu_a}.
    $$
\end{theorem}

The error guarantee given by the theorem consists of two parts: First, the oracle rate $s\frac{\log k}{n}$ that one would obtain when using $Y^a$ as outcome in the regression. Second, the nuisance estimation error that arises from estimating the propensity score and regression function, which is the price we pay for not observing the potential outcomes and having to estimate them.\\

To achieve the oracle rate $s\frac{\log k}{n}$, we require the nuisance estimation error to be of a smaller or equal order. Since the nuisance estimation error is a second-order term, i.e., a product of propensity score and regression estimation error, each nuisance can be estimated at a slower rate as long as the product of the errors is fast enough. Importantly, when the oracle rate is achieved, the prediction error of the mean potential outcomes depends on $k$ only logarithmically, showing that consistent estimation is possible even when $k$ is large in proportion to $n$. \\

In the following, we provide remarks on alternative upper bounds and the assumptions of the previous theorem. For a discussion of exact versus approximate sparsity and the setting where the intercept is unknown, we refer to Remarks \ref{remark:exact-versus-approximate-sparsity} and \ref{remark:intercept}.\\

\begin{remark}[Alternative upper bound]
    It is possible to show a more precise upper bound of the form
    $$
        \E\left( \left\{ A^T(\widehat\psi - \psi) \right\}^2 \right) \lesssim s\frac{\log k}{n} + s\E\left( \left( \sum_{a\in \mathcal{A}} \widehat\varpi_a \delta_n(a)\epsilon_n(a) \right)^2\right).
    $$
    This upper bound is established during the proof of the theorem. It shows that the nuisance error at a certain level is weighted by its empirical proportion on $D_1$ (the sample used for nuisance estimation). This particularly shows that the overall error does not suffer from treatment levels that are potentially unseen (i.e., when $\widehat\varpi_a=0$), as nuisance estimation can be arbitrarily inaccurate at such levels. \\
\end{remark}

\begin{remark}[Restricted eigenvalues and sample covariance condition]
The restricted eigenvalue assumption and sample covariance condition are satisfied, e.g., when $A=(A_1,\dots, A_k)\in\{0,1\}^k$ where the components are independent and $\Pb(A_j=1)=1/2$ for all $j$, as long as $n\gtrsim s^3\log k$. For details on the verification of the sample covariance condition, we refer to Appendix \ref{sec:sample-covariance}.
\end{remark}

\subsubsection{Continuous Vector Treatments}\label{sec:continuous-vector-treatment}

In this section, we propose estimators for mean potential outcomes for continuous vector treatments and analyze their error under sparsity.\\

To estimate $\psi$ in (\ref{eq:marginal-structural-model}), we leverage the pseudo-outcome regression framework, more specifically the Lasso estimator defined in (\ref{def:lasso}) and best subset selection estimator defined in (\ref{def:bestsubset}) for specific choices of $\widehat\varphi_1(Z), \widehat\varphi_{2,a}(Z)$ and the measures $\widehat\Pb, \widetilde\Pb$ in the empirical risk (\ref{eq:emp-risk}) that is minimized. For nuisance parameters $\eta=(\mu,\pi)$, define the pseudo-outcomes
\begin{align*}
    \varphi_1(Z;\eta, p) &= \frac{Y-\mu_A(X)}{\pi_A(X)/ p(A)},\\
    \varphi_{2,a}(Z;\eta) &= \mu_a(X) - \psi_0.
\end{align*}
for all $a\in\mathcal{A}$ where $p$ denotes the probability density function of $A$. Note that $\E\{\varphi_1(Z;\eta)\}=0$ and $\E\{\varphi_{2,a}(Z;\eta)\}=\E\{\mu_a(X)\} - \psi_0=a^T\psi$. Then, set
\begin{equation}\label{def:estimated-pseudo-outcome-continuous-vector-treatment}
\begin{split}
\widehat\varphi_1(Z)&=\varphi_1(Z,\widehat\eta, \widehat p)=\varphi_1(Z,\widehat\mu, \widehat\pi, \widehat p),\\
\widehat\varphi_{2,a}(Z)&=\varphi_{2,a}(Z,\widehat\eta)=\varphi_{2,a}(Z,\widehat\mu)
\end{split}
\end{equation}
to be the plug-in estimators of the pseudo-outcomes based on estimated $\widehat\mu, \widehat\pi$, and an estimated treatment density $\widehat p$ that are constructed on a separate sample $D_1$ that is independent of our sample $D_2=(Z_1,\dots,Z_n)$. Next, set $\widehat\Pb = \widetilde \Pb = \widehat p\cdot \lambda$ to be the measure given by the estimated treatment density $\widehat p$. Moreover, set the features $V_a$ to be the treatment combinations themselves, i.e., $V_a=a$. Plugging these expressions into (\ref{eq:emp-risk}), our proposed Lasso and best subset selection estimators in (\ref{def:lasso}) and (\ref{def:bestsubset}) then minimize the empirical risk
\begin{equation}\label{eq:emp-risk-continuous-vector-trt}
    \widehat R(\beta) = -2\Pb^2_{n}\left\{ \frac{Y-\widehat\mu_A(X)}{\widehat\pi_A(X)/ \widehat p(A)}A^T\beta + \int_\mathcal{A}\left( \widehat\mu_a(X) - \psi_0 \right)a^T\beta \widehat p(a)da\right\} + \int_\mathcal{A} (a^T\beta)^2\widehat p(a)da
\end{equation}
subject to an $L_1$- or $L_0$-constraint, respectively. This risk can be motivated similarly as in Remark \ref{remark:motivation-risk-function-single-trt}. In Remark \ref{remark:risk-in-continuous-case}, we discuss an alternative, more complex risk function tailored to the continuous treatment case that has certain advantages over this simpler risk function.\\

Next, we present the error guarantee for these proposed estimators under exact sparsity, i.e., when $\norm{\psi}_0\leq s$. We show that a fast Lasso rate can be achieved and that fast rates are possible even when the vector treatment is high-dimensional. For an analysis under approximate sparsity, we refer to Appendix \ref{section:vector-treatment-approximate-sparsity}.

\begin{theorem}\label{thm:continuous-vector-trt}
    Let $\widehat\psi\in\{\widehat\psi_\text{lasso}, \widehat\psi_\text{subset}\}$ be either the Lasso estimator defined in (\ref{def:lasso}) or the best subset selection estimator defined in (\ref{def:bestsubset}) minimizing the empirical risk in (\ref{eq:emp-risk-continuous-vector-trt}) for $\widehat\varphi_1,\widehat\varphi_{2,a}$ as in (\ref{def:estimated-pseudo-outcome-continuous-vector-treatment}). Assume the model defined in (\ref{eq:marginal-structural-model}) and suppose exact sparsity, i.e., $\norm{\psi}_0\leq s$ for some sparsity constraint $s$. Moreover, assume
    \begin{itemize}
        \item[(i)] (Strong positivity) $\widehat\pi_a(X)\geq\epsilon$ for all $a\in\mathcal{A}$ almost surely,
        \item[(ii)] (Restricted eigenvalues) the eigenvalues of $\int_\mathcal{A} aa^T\widehat p(a)da$ are lower bounded by $C$, and
        \item[(iii)] (Boundedness) $|Y|, |\widehat\mu|, |\psi_0|, \norm{A}_\infty, \widehat p, \int_\mathcal{A} \widehat p(a)da \leq B$ almost surely
    \end{itemize}
    for constants $B,C,\epsilon>0$ that are independent of $k$ and $n$.
    Then,
    \begin{equation*}
        \E\left( \int_\mathcal{A}\left\{ a^T(\widehat\psi - \psi) \right\}^2 \widehat p(a)da\right) \lesssim s\frac{\log k}{n} + s(\delta_n\epsilon_n)^2
    \end{equation*}
    for $\log k \lesssim n$, when $\delta_n(a)\leq\delta_n$ and $\epsilon_n(a)\leq\epsilon_n$ for all $a\in\mathcal{A}$ where
    $$
        \delta_n(a) = \anorm{\frac{\pi_a}{\widehat\pi_a} - 1} \quad\text{and}\quad \epsilon_n(a) = \anorm{\mu_a - \widehat\mu_a}.
    $$
\end{theorem}

For a discussion of the theorem, we refer to Section \ref{sec:binary-vector-treatment} and the discussion of the error rate for binary vector treatments in Theorem \ref{thm:binary-vector-trt} therein, as most of the discussion is analogous. We only point out a few remarks that differ from the binary setting.\\

\begin{remark}[Alternative upper bound]
    It is possible to show a more precise upper bound of the form
    $$
        \E\left( \int_\mathcal{A}\left\{ a^T(\widehat\psi - \psi) \right\}^2 \widehat p(a)da\right) \lesssim s\frac{\log k}{n} + s\E\left( \left( \int_\mathcal{A} \delta_n(a)\epsilon_n(a)\widehat p(a) da \right)^2\right).
    $$
    This upper bound is established during the proof of the theorem. It shows that the nuisance error at a certain level is weighted by its estimated density. This is beneficial as nuisance estimation can be less accurate for low-density regions. \\
\end{remark}

\begin{remark}[Restricted eigenvalue assumption]
    We demonstrate that the restricted eigenvalue assumption can be satisfied by providing an example. Assume that $A\in[0,1]^k$ and $\widehat p=1$. Then, a simple calculation verifies that the eigenvalues of $\int_\mathcal{A} aa^T \widehat p(a)da=\int_{[0,1]^k} aa^Tda$ are lower bounded by $1/12$, satisfying the restricted eigenvalue condition with $C=1/12$.\\
\end{remark}

\begin{remark}[Weights in the risk]\label{remark:risk-in-continuous-case}
    We note that the empirical risk $\widehat R(\beta)$ in \eqref{eq:emp-risk-continuous-vector-trt} can be motivated analogously to discrete treatments, such as in Remark \ref{remark:motivation-risk-function-single-trt}. In particular, in this derivation, we used only parts of the oracle risk's influence function for simplicity. In the discrete treatment setting, the price we pay for omitting some terms in the influence function when constructing the empirical risk is of smaller order than the oracle rate of the estimator, since treatment probabilities can be estimated using empirical proportions, yielding a negligible empirical process term. However, in the continuous case, this approach comes with the caveat that we can only obtain an upper bound on the risk weighted by the estimated density, as provided in Theorem \ref{thm:continuous-vector-trt}. This is because bounding the risk weighted by the true treatment density yields a first-order dependence on the density estimate. To overcome this issue in the continuous treatment case, the empirical risk must be based on the full influence function of the oracle risk. Then, assuming $\psi_0=0$ for simplicity, the empirical risk becomes
    $$
        \widehat{R}(\beta) = \mathbb{P}_n^2 \left\{ -2\widehat{\Psi}_AA^T\beta -2 \frac{Y - \widehat{\mu}_A(X)}{\widehat{\pi}_A(X)/\widehat{p}(A)} A^T\beta - 2\int_{\mathcal{A}} (\widehat{\mu}_a(X)-\widehat\Psi_a) a^T\beta \widehat p(a)da + (A^T\beta)^2 \right\}
    $$
    where $\widehat\Psi_a$ is an estimate of $\E(Y^a)$. Using this alternative risk for our estimators, it is now possible to upper bound the risk $\E( \{ A^T(\widehat\psi - \psi) \}^2 )$ weighted by the true density. Specifically, the error guarantee then contains an additional second-order nuisance dependence involving the error of $\widehat p$ and $\widehat\Psi$. Without providing further details, we note that this can be proved using similar decompositions as in the proof of Theorem \ref{thm:master}.
\end{remark}

\section{Empirical Data Analysis}\label{section:empirical-data-analysis}

In this section, we illustrate our proposed estimators of mean potential outcomes for high-dimensional single multivalued treatments in the setting of Section \ref{section:single-treatment}. We apply our method both to simulated and real data.

\subsection{Simulation Study}\label{section:simulation}

In this section, we illustrate the error guarantees of our proposed estimators on simulated data. While our simulations focus on estimating mean potential outcomes for single multi-valued treatments (i.e., applying the proposed estimators of Section \ref{section:single-treatment}), we note that simulations for vector treatments would produce similar insights. Further, we omit simulations that illustrate the typical story of high-dimensional statistics. For instance, \cite{wainwright2019highdimensional} conducts simulations that illustrate the importance of considering the dependence of the error on the dimensionality of the data for the case of binary hypothesis testing (refer to Figure 1.1 in \cite{wainwright2019highdimensional}), showing that we would otherwise severely overstate a test's accuracy. Then, they simulate how the low-dimensional theory asymptotically still applies under sparsity (refer to Figure 1.4 in \cite{wainwright2019highdimensional}). We refer to those simulations for an illustration that low-dimensional theory (not considering the dependence on the number of treatments) would severely overstate an estimator's accuracy when the treatment is high-dimensional, and that we can asymptotically achieve the low-dimensional prediction for accuracy once we assume sparsity. Instead, we aim to conduct a simulation that ties high-dimensionality specifically to the causal inference setting. The following subsection outlines our simulation setup.

\subsubsection{Setup and Motivation}

We simulate data according to the following data generating process. We sample the covariates from a uniform distribution, i.e., $X\sim\mathrm{Unif}[-1,1]$. We set the propensity score to be perfectly uniform, i.e., $\pi_a(x)=1/k$ for all $a\in\{1,\dots,k\}$. Set the intercept to be $\psi_0=2$ in the model assumption in \eqref{eq:saturated-model}. Moreover, let $\beta_a\sim\mathrm{Unif}[0.25,2], a=1,\dots,s$, be randomly drawn effects for all active treatment levels. Then, we set the outcome regression to be
$$
    \mu_a(x)=\begin{cases}
        \psi_0 & 0\leq a\leq k-s \\
        \psi_0 + \beta_{a-k+s} & k-s+1\leq a \leq k
    \end{cases}
$$
where $s$ is the sparsity constraint. This ensures exact sparsity $\norm{\psi}_0\leq s$, as
$$
    \psi_a = \E(Y^a)-\psi_0 = \E\{\mu_a(X)\}-\psi_0 = \begin{cases}
        0 & 0\leq a\leq k-s \\
        \beta_{a-k+s} & k-s+1\leq a \leq k.
    \end{cases}
$$
Finally, we generate the outcomes according to $Y|X=x,A=a \sim N(\mu_a(x), 1)$. We simulate the following two settings:
\begin{itemize}
    \item Setting 1 (high-dimensional treatment with sparsity): $n=1000, k=100, s=2$
    \item Setting 2 (high-dimensional treatment without sparsity): $n=1000, k=100, s=100$
\end{itemize}
On the one hand, Setting 1 serves to illustrate how our proposed estimator leverages sparsity and achieves a high accuracy. On the other hand, Setting 2 aims to reveal the cost of using our proposed estimator when the effects are not sparse.\\

In our simulation, we compare the following estimators with each other:
\begin{itemize}
    \item Thresholded DR estimator: This is our proposed Lasso estimator in \eqref{def:lasso-single-treatment} of Section \ref{section:single-treatment}, implemented using the formulation as a soft-thresholding estimator (see Appendix \ref{sec:gaussian-sequence-model}).
    \item Thresholded IPW estimator: This is the estimator in \eqref{def:lasso-single-treatment}, but instead we use $\widehat\varphi_a(Z)=\frac{\mathds{1}(A=a)}{\widehat\pi_a(X)}Y-\psi_0$ as estimated pseudo-outcome.
    \item Thresholded plug-in estimator: This is the Lasso estimator in \eqref{def:lasso-single-treatment}, but instead we use $\widehat\varphi_a(Z)=\widehat\mu_a(X)-\psi_0$ as estimated pseudo-outcome.
    \item DR estimator $\Pb_n\left\{\frac{\mathds{1}(A=a)}{\widehat\pi_a(X)}(Y - \widehat\mu_a(X)) + \widehat\mu_a(X)\right\}$
    \item IPW estimator $\Pb_n\left\{\frac{\mathds{1}(A=a)}{\widehat\pi_a(X)}Y\right\}$
    \item Plug-in estimator $\Pb_n\left\{\widehat\mu_a(X)\right\}$
\end{itemize}
We construct the nuisance estimators $\widehat\pi$ and $\widehat\mu$ based on the true nuisance functions to satisfy given (known) $L_2(\Pb)$ error guarantees controlled by the parameters $\alpha$ and $\beta$. This allows us to draw conclusions about the impact of nuisance estimation on the overall error guarantee of the above estimators. Specifically, we set
$$
    \widehat\pi_a(x)=\expit\left(\logit\left(\pi_a(x)\right) + \varepsilon\right) \quad \text{and} \quad \widehat\mu_a(x)=\mu_a(x) + \varepsilon'
$$
where $\varepsilon\sim N\left( \left(\frac{n}{k}\right)^{-\alpha}, \left(\frac{n}{k}\right)^{-2\alpha} \right)$ and $\varepsilon'\sim N\left( \left(\frac{n}{k}\right)^{-\beta}, \left(\frac{n}{k}\right)^{-2\beta} \right)$. This ensures
$$
    \anorm{\frac{\pi_a}{\widehat\pi_a} - 1}\lesssim \left(\frac{n}{k}\right)^{-\alpha} \quad \text{and} \quad \anorm{\widehat\mu_a - \mu_a}\lesssim \left(\frac{n}{k}\right)^{-\beta}.
$$
We further use cross-validation to choose the sparsity constraint $s$ in our thresholded estimators (refer to Appendix \ref{sec:cross-validation} for details).\\

For both Setting 1 and 2 described above, we can now plot the average mean squared error
$$
    \E\left(\frac{1}{k}\sum_{a=1}^k (\widehat\psi_a - \psi_a)^2\right)
$$
of each estimator (estimated across $N=500$ simulation runs) against different choices of $\alpha$ and $\beta$, where $\widehat\psi_a$ and $\psi_a$ are the estimated and true mean potential outcome, respectively. First, we fix the estimation accuracy of the propensity score with $\alpha=0.25$ and vary the estimation accuracy of the regression function $\beta\in\{0, 0.05, 0.1, 0.15, 0.2, 0.25, 0.3, 0.4, 0.5\}$. Secondly, we fix $\beta=0.25$ and vary $\alpha\in\{0, 0.05, 0.1, 0.15, 0.2, 0.25, 0.3, 0.4, 0.5\}$. Therefore, we overall obtain two plots each for Setting 1 and 2.

\subsubsection{Results and Interpretation}

\begin{figure}[ht!]
     \centering
     
     \begin{subfigure}[b]{0.49\textwidth}
         \centering
         \includegraphics[width=\textwidth]{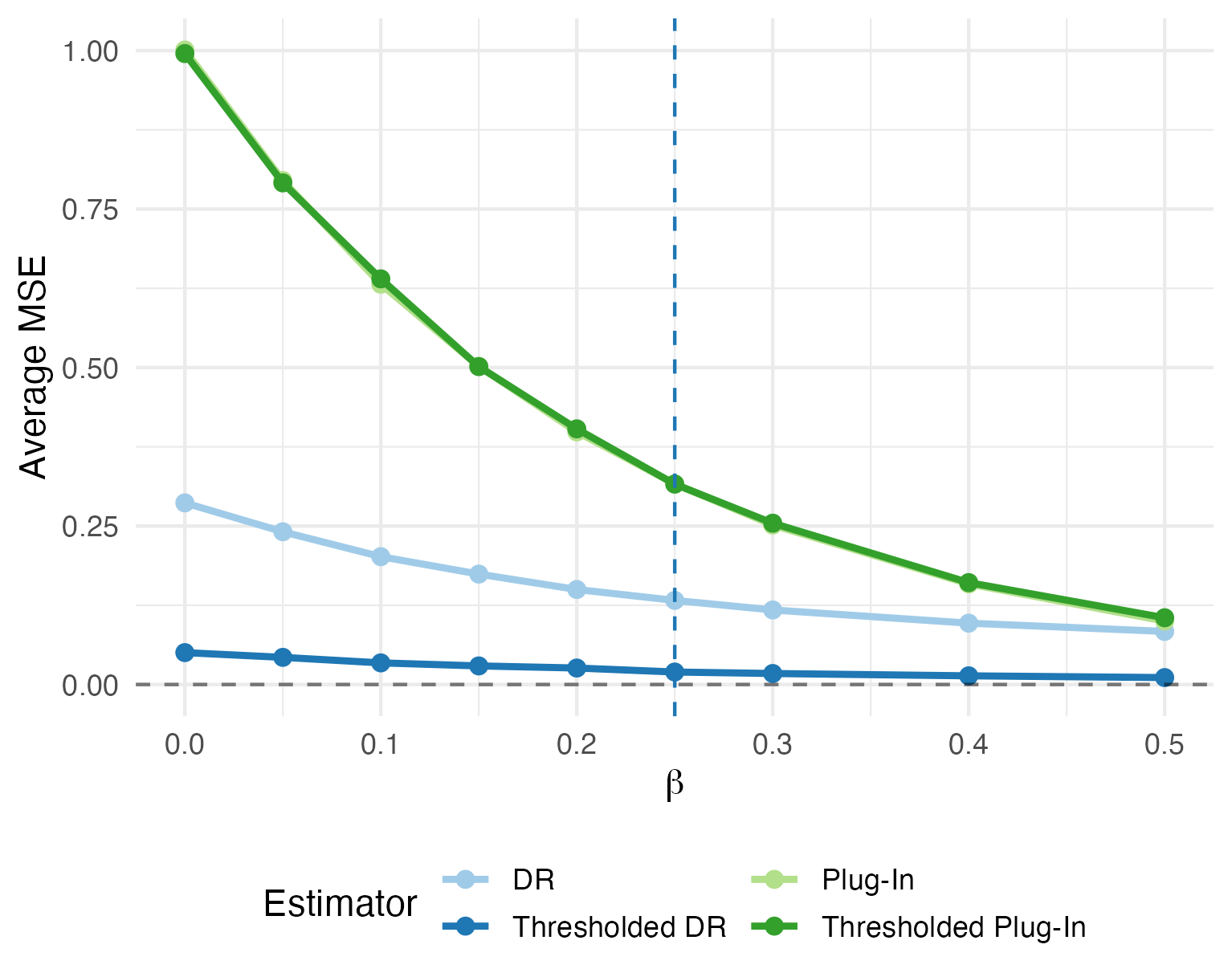}
         \caption{Setting 1 ($n=1000, k=100, s=2$): Plot of MSE against $\beta$ with fixed $\alpha=0.25$.}
         \label{fig:simulation_results_k100_s2_regression}
     \end{subfigure}
     \hfill
     \begin{subfigure}[b]{0.49\textwidth}
         \centering
         \includegraphics[width=\textwidth]{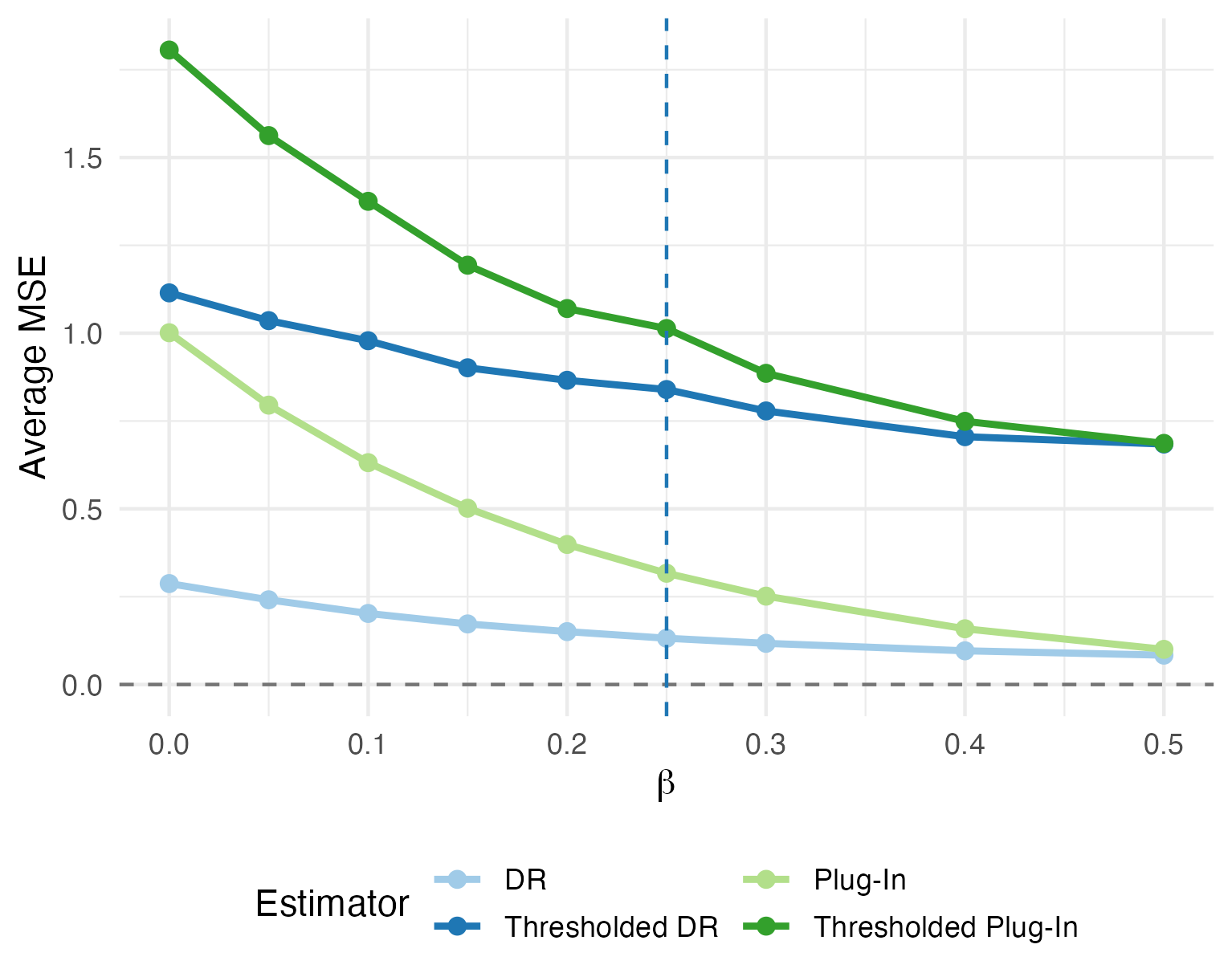}
         \caption{Setting 2 ($n=1000, k=100, s=100$): Plot of MSE against $\beta$ with fixed $\alpha=0.25$.}
         \label{fig:simulation_results_k100_s100_regression}
     \end{subfigure}

     \vspace{10pt}

     \begin{subfigure}[b]{0.49\textwidth}
         \centering
         \includegraphics[width=\textwidth]{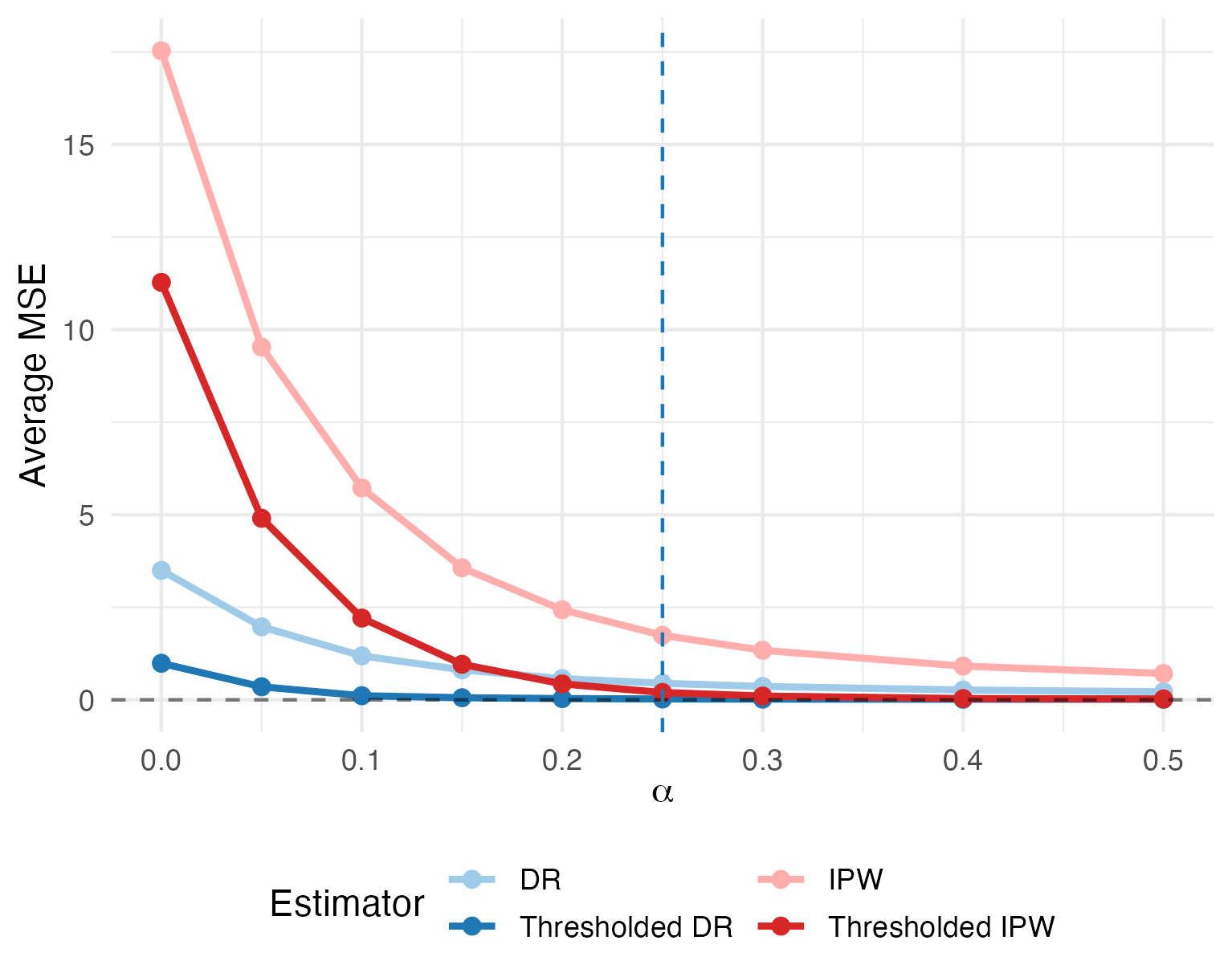}
         \caption{Setting 1 ($n=1000, k=100, s=2$): Plot of MSE against $\alpha$ with fixed $\beta=0.25$.}
         \label{fig:simulation_results_k100_s2_propensity_score}
     \end{subfigure}
     \hfill
     \begin{subfigure}[b]{0.49\textwidth}
         \centering
         \includegraphics[width=\textwidth]{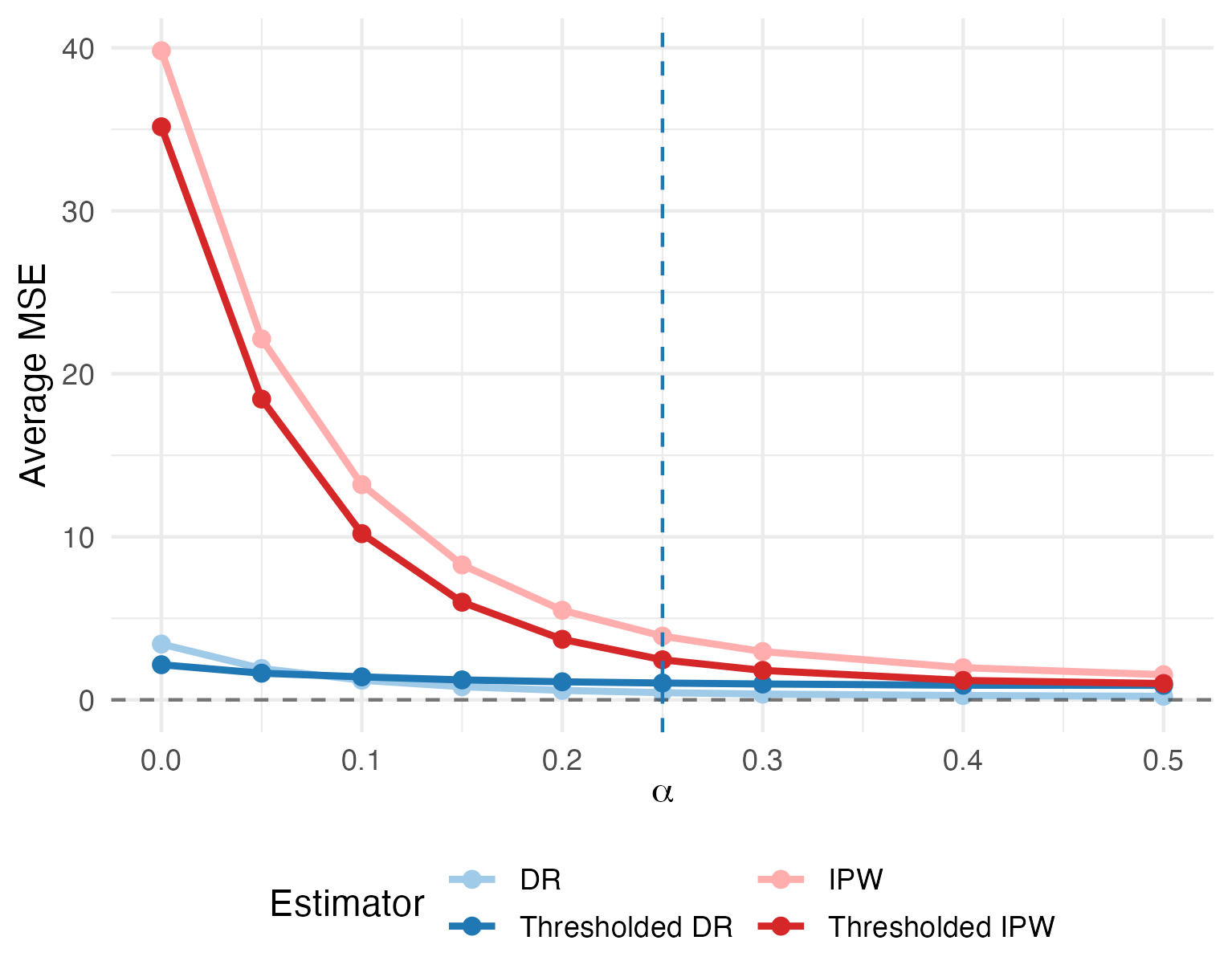}
         \caption{Setting 2 ($n=1000, k=100, s=100$): Plot of MSE against $\alpha$ with fixed $\beta=0.25$.}
         \label{fig:simulation_results_k100_s100_propensity_score}
     \end{subfigure}

     \caption{Simulation results. Plots of the average mean squared error for each estimator against the parameter that controls the nuisance estimation error. $\alpha$ controls the accuracy of the propensity score estimation and $\beta$ of the regression function estimation.}
     \label{fig:simulation_results}
\end{figure}

Figure \ref{fig:simulation_results} shows our simulation results. First of all, we can see from Plot \ref{fig:simulation_results_k100_s2_regression} and \ref{fig:simulation_results_k100_s2_propensity_score} that the thresholded DR estimator leverages the sparsity and achieves a better accuracy than the classical DR estimator across all choices of $\alpha$ and $\beta$. We can see that the classical DR estimator achieves the oracle rate of order $k/n$ once $\alpha=0.25$ and $\beta=0.25$, which is larger than the oracle rate $s\frac{\log k}{n}$ of the thresholded DR estimator in the sparse setting $s=2$.\\

On the contrary, Plots \ref{fig:simulation_results_k100_s100_regression} and \ref{fig:simulation_results_k100_s100_propensity_score} reveal that we pay a small price in terms of accuracy when using our proposed thresholded DR estimator in a high-dimensional setting when there is no sparsity at all (due to choosing the sparsity constraint based on cross-validation). However, we note that this appears not concerning in practice, as meaningful estimation in this setting without any sparsity is impossible regardless of the used procedure. Furthermore, we can see that the overall magnitude of the errors of the thresholded estimators are larger compared to Setting 1 with sparsity, which illustrates the dependence of the oracle rate $s\frac{\log k}{n}$ on the sparsity $s$.\\

Finally, all four plots in Figure \ref{fig:simulation_results} show that our proposed thresholded DR estimator is doubly robust, i.e., the nuisance penalty depends on the squared product of the error from estimating the regression function and propensity score. For example, from Plot \ref{fig:simulation_results_k100_s2_regression} and \ref{fig:simulation_results_k100_s100_regression}, we can see that even when the regression function is highly mispecified (i.e., $\beta$ is small), the thresholded DR estimator performs much better than the thresholded plug-in estimator (which suffers severely from the large estimation error of the regression function). Plots \ref{fig:simulation_results_k100_s2_propensity_score} and \ref{fig:simulation_results_k100_s100_propensity_score} tell an analogous story for the case when the propensity score is estimated poorly.\\

Overall, our simulation illustrates that in a high-dimensional setting with sparsity, the proposed thresholded DR estimator efficiently leverages sparsity and outperforms classical DR, IPW, and plug-in estimators, as well as thresholded versions of the plug-in and IPW estimators, which aligns with our presented theoretical error guarantee in Section \ref{section:single-treatment} and no free lunch results in Appendix \ref{sec:no-free-lunch}.

\subsection{Estimating Effects of Schools on Standardized Test Scores}\label{sec:data-analysis}

We apply the proposed method of Section \ref{section:single-treatment} to estimate the effect of public schools in North Carolina on standardized End-of-Grade test scores (\citet{eogtests}) using a dataset from the North Carolina Education Research Data Center. Full analysis details and results will appear in an upcoming version of the paper.

\section{Discussion}\label{section:discussion}

In this work, we studied optimal estimation of mean potential outcomes for high-dimensional treatments. We proposed estimators that achieve fast convergence rates and allow for a consistent estimation under sparsity assumptions even when the number of treatments is large relative to the sample size. Furthermore, we proved minimax optimality of the proposed estimators in the sparse, structure-agnostic regime for single multi-valued treatments.\\

We also shed light on the differences between various high-dimensional treatment settings, including single multi-valued treatments, as well as binary and continuous vector treatments. We identified fundamental challenges inherent in these treatment settings: First, our target parameter is a high-dimensional vector, making estimation increasingly difficult for many treatments. Second, for discrete high-dimensional treatments, we additionally encounter inevitable positivity violations. Finally, for single multi-valued treatments, we also face the challenge that the effective sample size is $n/k$ instead of $n$, as we only have approximately $n/k$ observations at each treatment level, even when the treatment is perfectly uniform. Despite these three fundamental challenges, our proposed estimators achieve fast convergence rates under sparsity.\\

We unified the estimation procedure across all treatment settings into a single general sparse pseudo-outcome regression framework, which can be applied to estimate mean potential outcomes as a special case. This proposed pseudo-outcome regression framework is applicable in general beyond our work here. It can be used to estimate arbitrary sparse high-dimensional statistical functionals by regressing an appropriately constructed estimated pseudo-outcome on features subject to a sparsity constraint, enabling efficient estimation with fast convergence rates. In future work, our framework could be applied to other high-dimensional target parameters, e.g., indirectly standardized excess risk or outcome ratio \citep{susmann2024doublyrobust}. This can be achieved by formulating a sparse model for the excess risk or outcome ratio, and utilizing the efficient influence functions of these parameters to construct suitable pseudo-outcomes. Finally, the assumptions of Theorem~\ref{thm:master} can be verified, confirming that our framework yields a doubly robust efficient estimator in this high-dimensional indirect standardization setting.\\

For future work, our results have numerous other possible extensions: First, our paper proves the optimality of the proposed estimators for single multi-valued treatments under exact sparsity. Similarly, one could verify optimality under approximate sparsity and in the vector treatment setting. Second, one could investigate whether our proposed pseudo-outcome regression framework can be extended to sparse generalized linear models, thereby enabling the modeling of nonlinear relationships in the marginal structural model used in the vector treatment setting. Finally, valid inference and hypothesis testing could be of future interest.

\section*{Acknowledgements}

EHK was supported by NSF CAREER Award 2047444.

\section*{References}
\vspace{-1cm}
\bibliographystyle{abbrvnat}
\bibliography{bibliography}

\pagebreak

\setlength{\parindent}{0cm}
\appendix

\section{Additional Results}

\begin{figure*}[ht]
    \centering
    \includegraphics[width=0.75\linewidth]{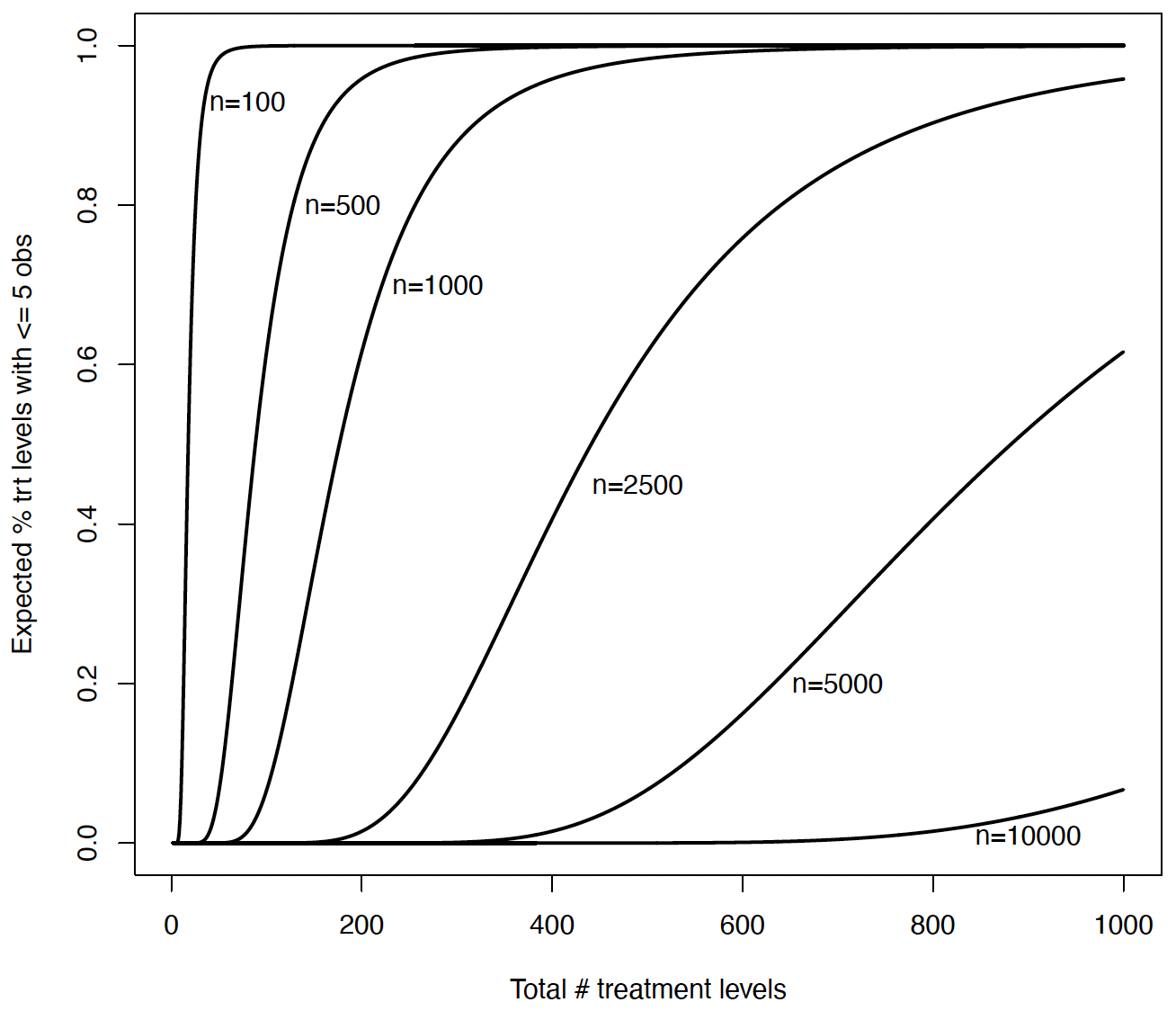}
    \caption{Plot of the expected percentage of treatment levels with less than 5 observations against the sample size $n$ and number of treatments, assuming that the single multi-valued treatment is uniformly distributed.}
    \label{fig:treatment-levels-with-few-obs}
\end{figure*}

\subsection{No Free Lunch}\label{sec:no-free-lunch}

In this subsection, we study the fundamental limits on how well mean potential outcomes can possibly be estimated for discrete high-dimensional treatments. Specifically, we provide two minimax lower bounds for single multi-valued treatments and vector treatments, respectively, under no further structure.\\

First, we present the minimax lower bound in the single multi-valued treatment setting for estimating mean potential outcomes in terms of mean squared error when assuming no additional structure. Intuitively, the minimax risk is the error rate that the best estimator can achieve (in the worst case). Therefore, a lower bound on the minimax risk indicates that no estimator can perform better than the given risk in this minimax sense, i.e., every estimator must achieve the same or a worse error rate. Minimax rates have crucial practical and theoretical implications: On the one hand, they provide a benchmark for the best possible estimation. On the other hand, they precisely determine the statistical difficulty of a certain estimation problem.

\begin{proposition}[Minimax risk in dense regime for single multi-valued treatments] \label{prop:classiclowerbd}
Let $A\in\{1,\dots,k\}$ be a single-multivalued treatment and define $\psi(P)=(\E_P(\mu_1(X)), \dots, \E_P(\mu_k(X)))\in\mathbb{R}^k$. Let $\mathcal{P}$ denote the set of distributions for which $\pi_a(x)\geq C'/k$ for a constant $C'$, $Y$ is binary, and $(A,Y) \ind X$.
If $k \geq 32$, then we have
$$ \inf_{\widehat\psi} \sup_{P \in \mathcal{P}} \E_P\left( \sum_{a=1}^k \varpi_a(P)(\widehat\psi_a - \psi_a(P))^2 \right) \geq   C \min\left( \frac{k}{n}, 1 \right)$$
for $\varpi_a(P)=P(A=a)$ and $C>0$ a universal constant.
\end{proposition}

Essentially, this statement shows that we suffer from very slow convergence rates in this high-dimensional treatment regime, and efficient estimation is hopeless for large $k$. For instance, if $k\asymp n$, then the minimax rate is lower bounded by a rate of constant order; hence, meaningful estimation is impossible, as even with increasingly more data our estimates of the mean potential outcomes would still not improve. This illustrates the challenge of reliable estimation in this high-dimensional regime.\\

Intuitively, the slow convergence rate stated above can be explained as follows: We have inevitable near-violations of positivity (refer to Section \ref{section:positivity-violation}). Even when assuming that each treatment is received with nearly the same probability, we have only around $n/k$ observations per treatment level; thus, the effective sample size for estimating $\E(Y^a)$ is $n/k$ instead of $n$. This explains a lower bound that scales with $k$.\\

Furthermore, note that the stated minimax lower bound is valid even when we assume the propensity scores to be roughly ``perfectly'' distributed in the sense that the propsensity score is nearly uniform, i.e., $C'/k\leq \pi_a(X)\leq C''/k$ for all treatments $a$ with probability one for some constants $C',C''$ that do not depend on $n$ and $k$. This means that the slow convergence rate is a result of the high-dimensionality in the treatment, rather than an ``unfavorable'' distribution of the treatments.\\

Next, we state the minimax lower bound for estimating mean potential outcomes for binary vector treatments, i.e., $A=(A_1,\dots,A_k)\in\{0,1\}^k$ under a linear marginal structural model.

\begin{proposition}[Minimax risk in dense regime for binary vector treatments]\label{prop:classiclowerbd-vector}
Let $A\in\{0,1\}^k$ be a binary vector treatment and define $\psi$ via the linear marginal structural model
\begin{equation*}
    \E_X(\E(Y\mid X, A_1=a_1,\dots, A_k=a_k)) = \psi_0 + \sum_{j=1}^k a_j\psi_j = \psi_0 + a^T\psi.
\end{equation*}
Let $\mathcal{P}$ denote the set of distributions for which $\E(Y\mid A_1=a_1,\dots,A_k=a_k)$ is a linear combination of $a_1,\dots,a_k$ and $(A,Y) \ind X$.
Then we have
$$ \inf_{\widehat\psi} \sup_{P \in \mathcal{P}} \E_P\left( \left\{ A^T(\widehat\psi - \psi(P))\right\}^2 \right) \geq   C \min\left( \frac{k}{n}, k \right)$$
for $C>0$ a universal constant.
\end{proposition}

The above minimax risk indicates that we incur slow error rates, and consistent estimation is impossible when $k$ is large relative to $n$, i.e., when the treatment vector is high-dimensional.\\

Note that under the model $\mathcal{P}$ defined in Proposition \ref{prop:classiclowerbd-vector}, the estimation of $\psi$ reduces to a classical high-dimensional linear regression problem of $Y$ on the high-dimensional features $A=(A_1,\dots,A_k)$. Consequently, Proposition \ref{prop:classiclowerbd-vector} can be verified by using classical minimax results for high-dimensional linear regression in the dense regime when no further structure is assumed, such as Theorem 3 of \citet{tsybakov2003optimal}.\\

In summary, Proposition \ref{prop:classiclowerbd} and Proposition \ref{prop:classiclowerbd-vector} demonstrate that consistent effect estimation for discrete high-dimensional treatments is hopeless without any further assumptions. Consequently, we must impose additional structural assumptions to achieve faster convergence rates.

\subsection{Connection to Gaussian Sequence Model for Single Multi-Valued Treatments} \label{sec:gaussian-sequence-model}

Assume the setting of Section \ref{section:single-treatment}, and, specifically recall the risk in \eqref{eq:emp-risk-single-trt}. We show that minimizing this risk with some $L_1$-penalty (i.e., the Lasso estimator of Section \ref{section:single-treatment}) is equivalent to a soft-thresholding estimator. Let
$$
    \widehat R(\beta) = \sum_{a=1}^k \left(-2\Pb_n^2\left\{\widehat\varphi_a(Z)\right\}\widehat\varpi_a\beta_a + \beta_a^2\Pb_n^2\left\{\mathds{1}(A=a)\right\}\right) + \lambda \sum_{a=1}^k |\beta_a|
$$
be the penalized version of the Lasso estimator (which is equivalent to the thresholded version) where $\widehat\varphi_a(Z)=\varphi_a(Z;\widehat\eta)$ for $\varphi_a(Z;\eta)$ as in Section \ref{section:single-treatment}. Then, to obtain
\begin{align*}
    \widehat\psi = \argmin_\beta\; \widehat R(\beta),
\end{align*}
we can minimize for each component individually, i.e.,
\begin{align*}
    \widehat\psi_a &= \argmin_{\beta_a} \left(-2\Pb_n^2\left\{\widehat\varphi_a(Z)\right\}\widehat\varpi_a\beta_a + \beta_a^2\Pb_n^2\left\{\mathds{1}(A=a)\right\} + \lambda|\beta_a|\right)\\
    &=\argmin_{\beta_a} \left(-2\Pb_n^2\left\{\widehat\varphi_a(Z)\right\}\frac{\widehat\varpi_a}{\Pb_n^2\left\{\mathds{1}(A=a)\right\}}\beta_a + \beta_a^2 + \frac{\lambda}{\Pb_n^2\left\{\mathds{1}(A=a)\right\}} |\beta_a|\right) \\
    &=\argmin_{\beta_a} \left( \left(\beta_a - \Pb_n^2\left\{\widehat\varphi_a(Z)\right\}\frac{\widehat\varpi_a}{\Pb_n^2\left\{\mathds{1}(A=a)\right\}}\right)^2 + \frac{\lambda}{\Pb_n^2\left\{\mathds{1}(A=a)\right\}} |\beta_a|\right).
\end{align*}
The closed-form solution of the previous line is the soft-thresholding estimator, i.e.,
$$
    \widehat\psi_a = \mathrm{sign}\left(\Pb_n^2\left\{\widehat\varphi_a(Z)\right\}\right)\max\left\{ \left| \frac{\widehat\varpi_a}{\Pb_n^2\left\{\mathds{1}(A=a)\right\}}\Pb_n^2\left\{\widehat\varphi_a(Z)\right\} \right|- \frac{\lambda}{2\Pb_n^2\left\{\mathds{1}(A=a)\right\}}, 0 \right\}.
$$

\subsection{Alternative Risk Function}\label{sec:alternative-risk}

Assume the setup and notation of Section \ref{section:single-treatment}. We wish to discuss another risk function given by
\begin{equation}\label{eq:alternative-risk-single-trt}
    \widehat R'(\beta) = \Pb^2_{n}\left\{ \sum_{a=1}^k w(a)\left( \widehat\varphi_a(Z) - \beta_a\right)^2 \right\}
\end{equation}
which uses a fixed weight function instead of empirical proportions as in \eqref{eq:emp-risk-single-trt}. For simplicity, we consider the case where $w=1$: Using a similar analysis as in the proofs of Theorem~\ref{thm:single-treatment} and Theorem~\ref{thm:master}, we can show that the Lasso and best subset selection estimator under this alternative risk function yield an error guarantee where the nuisance penalty depends on $\max_{a\in\{1,\dots,k\}} \anorm{\frac{\pi_a}{\widehat\pi_a} - 1}\norm{\mu_a - \widehat\mu_a}$. If we have unobserved treatment levels, then accurate estimation of the nuisance parameters $\pi_a$ and $\mu_a$ is impossible, and all we can do is guess. However, this is problematic because the nuisance penalty incurs the worst-case nuisance estimation error, as it is the maximum over all treatment levels. Our risk function in \eqref{eq:emp-risk-single-trt} avoids this by weighting the treatments according to the empirical weights on the data used for nuisance estimation (see Remark \ref{remark:weighted-nuisance-error} for details), yielding a nuisance penalty that depends only on observed treatment levels.\\

\begin{remark}[Connection to the Gaussian Sequence Model under this Alternative Risk]
Under the alternative risk outlined in (\ref{eq:alternative-risk-single-trt}), we want to mention that it is possible to establish a connection to the Gaussian sequence model rigorously (analogously as in Appendix \ref{sec:gaussian-sequence-model}). In particular, we can show the equivalence of the Lasso estimator minimizing the alternative risk in (\ref{eq:alternative-risk-single-trt}) and the soft-thresholding estimator, which is given by $\widehat\psi_S=(\widehat\psi_{S,1},\dots,\widehat\psi_{S,k})$ where
\begin{equation*}
    \widehat\psi_{S,a} = \widehat\psi_{S,a}(\lambda) = \mathrm{sign}\left(\Pb_n^2\left\{\widehat\varphi_a(Z)\right\}\right)\left(\left| \Pb_n^2\left\{\widehat\varphi_a(Z)\right\} \right| - \lambda\right)_+.
\end{equation*}
This can be explained as follows: The constrained version of the Lasso, as defined in (\ref{def:lasso}), minimizing the risk in (\ref{eq:alternative-risk-single-trt}), is equivalent to the penalized version of the Lasso. Next, we note that the soft-thresholding estimator can also be written in the penalized form as
\begin{equation*}
    \widehat\psi_{S} = \argmin_{\psi} \frac{1}{2}\sum_{a=1}^k \left( \Pb_n^2\left\{\widehat\varphi_a(Z)\right\} - \psi_a \right)^2 + \lambda\norm{\psi}_1.
\end{equation*}
Then, the equivalence follows by setting the weights $w$ in (\ref{eq:alternative-risk-single-trt}) to be $1/2$, and showing with simple transformations that the minimizer of $\sum_{a=1}^k \left( \Pb_n^2\left\{\widehat\varphi_a(Z)\right\} - \psi_a \right)^2$ is the minimizer of the expression $\Pb_n^2\left\{\sum_{a=1}^k \left( \widehat\varphi_a(Z) - \psi_a \right)^2 \right\}$.
\end{remark}

\subsection{Discrete Version of the Master Theorem}\label{sec:discrete-master-thm}

In this section, we discuss the special cases of Theorem \ref{thm:master} when the treatment variable is discrete. In this case, we can set $\widehat\varphi_1(Z)=0$, $\widehat\varphi_{2,a}(Z)=\widehat\varphi_a(Z)$ for some estimated pseudo-outcome $\widehat\varphi_a(Z)$, and $\widehat\Pb=\Pb_n^1$, $\widetilde\Pb = \Pb_n^2$ in the general empirical risk in \eqref{eq:emp-risk}. Then, the empirical risk we minimize subject to an $L_0$ or $L_1$ constraint takes the form
\begin{equation} \label{eq:emp-risk-discrete-2}
    \widehat R(\beta) =  \Pb^2_{n}\left\{ -2\sum_{a\in\mathcal{A}}\widehat\varphi_a(Z)V_a^T\beta \widehat\varpi_a+ (V_A^T\beta)^2\right\}
\end{equation}
where $\widehat\varpi_j=\Pb^1_{n}\{\mathds{1}(A=j)\}$. Now, we formulate the discrete version of Theorem \ref{thm:master} for our proposed estimators minimizing this specialized risk function.

\begin{theorem}[Discrete version of Theorem \ref{thm:master}] \label{thm:master-discrete}
Let $\widehat\psi\in\{\widehat\psi_\text{lasso}, \widehat\psi_\text{subset}\}$ be either the Lasso estimator defined in (\ref{def:lasso}) or the best subset selection estimator defined in (\ref{def:bestsubset}) minimizing the empirical risk in (\ref{eq:emp-risk-discrete-2}). Suppose exact sparsity, i.e., $\norm{\psi}_0\leq s$ for some sparsity constraint $s$. Moreover, assume for all $j,l\in\{1,\dots,k\}$ and $a\in\mathcal{A}$ that
\begin{enumerate}
    \item[(i)] (Almost sure boundedness) $\left| \sum_{a\in\mathcal{A}} \widehat\varphi_{a}(Z) V_{a,j}\widehat\varpi_a \right|, |V_A^T\psi V_{A,j}|, |V_{A,j}V_{A,l}|\leq C_1r_1(k)$ almost surely,
    \item[(ii)] (Second moment boundedness) $\Pb\left\{\left(\sum_{a\in\mathcal{A}} \widehat\varphi_{a}(Z) V_{a,j}\widehat\varpi_a\right)^2\right\}, (\sum_{a\in\mathcal{A}}V_a^T\psi V_{a,j}\widehat\varpi_a)^2 \leq C_1\widehat r_1(k)$ and $\Pb\left\{(V_A^T\psi V_{A,j})^2\right\}, \Pb(V_{A,j}^2V_{A,l}^2)\leq C_1r_1(k)$ almost surely where $\widehat r_1(k)$ can depend on $D_1$ and $\E_{D_1}\{\widehat r_1(k)^p\} \leq Rr_1(k)^p$ for $p=1,2$,
    \item[(iii)] (Nuisance estimation error rate) $\left| \Pb\left\{ \widehat\varphi_{a}(Z) - V_a^T\psi \right\} \right|\leq C_2r_{2,a}(n,k)$,
    \item[(iv)] (Restricted eigenvalue condition)
    $C_3\norm{v}_2^2 \leq v^T\Sigma v$ for $\Sigma = \E(V_AV_A^T)$, and
    \item[(v)] (Sample covariance condition) $\sup_{\norm{v}_2=1, \norm{v}_1\leq 2\sqrt{s}} \left| v^T\widehat\Sigma v - v^T\Sigma v \right| \leq C_3/2$ with probability at least $1-\alpha$ for $\alpha\leq C_4\frac{\log k}{n}$ where $\widehat\Sigma = \frac{1}{n}\mathbb{A}^T\mathbb{A}$ and $\mathbb{A} = \begin{pmatrix}
    V_{A_1}^T \\ \vdots \\ V_{A_n}^T\end{pmatrix}$
\end{enumerate}
for rates $r_1 = r_1(k)\gtrsim 1, r_2=r_2(n,k)$, and constants $R,C_1,C_2,C_3,C_4>0$ that are independent of $n$ and $k$.
Then,
\begin{align*}
    \E\left( \left\{ V_A^T(\widehat\psi - \psi) \right\}^2 \right) &\lesssim s\frac{r_1(k)\log k}{n} + s\E\left( \max_{j\in\{1,\dots,k\}} \left(\sum_{a\in\mathcal{A}} r_{2,a}(n,k)|V_{a,j}| \widehat\varpi_a \right)^2\right)
\end{align*}
for $\max\{s^4/r_1(k)^2,r_1(k)\}\log k \lesssim n$.
\end{theorem}

\begin{remark}[Error metric]
    The previous theorem weights the error according to the marginal distribution of the treatment variable. However, it is possible to state the same error bound for the error weighted according to the empirical proportions of the treatment, i.e., $$\E\left(\Pb_n\left[ \left\{ V_A^T(\widehat\psi - \psi) \right\}^2 \right]\right)\lesssim s\frac{r_1(k)\log k}{n} + s\E\left( \max_{j\in\{1,\dots,k\}} \left(\sum_{a\in\mathcal{A}} r_{2,a}(n,k)|V_{a,j}| \widehat\varpi_a \right)^2\right).$$ This is shown along the way in the proof of the theorem. We wish to point out that this bound on the in-sample error does no longer require the assumptions $|V_{A,j}V_{A,l}|\leq C_1r_1(k)$ and $\Pb(V_{A,j}^2V_{A,l}^2)\leq C_1r_1(k)$. Moreover, we only need to assume $n\gtrsim r_1(k)\log k$ instead of $n\gtrsim \max\{s^4/r_1(k)^2,r_1(k)\}\log k$.
\end{remark}

\subsection{Approximate Sparsity}

Depending on the application, assuming exact sparsity of the statistical functional $\psi$ might appear unreasonable. Instead, one might resort to approximate sparsity, where we do not have to assume that the majority of the entries of $\psi$ are exactly zero. For approximate sparsity, we assume that $\norm{\psi}_1\leq s$.

\subsubsection{Master Theorem under Approximate Sparsity}\label{section:master-theorem-under-approximate-sparsity}

In this subsection, we present an error guarantee of our general Lasso estimator of Section \ref{section:pseudo-outcome-regression} under approximate sparsity. Therefore, assume the setup of Section \ref{section:pseudo-outcome-regression}.

\begin{theorem}[Error rate under approximate sparsity]\label{thm:master-approximate-sparsity}
Let $\widehat\psi$ be the Lasso estimator defined in (\ref{def:lasso}), minimizing the empirical risk in (\ref{eq:emp-risk}). Suppose approximate sparsity, i.e., $\norm{\psi}_1\leq s$ for some sparsity constraint $s$. Moreover, assume for all $j\in\{1,\dots,k\}$ and $a\in\mathcal{A}$ that
\begin{enumerate}
    \item[(i)] (Almost sure boundedness) $|\widehat\varphi_1(Z)V_{A,j}|$, $\left| \int_\mathcal{A} \widehat\varphi_{2,a}(Z)V_{a,j} d\widehat\Pb(a) \right|, |V_A^T\psi V_{A,j}|\leq C_1r_1(k)$ almost surely,
    \item[(ii)] (Second moment boundedness) $\Pb\left\{ (\widehat\varphi_1(Z)V_{A,j})^2 \right\}$, $\Pb\left\{\left(\int_\mathcal{A} \widehat\varphi_{2,a}(Z)V_{a,j} d\widehat\Pb(a)\right)^2\right\}$,\\ $(\int_\mathcal{A}V_a^T\psi V_{a,j}d\widehat\Pb(a))^2 \leq C_1\widehat r_1(k)$ almost surely where $\widehat r_1(k)$ can depend on $D_1$ and\\ $\E_{D_1}\{\widehat r_1(k)\} \leq Rr_1(k)$,
    \item[(iii)] (Nuisance penalty) $\left| \Pb\left\{ \E(\widehat\varphi_1(Z)\mid X, A=a)p_{A\mid X}(a\mid X) + \widehat\varphi_{2,a}(Z) - V_a^T\psi \right\} \right|\leq C_2r_{2,a}(n,k)$ where $p_{A\mid X}$ is the conditional density of $A$ given $X$ with respect to $\widehat\Pb$,
    \item[(iv)] (Difference in measures) $\E\left[\max_{j\in\{1,\dots,k\}} \left| (\widetilde\Pb - \widehat\Pb)\left\{V_A^T\psi V_{A,j}\right\} \right| \right]\leq C_3\sqrt{\frac{r_1(k)\log k}{n}}$, and
    \item[(v)] (Integrability) $(\mathcal{A}, \widehat\Pb)$ is $\sigma$-finite and the map
    $$
        (a, Z) \mapsto \left( \E(\widehat\varphi_1(Z)\mid X, A=a)p_{A\mid X}(a\mid X) + \widehat\varphi_{2,a}(Z) - V_a^T\psi \right) V_{a,j}
    $$
    is $\widehat\Pb\otimes\Pb$-integrable
\end{enumerate}
for rates $r_1 = r_1(k)\gtrsim 1, r_2=r_2(n,k)$, and constants $R,C_1,C_2,C_3>0$ that are independent of $n$ and $k$.
Then,
\begin{align*}
    \E\left( \int_\mathcal{A}\left\{ V_a^T(\widehat\psi - \psi) \right\}^2 d\widetilde\Pb(a) \right) &\lesssim s\sqrt{\frac{r_1(k)\log k}{n}} + s\E\left( \max_{j\in\{1,\dots,k\}} \int_\mathcal{A} r_{2,a}(n,k)|V_{a,j}| d\widehat\Pb(a) \right)
\end{align*}
for $r_1(k)\log k \lesssim n$.
\end{theorem}

First, completely analogously to the error rate under exact sparsity, the given error rate consists of an oracle rate and a nuisance estimation error rate.\\

Second, for assuming only approximate sparsity, we pay the price of slower convergence rates. The obtained rate under approximate sparsity, containing a term of order $\sqrt{\log (k)/n}$, is also called \textit{slow Lasso rate}. In contrast, the rate obtained under exact sparsity containing a term of the faster order $\log (k)/n$ is called \textit{fast Lasso rate}. This follows a well-known pattern in the high-dimensional regression literature, e.g., refer to \citet{raskutti2011minimax}.\\

Third, although the assumptions are largely identical to those in Theorem \ref{thm:master}, we do not assume restricted eigenvalues here. The reason is that this assumption, in combination with exact sparsity, allowed us to achieve the fast Lasso rate. As we no longer assume exact sparsity, the restricted eigenvalue condition is consequently also no longer necessary.\\

\begin{remark}[Analysis of Lasso vs. best subset selection]
    Note that the above theorem for the approximate sparse case only provides an error guarantee for the Lasso estimator, not for the best subset selection estimator. This can be explained by the fact that the proof relies on bounding $\norm{\widehat\psi - \psi}_1\leq 2s$. However, $\norm{\widehat\psi_{\mathrm{subset}}}_1$ does not have to be bounded by the order $s$: The best subset selection constraints the $L_0$-norm, where the $L_0$-constraint tuning parameter has to be respected by the true $\psi$, i.e., we choose the true sparsity $\norm{\psi}_0$ as tuning parameter, so that the basic Lasso inequality can be applied in the proof, i.e., $\widehat R(\widehat\psi_{\mathrm{subset}})\leq \widehat R(\psi)$. However, since we assume only approximate sparsity, we might have $\norm{\psi}_0=k$. Hence, in this case, the best subset selection estimator does not incorporate a constraint in the least squares regression, and, in particular, $\norm{\widehat\psi_{\mathrm{subset}}}_1$ does not have to be bounded by $s$.\\
\end{remark}

As before in Section \ref{section:pseudo-outcome-regression}, we can formulate a discrete version of the previous theorem when the treatment variable is discrete. For details on how the discrete version of the empirical risk is derived, we refer to Remark \ref{remark:discrete-case} and Appendix \ref{sec:discrete-master-thm}.

\begin{theorem}[Discrete version of Theorem \ref{thm:master-approximate-sparsity}]
\label{thm:master-approximate-sparsity-discrete}
Let $\widehat\psi=\widehat\psi_\text{lasso}$ be the Lasso estimator defined in (\ref{def:lasso}), minimizing the empirical risk in (\ref{eq:emp-risk-discrete}). Suppose approximate sparsity, i.e., $\norm{\psi}_1\leq s$ for some sparsity constraint $s$. Moreover, assume for all $j\in\{1,\dots,k\}$ and $a\in\mathcal{A}$ that
\begin{enumerate}
    \item[(i)] (Almost sure boundedness) $\left| \sum_{a\in\mathcal{A}} \widehat\varphi_{a}(Z) V_{a,j}\widehat\varpi_a \right|, |V_A^T\psi V_{A,j}|\leq C_1r_1(k)$ almost surely,
    \item[(ii)] (Second moment boundedness) $\Pb\left\{\left(\sum_{a\in\mathcal{A}} \widehat\varphi_{a}(Z) V_{a,j}\widehat\varpi_a\right)^2\right\}, (\sum_{a\in\mathcal{A}}V_A^T\psi V_{A,j}\widehat\varpi_a)^2 \leq C_1\widehat r_1(k)$ and $\Pb\left\{(V_A^T\psi V_{A,j})^2\right\}\leq C_1r_1(k)$ almost surely where $\widehat r_1(k)$ can depend on $D_1$ and $\E_{D_1}\{\widehat r_1(k)\} \leq Rr_1(k)$, and
    \item[(iii)] (Nuisance estimation error rate) $\left| \Pb\left\{ \widehat\varphi_{a}(Z) - V_a^T\psi \right\} \right|\leq C_2r_{2,a}(n,k)$
\end{enumerate}
for rates $r_1 = r_1(k)\gtrsim 1, r_2=r_2(n,k)$, and constants $R,C_1,C_2>0$ that are independent of $n$ and $k$.
Then,
\begin{align*}
    \E\left( \Pb_n^2\left[ \left\{ V_A^T(\widehat\psi - \psi) \right\}^2 \right] \right) &\lesssim s\sqrt{\frac{r_1(k)\log k}{n}} + s\E\left( \max_{j\in\{1,\dots,k\}} \sum_{a\in\mathcal{A}} r_{2,a}(n,k)|V_{a,j}| \widehat\varpi_a \right)
\end{align*}
for $r_1(k)\log k \lesssim n$ where $\widehat\varpi_j=\Pb_n^1\{\mathds{1}(A=j)\}$.
\end{theorem}

\begin{remark}[In-sample vs. out-of-sample error]\label{remark:out-of-sample-error}
    It is possible to state the above error bound when the prediction error is weighted by the marginal distribution of the treatments instead of its empirical proportions at the cost of an $s^2$-dependence instead of $s$-dependence in the oracle rate. That is,
    \begin{equation*}
        \E\left( \left\{ V_A^T(\widehat\psi - \psi) \right\}^2 \right) \lesssim s^2\sqrt{\frac{r_1(k)\log k}{n}} + s\E\left( \max_{j\in\{1,\dots,k\}} \sum_{a\in\mathcal{A}} r_{2,a}(n,k)|V_{a,j}| \widehat\varpi_a \right),
    \end{equation*}
    where we need to additionally assume that $|V_{A,j}V_{A,l}|, \Pb\{V_{A,j}^2V_{A,l}^2\}\leq C_1r_1(k)$. This can be verified analogously to the proof of Theorem \ref{thm:master-discrete}, where the error is decomposed into the expected in-sample error and an empirical process term. The expected in-sample error can be bounded as stated by the above theorem, and the empirical process term can be analyzed using the calculation in \eqref{eq:emp-process-term} and the additionally assumed almost-sure and second-moment boundedness.
\end{remark}

\subsubsection{Estimating Mean Potential Outcomes for Single Multi-Valued Treatments under Approximate Sparsity}\label{section:single-treatment-approximate-sparsity}

In this section, we provide an error guarantee of the proposed Lasso estimator of mean potential outcomes for single multi-valued treatments under approximate sparsity, i.e., when assuming $\norm{\psi}_1\leq s$. Therefore, assume the setup of Section \ref{section:single-treatment}. We show that a slow Lasso rate can be achieved and that fast convergence rates are possible even when there are many treatments in proportion to the sample size.\\

For some practical applications, exact sparsity may not appear to be a realistic assumption, as it only allows $s$ mean potential outcomes to differ from some intercept value. At the same time, all other mean potential outcomes must be exactly equal. For that reason, some applications must resort to approximate sparsity. Although this assumption also requires that the mean potential outcomes do not deviate too much from the intercept, it imposes no constraint on the number of mean potential outcomes that may differ from the intercept. For example, all mean potential outcomes would be allowed to differ from the intercept as long as these deviations are small enough.

\begin{theorem}\label{thm:single-treatment-approximate-sparsity}
    Let $\widehat\psi=\widehat\psi_\text{lasso}$ be the Lasso estimator defined in (\ref{def:lasso-single-treatment}) minimizing the empirical risk in (\ref{eq:emp-risk-single-trt}). Assume the model in (\ref{eq:saturated-model}) with approximate sparsity, i.e., $\norm{\psi}_1\leq s$ for some sparsity constraint $s$. Moreover, assume that for all $a\in\{1,\dots,k\}$
    \begin{itemize}
        \item[(i)] (Nearly uniform treatment probabilities) $\frac{c}{k}\leq\varpi_a=\Pb(A=a)\leq \frac{C}{k}$,
        \item[(ii)] (Estimated propensity score close to empirical weights) $|\widehat\varpi_a/\widehat\pi_a(X)| \leq B$ almost surely, and
        \item[(iii)] (Boundedness)
            $\left|\widehat\mu_a(X)\right|, |Y|, |\psi_0| \leq B$ almost surely
    \end{itemize}
    for constants $c,C,B>0$ that are independent of $k$ and $n$. Then,
    \begin{equation*}
        \E\left(\sum_{a=1}^k w(a)(\widehat\psi_a - \psi_a)^2\right) \lesssim s\sqrt{\frac{\log k}{kn}} + \frac{s}{k}\delta_n\epsilon_n
    \end{equation*}
    for $w(a)=\Pb_n^2\{\mathds{1}(A=a)\}$ and $n\gtrsim k\log k$, when $\delta_n(a)\leq\delta_n$ and $\epsilon_n(a)\leq\epsilon_n$ for all $a\in\{1,\dots,k\}$.
\end{theorem}

The price we pay for assuming approximate sparsity rather than exact sparsity is that the rate obtained above is slower. For the oracle part of the rate, we now obtain $\sqrt{\frac{\log k}{kn}}$ instead of $\frac{\log k}{n}$. Similarly, the nuisance penalty now contains $\delta_n\epsilon_n$ instead of $(\delta_n\epsilon_n)^2$. The requirement for the nuisance error to achieve the oracle rate is the same as for exact sparsity, as we need $\delta_n\epsilon_n \lesssim\sqrt{\frac{k\log k}{n}}$, which is, for instance, satisfied when $\delta_n,\epsilon_n\asymp \left( \frac{k\log k}{n} \right)^{1/4}$.\\

For a discussion of the assumptions, we refer to the remarks belonging to Theorem \ref{thm:single-treatment}.\\

\begin{remark} Similar to the corresponding result under exact sparsity in Theorem \ref{thm:single-treatment}, it is possible to derive a more informative upper bound
    \begin{equation*}
        \E\left(\sum_{a=1}^k w(a)(\widehat\psi_a - \psi_a)^2\right) \lesssim s\sqrt{\frac{\log k}{kn}} + s\E\left[ \max_{a\in\{1,\dots,k\}} \left\{\widehat\varpi_a \delta_n(a)\epsilon_n(a) \right\} \right],
    \end{equation*}
    which reveals that the nuisance penalty is weighted by the empirical proportion at each treatment level. This result is shown as part of the above theorem's proof. Due to Remark \ref{remark:out-of-sample-error} on the master theorem under approximate sparsity (which is used to prove this result), it is also possible to state the above upper bound for the out-of-sample error instead of the expected in-sample error at the cost of an additional $s^2$-dependence, i.e.,
    \begin{equation*}
        \E\left(\sum_{a=1}^k \varpi_a(\widehat\psi_a - \psi_a)^2\right) \lesssim s^2\sqrt{\frac{\log k}{kn}} + s\E\left[ \max_{a\in\{1,\dots,k\}} \left\{\widehat\varpi_a \delta_n(a)\epsilon_n(a) \right\} \right] \lesssim s^2\sqrt{\frac{\log k}{kn}} + \frac{s}{k}\delta_n\epsilon_n
    \end{equation*}
    for $\varpi_a=\Pb(A=a)$.
\end{remark}

\subsubsection{Estimating Mean Potential Outcomes for Vector Treatments under Approximate Sparsity}\label{section:vector-treatment-approximate-sparsity}

In this section, we provide an error guarantee of the Lasso estimator for vector treatments under approximate sparsity, i.e., when assuming $\norm{\psi}_1\leq s$. Therefore, assume the setup of Section \ref{section:vector-treatment}. We demonstrate that a slow Lasso rate can be achieved, and that fast convergence rates are possible even when the dimension of the vector treatment is large relative to the sample size.\\

First we assume the setup of Section \ref{sec:binary-vector-treatment}, i.e., when the vector treatment is binary. Then, we can state the following error rate under approximate sparsity.

\begin{theorem}[Binary vector treatments under approximate sparsity]\label{thm:binary-vector-trt-approximate-sparsity}
    Let $\widehat\psi=\widehat\psi_\text{lasso}$ be the Lasso estimator defined in (\ref{def:lasso}) minimizing the empirical risk in (\ref{eq:emp-risk-binary-vector-trt}). Assume the model defined in (\ref{eq:marginal-structural-model}) and suppose approximate sparsity, i.e., $\norm{\psi}_1\leq s$ for some sparsity constraint~$s$. Moreover, assume that for all $a\in\{0,1\}^k$
    \begin{itemize}
        \item[(i)] (Estimated propensity score close to empirical weights) $|\widehat\varpi_a/\widehat\pi_a(X)| \leq B$, and
        \item[(ii)] (Boundedness) $|Y|, |\widehat\mu|, |\psi_0|\leq B$
    \end{itemize}
    almost surely for a constant $B>0$ that is independent of $k$ and $n$.
    Then,
    \begin{equation*}
        \E\left( \Pb_n^2\left[ \left\{ A^T(\widehat\psi - \psi) \right\}^2 \right] \right) \lesssim s\sqrt{\frac{\log k}{n}} + s\delta_n\epsilon_n
    \end{equation*}
    for $\log k \lesssim n$, when $\delta_n(a)\leq\delta_n$ and $\epsilon_n(a)\leq\epsilon_n$ for all $a\in\{0,1\}^k$.
\end{theorem}

The obtained rate consists again of two parts, analogously to Theorem \ref{thm:binary-vector-trt}. The price we pay for only assuming approximate sparsity is that we obtain a slow Lasso rate. For further discussion of the assumptions and implications, we refer to the remarks of Theorem \ref{thm:binary-vector-trt}.\\

\begin{remark}
    Similar to the result under exact sparsity in Theorem \ref{thm:binary-vector-trt}, it is possible to show the more informative upper bound where the nuisance errors are weighted individually at each treatment level, i.e.,
    \begin{equation*}
        \E\left( \Pb_n^2\left[ \left\{ A^T(\widehat\psi - \psi) \right\}^2 \right] \right) \lesssim s\sqrt{\frac{\log k}{n}} + s\E\left( \sum_{a\in \mathcal{A}} \widehat\varpi_a \delta_n(a) \epsilon_n(a) \right).
    \end{equation*}
    This result is proved as part of the above theorem's proof. Moreover, based on Remark \ref{remark:out-of-sample-error}, we can state the above error for the out-of-sample error at the cost of an additional $s^2$-dependence, i.e.,
    \begin{equation*}
        \E\left( \left\{ A^T(\widehat\psi - \psi) \right\}^2 \right) \lesssim s^2\sqrt{\frac{\log k}{n}} + s\E\left( \sum_{a\in \mathcal{A}} \widehat\varpi_a \delta_n(a) \epsilon_n(a) \right) \lesssim s^2\sqrt{\frac{\log k}{n}} + s\delta_n\epsilon_n.
    \end{equation*}\\
\end{remark}

Finally, assume the setup of Section \ref{sec:continuous-vector-treatment}. Then, we can state the following error rate for continuous vector treatments under approximate sparsity.

\begin{theorem}[Continuous vector treatments under approximate sparsity]\label{thm:continuous-vector-trt-approximate-sparsity}
    Let $\widehat\psi=\widehat\psi_\text{lasso}$ be the Lasso estimator defined in (\ref{def:lasso}), minimizing the empirical risk in (\ref{eq:emp-risk-continuous-vector-trt}). Assume the model defined in (\ref{eq:marginal-structural-model}) and suppose approximate sparsity, i.e., $\norm{\psi}_1\leq s$ for some sparsity constraint~$s$. Moreover, assume
    \begin{itemize}
        \item[(i)] (Strong positivity) $\widehat\pi_a(X)\geq\epsilon$ for all $a\in\mathcal{A}$ almost surely, and
        \item[(ii)] (Boundedness) $|Y|, |\widehat\mu|, |\psi_0|, \norm{A}_\infty, \widehat p, \int_\mathcal{A} \widehat p(a)da \leq B$ almost surely
    \end{itemize}
    for constants $B,\epsilon>0$ that are independent of $k$ and $n$.
    Then,
    \begin{equation*}
        \E\left( \int_\mathcal{A}\left\{ a^T(\widehat\psi - \psi) \right\}^2 \widehat p(a)da\right) \lesssim s\sqrt{\frac{\log k}{n}} + s\delta_n\epsilon_n
    \end{equation*}
    for $\log k \lesssim n$, when $\delta_n(a)\leq\delta_n$ and $\epsilon_n(a)\leq\epsilon_n$ for all $a\in\mathcal{A}$.
\end{theorem}

\begin{remark}
    Similar to the result under exact sparsity in Theorem \ref{thm:continuous-vector-trt}, it is possible to show the more informative upper bound where the nuisance errors are weighted according to the estimated treatment density, i.e.,
    $$
        \E\left( \int_\mathcal{A}\left\{ a^T(\widehat\psi - \psi) \right\}^2 \widehat p(a)da\right) \lesssim s\sqrt{\frac{\log k}{n}} + s\E\left( \int_\mathcal{A} \delta_n(a)\epsilon_n(a)\widehat p(a) da \right).
    $$
\end{remark}

\subsection{Verification of the Sample Covariance Condition for Binary Vector Treatments} \label{sec:sample-covariance}

In this subsection, we provide a detailed example when the sample covariance condition of Theorem \ref{thm:binary-vector-trt} for binary vector treatments is satisfied: Assume $A=(A_1,\dots, A_k)\in\{0,1\}^k$ with independent components and $\Pb(A_j=1)=1/2$ for all $j$. We want to verify that the sample covariance condition is satisfied in this special case, i.e., we want to show that 
$$\sup_{v\in\mathcal{V}} \left|v^T\Delta v \right| = \sup_{v\in\mathcal{V}} \left|v^T\left( \frac{1}{n}\sum_{i=1}^n A_iA_i^T - \E(AA^T) \right)v \right| > \frac{D}{2}$$
with probability less than $\alpha\leq\frac{\log k}{n}$ where $\mathcal{V}=\{v\in\mathbb{R}^k \mid \norm{v}_2=1, \norm{v}_1\leq 2\sqrt{s}\}$ and $\Delta = \frac{1}{n}\mathbb{A}^T\mathbb{A} - \E(AA^T)$. It can by verified with a straightforward calculation that $|v^T\Delta v|\leq \norm{\Delta}_\infty\norm{v}_1^2$. For $v\in\mathcal{V}$ we have $\norm{v}_1\leq 2\sqrt{s}$, so $\sup_{v\in\mathcal{V}} |v^T\Delta v|\leq 4s\norm{\Delta}_\infty$. Hence, it suffices to bound $\norm{\Delta}_\infty$. We start by bounding the entries of $\Delta_{ij}$. We have
$$
    \Delta_{lm} = \frac{1}{n}\sum_{i=1}^n \left\{ A_{i,l}A_{i,m} - \E(A_{i,l}A_{i,m}) \right\}
$$
which is a sum of independent mean-zero random variables. Using the fact that
$$
    |A_{i,l}A_{i,m} - \E(A_{i,l}A_{i,m})|\leq 2B^2
$$
and
$$
    \Pb\left\{(A_{i,l}A_{i,m} - \E(A_{i,l}A_{i,m}))^2\right\} \leq 4B^4
$$
under the assumption $\norm{A}_\infty\leq B$, we can apply Bernstein's inequality to obtain
$$
    \Pb(|\Delta_{lm}| \geq t)\leq 2\exp\left(-\frac{nt^2}{8B^4 + \frac{4}{3}B^2t}\right).
$$
By the union bound, we obtain
$$
    \Pb(\norm{\Delta}_\infty \geq t)\leq k^2\max_{l,m} \Pb(|\Delta_{lm}| \geq t) \leq 2k^2\exp\left(-\frac{nt^2}{8B^4 + \frac{4}{3}B^2t}\right).
$$
So,
$$
    \Pb\left( \sup_{v\in\mathcal{V}} \left|v^T\Delta v \right| \geq \frac{D}{2} \right) \leq \Pb\left( \norm{\Delta}_\infty \geq \frac{D}{8s} \right) \leq 2k^2\exp\left(-\frac{n(D/8s)^2}{8B^4 + \frac{D}{6s}B^2}\right).
$$
We need that the right-hand side is bounded by $\frac{\log k}{n}$, i.e., we require
$$
    2k^2\exp\left(-\frac{n(D/8s)^2}{8B^4 + \frac{D}{6s}B^2}\right) \leq \frac{\log k}{n},
$$
which is satisfied for $n\gtrsim s^2\log k$.

\subsection{Cross-Validation}\label{sec:cross-validation}

For the data analysis in Section \ref{section:empirical-data-analysis}, we use cross-validation based on an estimated risk to choose the sparsity constraint. More specifically, we split the data $D_n$ into three samples $D_1$, $D_2$, and $D_3$. On the first sample $D_1$, we estimate the nuisance parameters $\widehat\pi$ and $\widehat\mu$ to obtain the estimated pseudo-outcomes $\widehat\varphi_a(Z)$. On the second sample $D_2$, we minimize the empirical risk in \eqref{eq:emp-risk-single-trt} subject to different sparsity constraints, i.e., we construct $\widehat\psi_s = \argmin_{\beta: \norm{\beta}\leq s} \widehat R(\beta)$ for the choices $s\in\Lambda$. On the third sample $D_3$, we calculate the estimated risk for each $\widehat\psi_s$ and choose the one that yields the minimal risk, i.e., our final estimate is $\widehat\psi = \widehat\psi_{s^*}$ for $s^*=\argmin_{s\in\Lambda} \widehat R(\widehat\psi_s)$, where $\widehat R$ is defined as in \eqref{eq:emp-risk-single-trt} with the outer sample average now being over the sample $D_3$.

\section{Proofs of Upper Bounds}

\subsection{Proof of Theorem \ref{thm:single-treatment}}

We derive the rate using Theorem \ref{thm:master-discrete} since the Lasso and best subset selection estimator in \eqref{def:lasso-single-treatment} and \eqref{def:bestsubset-single-treatment}, minimizing the empirical risk in \eqref{eq:emp-risk-single-trt}, are special cases of the estimators defined in \eqref{def:lasso} and \eqref{def:bestsubset}, minimizing the discrete general risk in \eqref{eq:emp-risk-discrete}. This can be verified by setting
$$
    \widehat\varphi_a(Z)=\sqrt{k}\left(\frac{\mathds{1}(A=a)}{\widehat\pi_a(X)}\{Y-\widehat\mu_a(X)\}+\widehat\mu_a(X)-\psi_0\right),\quad V_a=\sqrt{k}(\mathds{1}(a=1),\dots,\mathds{1}(a=k)),
$$
so that the risk in \eqref{eq:emp-risk-discrete} turns into \eqref{eq:emp-risk-single-trt} up to scaling by $k$. Since $\argmin_\beta \widehat R(\beta) = \argmin_\beta k\widehat R(\beta)$, i.e., the additional scaling does not matter, this shows that we recover the Lasso and best subset selection estimators as special cases.\\

In the following sections, we now verify the necessary assumptions of Theorem \ref{thm:master-discrete} for those specific choices of $\widehat\varphi_a$ and $V_a$.

\subsubsection{Almost Sure Boundedness}\label{subsection:bounded-weighted-pseudo-outcomes}
We have
\begin{equation*}
    \left| \sum_{a=1}^k \widehat\varphi_a(Z)V_{a,j}\widehat\varpi_a \right| \leq k\widehat\varpi_j \left[ \frac{\mathds{1}(A=j)}{\widehat\pi_j(X)}\left| Y - \widehat\mu_j(X) \right| + |\widehat\mu_j(X)| + |\psi_0| \right] \leq k\frac{\widehat\varpi_j}{\widehat\pi_j(X)}\mathds{1}(A=j)2B + k\widehat\varpi_j2B \lesssim k
\end{equation*}
by using Assumption (ii),
\begin{equation*}
    \left| V_A^T\psi V_{A,j} \right| \leq k|\psi_A|\mathds{1}(A=j)\lesssim k
\end{equation*}
by using the outcome boundedness assumptions in (iii), and
$$
    |V_{A,j}V_{A,l}| = k\mathds{1}(A=j)\mathds{1}(A=l) \leq k.
$$
Therefore, Assumption (i) of Theorem \ref{thm:master-discrete} is satisfied with $r_1(k)=k$.

\subsubsection{Second Moment Boundedness}\label{subsection:bounded-weighted-second-moments}
We have
\begin{align*}
    \Pb\left\{ \left( \sum_{a=1}^k\widehat\varphi_a(Z)V_{a,j}\widehat\varpi_a \right)^2 \right\} &\leq \Pb\left\{ \left( k\frac{\widehat\varpi_j}{\widehat\pi_j(X)}\mathds{1}(A=j)2B + k\widehat\varpi_j2B \right)^2 \right\} \leq 2\Pb\left\{ k^2\mathds{1}(A=j)4B^4 + k^2\widehat\varpi_j^24B^2 \right\} \\
    &\leq 8B^4k^2\varpi_j + 8B^2k^2\widehat\varpi_j^2 \lesssim k+k^2\max_{j\in\{1,\dots,k\}}\widehat\varpi_j^2,
\end{align*}
\begin{align*}
    \left(\sum_{a=1}^k V_a^T\psi V_{a,j}\widehat\varpi_a\right)^2 = (k\psi_j\widehat\varpi_j)^2 \lesssim k^2\max_{j\in\{1,\dots,k\}} \widehat\varpi_j^2,
\end{align*}
as well as
\begin{equation*}
    \Pb\left\{ (V_A^T\psi V_{A,j})^2 \right\} = \Pb\left\{ k^2\psi_j^2\mathds{1}(A=j) \right\} = k^2\psi_j\varpi_j \lesssim k
\end{equation*}
and
\begin{equation*}
    \Pb\left(V_{A,j}^2V_{A,l}^2\right) = k^2\Pb\left\{\mathds{1}(A=j)\mathds{1}(A=l)\right\}\leq k^2\E\left\{\mathds{1}(A=j)\right\} = k^2\varpi_j \lesssim k.
\end{equation*}
Therefore, we can set $r_1(k)=k$ and $\widehat r_1(k)=k + k^2\max_{j\in\{1,\dots,k\}} \widehat\varpi_j^2$, as
\begin{equation*}
    \E\left(\max_{j\in\{1,\dots,k\}} \widehat\varpi_j^2\right) \leq \sum_{j=1}^k\E(\widehat\varpi_j^2) = \sum_{j=1}^k \frac{1}{n^2}\E\left\{\mathrm{Bin}(n,\varpi_j)^2\right\} \leq \frac{1}{n^2} \sum_{j=1}^k \left\{ n(n-1)\frac{C^2}{k^2} + n\frac{C}{k} \right\} \lesssim \frac{1}{k}
\end{equation*}
and
\begin{equation*}
    \E\left(\max_{j\in\{1,\dots,k\}} \widehat\varpi_j^4\right) \leq \frac{k}{n^4}\E\left\{ \mathrm{Bin}(n,C/k)^4 \right\} \lesssim \frac{k}{n^4}\cdot\frac{n^4}{k^4} \lesssim \frac{1}{k^2}
\end{equation*}
imply $\E(\widehat r_1(k)^j) \lesssim r_1(k)^j$ for $j=1,2$.

\subsubsection{Nuisance Estimation Error Rate}\label{subsection:nuisance-estimation-error-rate}

We calculate that
\begin{align*}
    \Pb\left\{ \widehat\varphi_a(Z) - V_a^T\psi \right\} &= \sqrt{k}\Pb\left\{ \frac{\mathds{1}(A=a)}{\widehat\pi_a(X)}\{\E(Y\mid X,A) - \widehat\mu_a(X)\} + \widehat\mu_a(X) -\mu_a(X)\right\}\\
    &= \sqrt{k}\Pb\left\{\frac{\mathds{1}(A=a)}{\widehat\pi_a(X)}\{\mu_a(X) - \widehat\mu_a(X)\} + \widehat\mu_a(X)-\mu_a(X)  \right\}\\
    &= \sqrt{k}\Pb\left\{ \frac{\E(\mathds{1}(A=a)\mid X)}{\widehat\pi_a(X)}\{\mu_a(X) - \widehat\mu_a(X)\} + \widehat\mu_a(X) -\mu_a(X) \right\}\\
    &= \sqrt{k}\Pb\left\{ \frac{\pi_a(X)}{\widehat\pi_a(X)}\{\mu_a(X) - \widehat\mu_a(X)\} + \widehat\mu_a(X) - \mu_a(X)\right\} \\
    &= \sqrt{k}\Pb\left\{ \left(\mu_a(X) - \widehat\mu_a(X)\right)\left( \frac{\pi_a(X)}{\widehat\pi_a(X)} - 1 \right) \right\} \\
    &\leq \sqrt{k}\anorm{\frac{\pi_a(\cdot)}{\widehat\pi_a(\cdot)} - 1}\anorm{\mu_a(\cdot) - \widehat\mu_a(\cdot)} = r_{2,a}(n,k).
\end{align*}
Hence, Assumption (iii) of Theorem \ref{thm:master-discrete} is satisfied with $r_{2,a}(n,k) = \sqrt{k}\anorm{\frac{\pi_a(\cdot)}{\widehat\pi_a(\cdot)} - 1}\anorm{\mu_a(\cdot) - \widehat\mu_a(\cdot)}$.

\subsubsection{Restricted Eigenvalue Condition}

We have that
$$
    v^T\E(V_AV_A^T)v =  k\cdot v^T\mathrm{diag}(\E(\mathds{1}(A=j)))v = k\cdot v^T\mathrm{diag}(\Pb(A=j))v = \sum_{j=1}^k k\varpi_j v_j^2 \geq c\norm{v}_2^2
$$
where the last inequality follows from the assumption that $k\varpi_j\geq c$ for all $j\in\{1,\dots,k\}$. Hence, this proves Assumption (iv) of Theorem \ref{thm:master-discrete} for $C_3=c$.

\subsubsection{Sample Covariance Condition}

We have
\begin{equation*}
    \mathbb{A} = \sqrt{k}\begin{pmatrix}
        \mathds{1}(A_1=1) & \dots & \mathds{1}(A_1=k) \\
        \vdots & & \vdots \\
        \mathds{1}(A_n=1) & \dots & \mathds{1}(A_n=k)
    \end{pmatrix},
\end{equation*}
so
\begin{equation*}
    v^T\widehat\Sigma v = v^T\left(\frac{1}{n}\mathbb{A}^T\mathbb{A}\right)v = k\cdot v^T\mathrm{diag}\left(\frac{N_j}{n}\right)v
\end{equation*}
for $N_j=\sum_{i=1}^n \mathds{1}(A_i=j)$. Further, we have
\begin{equation*}
    v^T\Sigma v = k\cdot v^T\mathrm{diag}(\Pb(A=j))v.
\end{equation*}
Hence, we must show that
\begin{equation*}
    \sup_{v\in S_1} \left| v^T\left\{\mathrm{diag}\left(\frac{N_j}{n}\right) - \mathrm{diag}(\Pb(A=j))\right\}v \right| \leq \frac{c}{2k}
\end{equation*}
with high probability. Since the matrix $\mathrm{diag}\left(\frac{N_j}{n}\right) - \mathrm{diag}(\Pb(A=j))$ is symmetric, the left-hand side equals the operator 2-norm of this matrix, hence it suffices to bound the operator norm by $c/2k$. Since we have a diagonal matrix, it further suffices to bound the largest diagonal entry with high probability. Consequently, we need to show that
\begin{equation}\label{eq:sample-variance-bound}
    \max_{j\in\{1,\dots,k\}} \left| N_j - n\Pb(A=j)\right| \leq \frac{cn}{2k}
\end{equation}
with probability at least $1-\alpha$ for $\alpha\leq \frac{\log k}{n}$. For a fixed $j$, we can use the assumption $\pi_j\leq C/k$ to obtain
$$
    \Pb\left(\left| N_j - n\Pb(A=j) \right| \geq \frac{cn}{2k}\right) \leq \Pb\left( |N_j - n\Pb(A=j)| \geq \frac{c}{2C}n\Pb(A=j) \right).
$$
We have that $N_j$ is a sum of iid Bernoulli random variables with mean $\mu=n\Pb(A=j)$. Hence, we can apply the multiplicative Chernoff bound to obtain, for $\delta := \frac{c}{2C}$,
$$
    \Pb\left( |N_j - n\Pb(A=j)| \geq \frac{c}{2C}n\Pb(A=j) \right) = \Pb\left( |N_j - \mu| \geq \delta\mu \right) \leq 2\exp\left( \frac{-\delta^2\mu}{3} \right) = 2\exp\left( -\frac{c^2}{12C^2}n\varpi_j \right).
$$
Now using $\varpi_j\geq c/k$, we obtain
$$
    \Pb\left( |N_j - n\Pb(A=j)| \geq \frac{c}{2C}n\Pb(A=j) \right) \leq 2\exp\left( -\frac{c^3}{12C^2}\frac{n}{k} \right).
$$
Using the union bound, we immediately obtain
$$
    \Pb\left(\exists j\in\{1,\dots,k\}: |N_j - n\Pb(A=j)| \geq \frac{c}{2C}n\Pb(A=j) \right) \leq 2k\exp\left( -\frac{c^3}{12C^2}\frac{n}{k} \right).
$$
To satisfy (\ref{eq:sample-variance-bound}) with probability at least $1-\alpha$, we need to bound the previous line by $\alpha$, i.e., we need
$$
    2k\exp\left( -\frac{c^3}{12C^2}\frac{n}{k} \right) \leq \alpha.
$$
Setting $\alpha=\frac{\log k}{n}$, and substituting $n=\gamma k\log k$, the previous condition is equivalent to
$$
    k^{2-\frac{c^3\gamma}{12C^2}}\leq \frac{1}{2\gamma}.
$$
Taking the log on both sides, we obtain
$$
    \left( 2-\frac{c^3\gamma}{12C^2} \right) \log k \leq -\log(2\gamma),
$$
which is equivalent to
$$
    \gamma \geq \frac{12C^2}{c^3}\left( 2 + \frac{\log(2\gamma)}{\log k}\right).
$$
Using that $k\geq 2$, the above is satisfied whenever we have
$$
    \gamma \geq \frac{12C^2}{c^3}\left( 2 + \frac{\log(2\gamma)}{\log 2}\right).
$$
This condition holds for a large enough constant $\gamma$.

\subsubsection{Final Rate}\label{subsection:final-rate}

Plugging in the obtained rates into the error guarantee of Theorem \ref{thm:master-discrete}, we obtain
\begin{align*}
    \E\left( \left\{ V_A^T(\widehat\psi - \psi) \right\}^2 \right) &\lesssim s\frac{k\log k}{n} + sk^2\E\left( \max_{a\in\{1,\dots,k\}} \widehat\varpi_a^2\anorm{\frac{\pi_a(\cdot)}{\widehat\pi_a(\cdot)} - 1}^2\anorm{\mu_a(\cdot) - \widehat\mu_a(\cdot)}^2 \right).
\end{align*}
The left-hand side equals
\begin{equation*}
    \E\left( k\sum_{a=1}^k \mathds{1}(A=a)(\widehat\psi_a - \psi_a)^2 \right) = k\cdot\E\left( \sum_{a=1}^k \varpi_a(\widehat\psi_a - \psi_a)^2 \right),
\end{equation*}
so after scaling by $1/k$ we obtain
\begin{align*}
    \E\left( \sum_{a=1}^k \varpi_a(\widehat\psi_a - \psi_a)^2 \right) \lesssim s\frac{\log k}{n} + sk\E\left( \max_{a\in\{1,\dots,k\}} \widehat\varpi_a^2\anorm{\frac{\pi_a(\cdot)}{\widehat\pi_a(\cdot)} - 1}^2\anorm{\mu_a(\cdot) - \widehat\mu_a(\cdot)}^2 \right).
\end{align*}
It remains to show that
$$
    sk\E\left( \max_{a\in\{1,\dots,k\}} \widehat\varpi_a^2\anorm{\frac{\pi_a(\cdot)}{\widehat\pi_a(\cdot)} - 1}^2\anorm{\mu_a(\cdot) - \widehat\mu_a(\cdot)}^2 \right) \leq \frac{s}{k}(\delta_n\epsilon_n)^2.
$$
Using that $\delta_n(a)\epsilon_n(a)\leq \delta_n\epsilon_n$, we obtain
\begin{equation}\label{eq:simplify-nuisance-error}
    sk\E\left[ \max_{a\in\{1,\dots,k\}} \left\{\widehat\varpi_a \delta_n(a)\epsilon_n(a) \right\}^2 \right] \leq s(\delta_n\epsilon_n)^2k\E\left( \max_{a\in\{1,\dots,k\}} \widehat\varpi_a^2 \right).
\end{equation}
We now compute that
$$
    \E\left( \max_{a\in\{1,\dots,k\}} \widehat\varpi_a^2 \right) = \frac{1}{n^2}\E\left( \max_{a\in\{1,\dots,k\}} N_a^2 \right)
$$
where $N_a\sim\mathrm{Bin}(n,\varpi_a)$. To compute the expectation on the right-hand side, we note that this is the balls-into-bins problem where we want to calculate the expected squared maximum number of balls across all bins. If we can show that $\E\left( \max_{a\in\{1,\dots,k\}} N_a^2 \right)\lesssim \frac{n^2}{k^2}$, we overall obtain the desired identity. We first calculate that
$$
    \E\left( \max_{a\in\{1,\dots,k\}} N_a^2 \right) = \int_{t>0} \Pb\left(\max_{a\in\{1,\dots,k\}} N_a^2 > t\right)dt = 2\int_{s>0} s\Pb\left(\max_{a\in\{1,\dots,k\}} N_a > s\right)ds
$$
where we substituted $s=\sqrt{t}$. Now define $T=\alpha\frac{n}{k}$ for some constant $\alpha$ that is independent of $n$ and $k$. Then, we can split up the integral in two terms
$$
    \int_{s>0} s\Pb\left(\max_{a\in\{1,\dots,k\}} N_a > s\right)ds = \int_0^T s\Pb\left(\max_{a\in\{1,\dots,k\}} N_a > s\right)ds + \int_T^\infty s\Pb\left(\max_{a\in\{1,\dots,k\}} N_a > s\right)ds.
$$
We first upper bound the first integral as follows:
$$
    \int_0^T s\Pb\left(\max_{a\in\{1,\dots,k\}} N_a > s\right)ds\leq \int_0^T sds = \frac{1}{2}T^2 = \frac{\alpha^2}{2}\frac{n^2}{k^2}.
$$
This upper bound is of the desired order. Hence, it remains to show that the tail integral can be bounded by an equal or smaller order. By using the union bound, we obtain
$$
    \Pb\left(\max_{a\in\{1,\dots,k\}} N_a > s\right) \leq \sum_{a=1}^k \Pb\left( N_a > s \right) \leq k \Pb(N>s)
$$
for some random variable $N\sim \mathrm{Bin}(n,C/k)$. Hence, we obtain
$$
    \int_T^\infty s\Pb\left(\max_{a\in\{1,\dots,k\}} N_a > s\right)ds \leq k\int_T^\infty s\Pb(N>s)ds.
$$
We want to apply the multiplicative Chernoff bound $\Pb(N>(1+\delta)\mu)\leq \left( \frac{\exp(\delta)}{(1+\delta)^{1+\delta}} \right)^\mu$ for a sum $N$ of independent Bernoulli random variables with mean $\E(N)=\mu$. Set $\delta:=\frac{s}{\mu} - 1$. Take $\alpha>C_2$, then $s>T=\alpha\frac{n}{k} > C_2\frac{n}{k}=\E(N)=\mu$, so $\delta>0$. Now,
$$
    \Pb(N>s)=\Pb(N>(1+\delta)\mu)\leq\left( \frac{\exp(\delta)}{(1+\delta)^{1+\delta}} \right)^\mu = \frac{\exp(s-\mu)}{(s/\mu)^s} = \exp\left( s-\mu -s\log\left(\frac{s}{\mu}\right) \right).
$$
Since $\mu\geq 0$ and $\log e=1$, we can upper bound the previous line by
$$
    \leq \exp\left( -s\log\left( \frac{s}{\mu e} \right) \right).
$$
Take $\alpha\geq e^2$, then $\log\left( \frac{s}{\mu e} \right)\geq \log\left( \frac{T}{\mu e} \right) = \log\left( \frac{\alpha\mu}{\mu e} \right)\geq 1$, which implies that
$$
    \Pb(N>s)\leq \exp(-s).
$$
Consequently,
$$
    \int_T^\infty s\Pb\left(\max_{a\in\{1,\dots,k\}} N_a > s\right)ds \leq k\int_T^\infty s\exp(-s)ds.
$$
This right-hand side equals
$$
    k\left[ -(s+1)\exp(-s) \right]_T^\infty = k(T+1)\exp(-T) = k\left(\alpha\mu + 1\right)\exp\left( -\alpha\mu \right).
$$
Finally, we show that this term is of equal or smaller order than $\mu^2$ by calculating that
$$
    \mu=C\frac{n}{k} \geq C\gamma\log k,
$$
so
$$
    k\left(\alpha\mu + 1\right)\exp\left( -\alpha\mu \right) \leq k\left(\alpha\mu + 1\right)\exp\left( -\alpha C\gamma\log k \right) = \left(\alpha\mu + 1\right)k^{1-\alpha C\gamma} \leq 2\alpha\mu k^{1-\alpha C\gamma}.
$$
The right-hand side is of order at most $\mu^2$ if $1-\alpha C\gamma <0$, i.e., we choose $\alpha>\frac{1}{C_2\gamma}$. Hence, the tail integral is negligible. In summary, this proves (\ref{eq:simplify-nuisance-error}) and we can conclude that
$$
    \E\left( \sum_{a=1}^k \varpi_a(\widehat\psi_a - \psi_a)^2 \right) \lesssim s\frac{\log k}{n} + \frac{s}{k}(\delta_n\epsilon_n)^2.
$$

\subsection{Proof of Theorem \ref{thm:master}}

The vast majority of the proof is identical for the Lasso estimator $\widehat\psi_\text{lasso}$ and the best subset selection estimator $\widehat\psi_\text{subset}$. In the following, $\widehat\psi$ can be replaced by either of the estimators, and the statements hold. In the parts of the proof where it matters which estimator is being used, we will explicitly highlight the difference in the analysis.

\subsubsection{Basic Lasso Inequality}

First, by applying the basic Lasso inequality, we obtain
\begin{align*}
    &\E\left( \int_\mathcal{A} \left\{ V_a^T(\widehat\psi - \psi) \right\}^2 d\widetilde\Pb(a) \right) \\
    &= \E\left[ \int_\mathcal{A} \left\{ (V_a^T\widehat\psi)^2 - 2(V_a^T\widehat\psi)(V_a^T\psi) + (V_a^T\psi)^2 \right\} d\widetilde\Pb(a) \right] \\
    &= \E\bigg[ -2\Pb^2_{n}\left\{ \widehat\varphi_1(Z)V_A^T\widehat\psi + \int_\mathcal{A}\widehat\varphi_{2,a}(Z)V_a^T\widehat\psi d\widehat\Pb(a) \right\} + \int_\mathcal{A}(V_a^T\widehat\psi)^2d\widetilde\Pb(a) \\
    &\quad +2\Pb^2_{n}\left\{ \widehat\varphi_1(Z)V_A^T\widehat\psi +\int_\mathcal{A}\widehat\varphi_{2,a}(Z)V_a^T\widehat\psi d\widehat\Pb(a) \right\} + \int_\mathcal{A} \left\{ - 2(V_a^T\widehat\psi)(V_a^T\psi) + (V_a^T\psi)^2 \right\} d\widetilde\Pb(a) \bigg] \\
    &\leq \E\bigg[ -2\Pb^2_{n}\left\{ \widehat\varphi_1(Z)V_A^T\psi + \int_\mathcal{A}\widehat\varphi_{2,a}(Z)V_a^T\psi d\widehat\Pb(a) \right\} + \int_\mathcal{A}(V_a^T\psi)^2d\widetilde\Pb(a) \\
    &\quad +2\Pb^2_{n}\left\{ \widehat\varphi_1(Z)V_A^T\widehat\psi +\int_\mathcal{A}\widehat\varphi_{2,a}(Z)V_a^T\widehat\psi d\widehat\Pb(a) \right\} + \int_\mathcal{A} \left\{ - 2(V_a^T\widehat\psi)(V_a^T\psi) + (V_a^T\psi)^2 \right\} d\widetilde\Pb(a) \bigg].
\end{align*}
Further simplifying the right-hand side gives
\begin{align*}
    &= \E\left[ 2\Pb^2_{n}\left\{ \widehat\varphi_1(Z)V_A^T(\widehat\psi-\psi)+\int_\mathcal{A}\widehat\varphi_{2,a}(Z)V_a^T(\widehat\psi - \psi) d\widehat\Pb(a)\right\} + 2\int_\mathcal{A}(V_a^T\psi)V_a^T(\psi - \widehat\psi)d\widetilde\Pb(a) \right] \\
    &= \E\left[ (\widetilde\Pb - \widehat\Pb)\left\{2V_A^T\psi V_A^T(\psi - \widehat\psi)\right\} + 2\Pb^2_{n}\left\{\widehat\varphi_1(Z)V_A^T(\widehat\psi-\psi) + \int_\mathcal{A} (\widehat\varphi_{2,a}(Z)-V_a^T\psi)V_a^T(\widehat\psi - \psi)d\widehat\Pb(a) \right\} \right].
\end{align*}
By writing the scalar products as sums and taking the maximum over all treatments $j$, we obtain
\begin{equation}\label{upper-bound}
    \leq 2\E\left[ \norm{\widehat\psi - \psi}_1 \max_{j\in\{1,\dots,k\}} \left| \Pb^2_{n}\left\{\widehat\varphi_1(Z)V_{A,j} + \int_\mathcal{A} (\widehat\varphi_{2,a}(Z)-V_a^T\psi)V_{a,j}d\widehat\Pb(a)\right\} + (\widetilde\Pb - \widehat\Pb)\left\{V_A^T\psi V_{A,j}\right\} \right| \right].
\end{equation}

\subsubsection{Restricted Eigenvalue Assumption}

By the restricted eigenvalue assumption, we have
\begin{equation}\label{eq:sample-re-condition}
    v^T\left(\int_\mathcal{A} V_aV_a^Td\widetilde\Pb(a)\right)v \geq C_3\norm{v}_2^2
\end{equation}
with probability at least $1-\alpha$ for all $v$ such that $\norm{v}_1 \leq 2\sqrt{s}\norm{v}_2$. We now leverage the exact sparsity assumption to show that the $L_1$-norm of $\widehat\psi-\psi$ can be accordingly bounded in terms of the $L_2$-norm, i.e., $\norm{\widehat\psi - \psi}_{1}\leq 2\sqrt{s}\norm{\widehat\psi - \psi}_{2}$. If we consider the best subset selection estimator, we know that $\norm{\widehat\psi_\text{subset}}_0\leq s$, and therefore $\norm{\widehat\psi_\text{subset} - \psi}_{1}\leq 2\sqrt{s}\norm{\widehat\psi_\text{subset} - \psi}_{2}$. The same inequality can be shown for the Lasso estimator $\widehat\psi_\text{lasso}$ with a bit of extra work: The fact that $\norm{\psi}_0\leq s$ and $\norm{\widehat\psi_\text{lasso}}_{1}\leq \norm{\psi}_{1}$ implies that $\norm{\widehat\psi_\text{lasso} - \psi}_{1, S^c}\leq \norm{\widehat\psi_\text{lasso} - \psi}_{1, S}$ where $S$ is the support of $\psi$; for details, we refer to the proof of Theorem 7.8 of \cite{wainwright2019highdimensional}. Using this property that the $L_1$-error off the support is dominated by the $L_1$-error on the support (analogously to its use in Theorem 7.13 of \cite{wainwright2019highdimensional}), we obtain
\begin{equation*}
    \norm{\widehat\psi_\text{lasso} - \psi}_{1}=\norm{\widehat\psi_\text{lasso} - \psi}_{1,S}+\norm{\widehat\psi_\text{lasso} - \psi}_{1,S^c}\leq 2\norm{\widehat\psi_\text{lasso} - \psi}_{1,S} \leq 2\sqrt{s}\norm{\widehat\psi_\text{lasso} - \psi}_{2}.
\end{equation*}
Consequently, we proved that
\begin{equation}\label{eq-bound-l1-by-l2}
    \norm{\widehat\psi - \psi}_{1}\leq 2\sqrt{s}\norm{\widehat\psi - \psi}_{2}
\end{equation}
irrespective of which of the two estimators is used. We can now proceed with the same analysis for both estimators again and use the restricted eigenvalue assumption established in (\ref{eq:sample-re-condition}) for $v=\widehat\psi - \psi$ to obtain
\begin{equation*}
    \norm{\widehat\psi - \psi}_2^2 \leq \frac{1}{C_3}\int_\mathcal{A} \left\{ V_a^T(\widehat\psi - \psi) \right\}^2 d\widetilde\Pb(a)
\end{equation*}
with probability at least $1-\alpha$. Then, we can again apply the basic Lasso inequality to the right-hand side, so
\begin{equation}\label{eq-bound-l2}
\begin{split}
    \norm{\widehat\psi - \psi}_{2}^2 \leq \frac{2}{C_3}\norm{\widehat\psi - \psi}_1 \max_{j\in\{1,\dots,k\}} \left| \Pb^2_{n}\left\{\widehat\varphi_1(Z)V_{A,j} + \int_\mathcal{A} (\widehat\varphi_{2,a}(Z)-V_a^T\psi)V_{a,j}d\widehat\Pb(a)\right\} + (\widetilde\Pb - \widehat\Pb)\left\{V_A^T\psi V_{A,j}\right\} \right|,
\end{split}
\end{equation}
and by combining (\ref{eq-bound-l1-by-l2}) and (\ref{eq-bound-l2}),
\begin{equation}\label{eq:bound-on-l1-norm}
    \norm{\widehat\psi - \psi}_1 \leq \frac{8s}{C_3}\cdot\max_{j\in\{1,\dots,k\}} \left| \Pb^2_{n}\left\{\widehat\varphi_1(Z)V_{A,j} + \int_\mathcal{A} (\widehat\varphi_{2,a}(Z)-V_a^T\psi)V_{a,j}d\widehat\Pb(a)\right\} + (\widetilde\Pb - \widehat\Pb)\left\{V_A^T\psi V_{A,j}\right\} \right|
\end{equation}
with probability at least $1-\alpha$. Now, denote by $\mathrm{SC}$ the high-probability event that the restricted eigenvalue condition holds. Conditional on this event, we showed that the previous inequality holds. The next step is to decompose the expectation in (\ref{upper-bound}) into
\begin{align*}
    &2\E\left[ \mathds{1}_{\mathrm{SC}}\cdot\norm{\widehat\psi - \psi}_1 \max_{j\in\{1,\dots,k\}} \left| \Pb^2_{n}\left\{\widehat\varphi_1(Z)V_{A,j} + \int_\mathcal{A} (\widehat\varphi_{2,a}(Z)-V_a^T\psi)V_{a,j}d\widehat\Pb(a)\right\} + (\widetilde\Pb - \widehat\Pb)\left\{V_A^T\psi V_{A,j}\right\} \right| \right]\\
    &+ 2\E\left[ \mathds{1}_{\mathrm{SC}^c}\cdot\norm{\widehat\psi - \psi}_1 \max_{j\in\{1,\dots,k\}} \left| \Pb^2_{n}\left\{\widehat\varphi_1(Z)V_{A,j} + \int_\mathcal{A} (\widehat\varphi_{2,a}(Z)-V_a^T\psi)V_{a,j}d\widehat\Pb(a)\right\} + (\widetilde\Pb - \widehat\Pb)\left\{V_A^T\psi V_{A,j}\right\} \right| \right].
\end{align*}
We show that the second expectation outside the event $\mathrm{SC}$ is negligible, i.e., that it suffices to bound the expectation on the event $\mathrm{SC}$. Outside the event $\mathrm{SC}$, we can bound
\begin{align*}
    &2\E\left[ \mathds{1}_{\mathrm{SC}^c}\cdot\norm{\widehat\psi - \psi}_1 \max_{j\in\{1,\dots,k\}} \left| \Pb^2_{n}\left\{\widehat\varphi_1(Z)V_{A,j} + \int_\mathcal{A} (\widehat\varphi_{2,a}(Z)-V_a^T\psi)V_{a,j}d\widehat\Pb(a)\right\} + (\widetilde\Pb - \widehat\Pb)\left\{V_A^T\psi V_{A,j}\right\} \right| \right] \\
    &\lesssim \E \left( \mathds{1}_{\mathrm{SC}^c}\cdot s\cdot r_1(k) \right) = sr_1(k)\Pb(\mathrm{SC}^c) \leq sr_1(k)\alpha \lesssim s\frac{r_1(k)\log k}{n}
\end{align*}
by using the almost sure boundedness assumption (i) to obtain the first inequality, and the assumption on $\alpha$ to obtain the last inequality. The remainder of the proof will show that this does not contribute to the overall error up to constants. It remains to bound the expectation on the event $\mathrm{SC}$: We use (\ref{eq:bound-on-l1-norm}) to obtain
\begin{align*}
    &2\E\left[ \mathds{1}_{\mathrm{SC}} \cdot \norm{\widehat\psi - \psi}_1 \max_{j\in\{1,\dots,k\}} \left| \Pb^2_{n}\left\{\widehat\varphi_1(Z)V_{A,j} + \int_\mathcal{A} (\widehat\varphi_{2,a}(Z)-V_a^T\psi)V_{a,j}d\widehat\Pb(a)\right\} + (\widetilde\Pb - \widehat\Pb)\left\{V_A^T\psi V_{A,j}\right\} \right| \right]\\
    &\leq \frac{16s}{C_3}\E\left( \max_{j\in\{1,\dots,k\}} \left[\Pb^2_{n}\left\{\widehat\varphi_1(Z)V_{A,j} + \int_\mathcal{A} (\widehat\varphi_{2,a}(Z)-V_a^T\psi)V_{a,j}d\widehat\Pb(a)\right\} + (\widetilde\Pb - \widehat\Pb)\left\{V_A^T\psi V_{A,j}\right\} \right]^2 \right).
\end{align*}

\subsubsection{Decomposition of the Error}

Using the fact that $(a+b)^2\leq 2a^2+2b^2$, we decompose the previous line into
\begin{align*}
    &\frac{64s}{C_3}\underbrace{\E\left[ \max_{j\in\{1,\dots,k\}} (\Pb^2_{n}-\Pb)\left\{ \widehat\varphi_1(Z)V_{A,j} + \int_\mathcal{A} (\widehat\varphi_{2,a}(Z)-V_a^T\psi)V_{a,j}d\widehat\Pb(a) \right\}^2 \right]}_{T_1} \\
    &+ \frac{64s}{C_3}\underbrace{\E\left[ \max_{j\in\{1,\dots,k\}} \Pb\left\{ \widehat\varphi_1(Z)V_{A,j} + \int_\mathcal{A} (\widehat\varphi_{2,a}(Z)-V_a^T\psi)V_{a,j}d\widehat\Pb(a)\right\}^2 \right]}_{T_2} \\
    &+ \frac{32s}{C_3}\underbrace{\E\left[\max_{j\in\{1,\dots,k\}}(\widetilde\Pb - \widehat\Pb)\left\{V_A^T\psi V_{A,j}\right\}^2\right]}_{T_3}.
\end{align*}
Using Assumption (v), the term $T_3$ can immediately be upper bounded by the desired rate. Hence, it remains to bound the terms $T_1$ and $T_2$, which is done in the following subsections.

\subsubsection{Bounding the Term $T_1$}\label{sec-t1}

We recall a useful lemma given in \cite{van1996weak} (refer to Lemma 2.2.10), which will be the main ingredient in obtaining an upper bound of the desired order.
\begin{lemma}\label{lemma-orlicz}
    Let $X_1,\dots, X_k$ be arbitrary random variables that satisfy the tail bound
    \begin{equation}\label{tailbound}
        \Pb(|X_i|\geq x)\leq 2\exp\left(-\frac{1}{2}\frac{x^2}{b+ax}\right)
    \end{equation}
    for all $x$ and fixed $a,b>0$. Then,
    \begin{equation*}
        \anorm{\max_{1\leq i\leq k} X_i}_{\Psi_1} \leq K\left( a\log(1+k) + \sqrt{b}\sqrt{\log(1+k)} \right)
    \end{equation*}
    for a universal constant $K$.
\end{lemma}

It is well-known (e.g., refer to Chapter 2.2 of \cite{van1996weak}) that
\begin{equation*}
    \norm{X}_p \leq p!\norm{X}_{\Psi_1},
\end{equation*}
i.e., that an upper bound on the exponential Orlicz norm will also provide an upper bound on the $L_p$-norm up to a constant. In particular, this allows us to reformulate the previous lemma in terms of the $L_2$-norm.
\begin{lemma}\label{lemma-l2}
    Under the notation and assumptions of Lemma \ref{lemma-orlicz}, we have
    \begin{equation*}
        \E\left(\max_{1\leq i\leq k} |X_i|^2\right) \leq 4K^2\left( a\log(1+k) + \sqrt{b}\sqrt{\log(1+k)} \right)^2.
    \end{equation*}
\end{lemma}

These results will be used to bound the empirical process term $T_1$ in our decomposition. It remains to justify the assumptions of Lemma \ref{lemma-l2}. In particular, we need to show that the random variables $(\Pb^2_{n}-\Pb)\left\{ \int_\mathcal{A} (\widehat\varphi_a(Z)-V_a^T\psi)V_{a,j}d\Pb^1_{n}(a) \right\}, j=1,\dots,k$ satisfy the tail bound in (\ref{tailbound}). This can be done using Bernstein's inequality, as the random variables are independent (conditionally on $D_1$) and mean-zero. By Assumption (i) and (ii), we have
\begin{equation*}
    \left| \widehat\varphi_1(Z)V_{A,j} + \int_\mathcal{A} (\widehat\varphi_{2,a}(Z)-V_a^T\psi)V_{a,j}d\widehat\Pb(a) \right|\leq 3C_1r_1(k)
\end{equation*}
and
\begin{equation*}
    \Pb\left\{ \left( \widehat\varphi_1(Z)V_{A,j} + \int_\mathcal{A} (\widehat\varphi_{2,a}(Z)-V_a^T\psi)V_{a,j}d\widehat\Pb(a) \right)^2 \right\} \leq 10C_1 \widehat r_1(k).
\end{equation*}
By using this almost sure and $L_2$ upper bounds (and scaling them by $1/n$), Bernstein inequality gives
\begin{equation*}
    \Pb\left(\left|(\Pb^2_{n}-\Pb)\left\{ \widehat\varphi_1(Z)V_{A,j} + \int_\mathcal{A} (\widehat\varphi_{2,a}(Z)-V_a^T\psi)V_{a,j}d\widehat\Pb(a)\right\}\right|\geq x\right)\leq 2\exp\left(-\frac{1}{2}\frac{x^2}{b+ax}\right)
\end{equation*}
with $a\asymp \max\{r_1(k), \widehat r_1(k)\}/n$ and $b\asymp \max\{r_1(k), \widehat r_1(k)\}/n$. Then, Lemma \ref{lemma-l2} provides us with the error rate of the desired order by plugging in $a$ and $b$, i.e.,
\begin{align*}
    T_1 &= \E\left(\max_{j\in\{1,\dots,k\}} (\Pb^2_{n}-\Pb)\left\{ \widehat\varphi_1(Z)V_{A,j} + \int_\mathcal{A} (\widehat\varphi_{2,a}(Z)-V_a^T\psi)V_{a,j}d\widehat\Pb(a) \right\}^2\right) \\
    &= \E\left\{ \E\left(\max_{j\in\{1,\dots,k\}} (\Pb^2_{n}-\Pb)\left\{ \widehat\varphi_1(Z)V_{A,j} + \int_\mathcal{A} (\widehat\varphi_{2,a}(Z)-V_a^T\psi)V_{a,j}d\widehat\Pb(a) \right\}^2 \vertline D_1 \right) \right\}\\
    &\lesssim \E\left\{ \frac{\max\{r_1(k), \widehat r_1(k)\}^2\log^2 k}{n^2} + \frac{\max\{r_1(k), \widehat r_1(k)\}\log k}{n} \right\}.
\end{align*}
Since $\E\left[ \max\{r_1(k), \widehat r_1(k)\} \right] \leq r_1(k) + \E\{\widehat r_1(k)\} \lesssim r_1(k)$ and $\E\left[ \max\{r_1(k), \widehat r_1(k)\}^2 \right] \leq r_1(k)^2 + \E\{\widehat r_1(k)^2\} \lesssim r_1(k)^2$ according to Assumption (ii), we obtain
$$
    T_1 \lesssim \frac{r_1(k)\log k}{n} + \frac{r_1(k)^2\log^2k}{n^2} \lesssim \frac{r_1(k)\log k}{n}
$$
using $r_1(k)\log k\lesssim n$.

\subsubsection{Bounding the Term $T_2$}\label{sec-t2}
We have
\begin{align*}
    T_2 &= \E\left[ \max_{j\in\{1,\dots,k\}} \Pb\left\{ \widehat\varphi_1(Z)V_{A,j} + \int_\mathcal{A} (\widehat\varphi_{2,a}(Z)-V_a^T\psi)V_{a,j}d\widehat\Pb(a)\right\}^2 \right] \\
    &= \E\left[ \max_{j\in\{1,\dots,k\}} \Pb\left\{ \int_\mathcal{A}\E(\widehat\varphi_1(Z)V_{A,j}\mid X, A=a)p_{A\mid X}(a\mid X)d\widehat\Pb(a) + \int_\mathcal{A} (\widehat\varphi_{2,a}(Z)-V_a^T\psi)V_{a,j}d\widehat\Pb(a)\right\}^2 \right] \\
    &= \E\left[ \max_{j\in\{1,\dots,k\}} \Pb\left\{ \int_\mathcal{A}\left( \E(\widehat\varphi_1(Z)\mid X, A=a)p_{A\mid X}(a\mid X) + \widehat\varphi_{2,a}(Z)-V_a^T\psi\right)V_{a,j} d\widehat\Pb(a)\right\}^2 \right] \\
    &= \E\left[ \max_{j\in\{1,\dots,k\}} \left(\int_\mathcal{A}\Pb\left\{ \E(\widehat\varphi_1(Z)\mid X, A=a)p_{A\mid X}(a\mid X) + \widehat\varphi_{2,a}(Z)-V_a^T\psi\right\}V_{a,j} d\widehat\Pb(a)\right)^2 \right] \\
    &\leq \E\left[ \max_{j\in\{1,\dots,k\}} \left(\int_\mathcal{A}r_{2,a}(n,k)|V_{a,j}| d\widehat\Pb(a)\right)^2 \right]
\end{align*}
using Assumption (vi) to exchange the integrals in the last equality and Assumption (iii) in the last inequality. This is the desired second part of the rate stated in the theorem. Combining all derived rates concludes the proof.

\subsection{Proof of Theorem \ref{thm:binary-vector-trt}}

We prove the statement by applying the general error rate given in Theorem \ref{thm:master-discrete} and verifying the necessary assumptions therein for
$$
    \widehat\varphi_a(Z) = \frac{\mathds{1}(A=a)}{\widehat\pi_a(X)}\left\{Y-\widehat\mu_a(X)\right\} + \widehat\mu_a(X) - \psi_0 \quad\text{and}\quad V_a=a.
$$

\subsubsection{Almost Sure Boundedness}

We have
\begin{equation*}
    \left| \sum_{a\in\mathcal{A}} \widehat\varphi_a(Z)V_{a,j}\widehat\varpi_a \right| \leq \sum_{a\in\mathcal{A}} \left( \frac{\mathds{1}(A=a)}{\widehat\pi_a(X)}2B + 2B \right)\widehat\varpi_a \leq 2B^2\sum_{a\in\mathcal{A}} \mathds{1}(A=a) + 2B\sum_{a\in\mathcal{A}} \widehat\varpi_a \lesssim 1
\end{equation*}
as well as
\begin{equation*}
    \left| V_A^T\psi V_{A,j} \right| = |A^T\psi A_j| \lesssim 1
\end{equation*}
and
\begin{equation*}
    |V_{A,j}V_{A,l}| = |A_jA_l|\lesssim 1
\end{equation*}
by using the boundedness assumptions. Therefore, Assumption (i) of Theorem \ref{thm:master-discrete} is satisfied with $r_1(k)=1$.

\subsubsection{Second Moment Boundedness}

The second-moment boundedness assumptions follow immediately from the almost sure boundedness proved above. Hence, Assumption (ii) of Theorem \ref{thm:master-discrete} is satisfied with $\widehat r_1(k), r_1(k)=1$.

\subsubsection{Nuisance Estimation Error Rate}

We have
\begin{align*}
    &\Pb\left\{ \widehat\varphi_a(Z) - V_a^T\psi \right\} \\
    &= \Pb\left\{ \widehat\varphi_a(Z) - a^T\psi \right\} \\
    &= \Pb\left\{ \frac{\mathds{1}(A=a)}{\widehat\pi_a(X)}\{\E(Y\mid X,A) - \widehat\mu_a(X)\} + \widehat\mu_a(X) -\mu_a(X)\right\}\\
    &= \Pb\left\{\frac{\mathds{1}(A=a)}{\widehat\pi_a(X)}\{\mu_a(X) - \widehat\mu_a(X)\} + \widehat\mu_a(X)-\mu_a(X)  \right\}\\
    &= \Pb\left\{ \frac{\E(\mathds{1}(A=a)\mid X)}{\widehat\pi_a(X)}\{\mu_a(X) - \widehat\mu_a(X)\} + \widehat\mu_a(X) -\mu_a(X) \right\}\\
    &= \Pb\left\{ \frac{\pi_a(X)}{\widehat\pi_a(X)}\{\mu_a(X) - \widehat\mu_a(X)\} + \widehat\mu_a(X) - \mu_a(X)\right\} \\
    &= \Pb\left\{ \left(\mu_a(X) - \widehat\mu_a(X)\right)\left( \frac{\pi_a(X)}{\widehat\pi_a(X)} - 1 \right) \right\} \\
    &\leq \anorm{\frac{\pi_a}{\widehat\pi_a} - 1}\anorm{\mu_a - \widehat\mu_a} = r_{2,a}(n,k).
\end{align*}

\subsubsection{Restricted Eigenvalue Condition}

By assumption, we have
$$
    v^T\E(V_AV_A^T)v = v^T\E(AA^T)v \geq D\norm{v}_2^2,
$$
so Assumption (iv) of Theorem \ref{thm:master-discrete} is satisfied with $C_3=D$.

\subsubsection{Sample Covariance Condition}

By assumption, we have
$$
    \sup_{\norm{v}_2=1, \norm{v}_1\leq 2\sqrt{s}} \left|v^T\Sigma v - v^T\widehat\Sigma v \right| = \sup_{\norm{v}_2=1, \norm{v}_1\leq 2\sqrt{s}} \left|v^T\left(\E(AA^T) - \frac{1}{n}\mathbb{A}^T\mathbb{A}\right)v \right|\leq \frac{D}{2}
$$
with high probability.

\subsubsection{Final Rate}

Plugging the obtained rates into Theorem \ref{thm:master-discrete}, we obtain
\begin{align*}
    \E\left( \Pb^2_{n}\left[ \left\{ A^T(\widehat\psi - \psi) \right\}^2 \right] \right) &\lesssim s\frac{\log k}{n} + s\E\left( \max_{j\in\{1,\dots,k\}} \left(\sum_{a\in\mathcal{A}} r_{2,a}(n,k)|V_{a,j}| \widehat\varpi_a \right)^2\right) \\
    &\lesssim s\frac{\log k}{n} + s\E\left( \left( \sum_{a\in \mathcal{A}} \widehat\varpi_a \anorm{\frac{\pi_a}{\widehat\pi_a} - 1}\anorm{\mu_a - \widehat\mu_a} \right)^2\right).
\end{align*}
Finally, we can simplify the upper bound by using that $\delta_n(a)\leq\delta_n$, $\epsilon_n(a)\leq\epsilon_n$, and the fact that $\sum_{a\in\mathcal{A}} \widehat\varpi_a=1$, which gives
$$
    \E\left( \left\{ A^T(\widehat\psi - \psi) \right\}^2 \right) \lesssim s\frac{\log k}{n} + s(\delta_n\epsilon_n)^2.
$$

\subsection{Proof of Theorem \ref{thm:continuous-vector-trt}}

We prove the theorem by verifying the conditions of Theorem \ref{thm:master} for $\widehat\varphi_1(Z)=\frac{Y-\widehat\mu_A(X)}{\widehat\pi_A(X)/\widehat p(A)}$, $\widehat\varphi_{2,a}(Z)=\widehat\mu_a(X)-\psi_0$, and $\widehat\Pb=\widetilde\Pb=\widehat p\cdot\lambda$.

\subsubsection{Almost Sure Boundedness}

Using the boundedness assumption, we obtain
$$
    |\widehat\varphi_1(Z)V_{A,j}|=\left|\frac{Y-\widehat\mu_A(X)}{\widehat\pi_A(X)/\widehat p(A)}\right||A_j| \leq \frac{2B^2}{\epsilon}B=\frac{2B^3}{\epsilon},
$$
$$
    \left|\int_\mathcal{A} \widehat\varphi_{2,a}(Z)V_{a,j}d\widehat\Pb(a)\right| \leq \int_\mathcal{A} |\widehat\mu_a(X)-\psi_0|\norm{A}_\infty\widehat p(a)da \leq 2B^3,
$$
and
$$
    |V_A^T\psi V_{A,j}| = |A^T\psi A_j| \leq |A^T\psi|\norm{A}_\infty \leq 2B^2.
$$
So Assumption (i) of Theorem \ref{thm:master} is satisfied with $r_1(k)=1$.

\subsubsection{Second Moment Boundedness}

Assumption (ii) of Theorem \ref{thm:master} with $\widehat r_1(k)= r_1(k)=1$ is satisfied, as the almost sure boundedness of the previous section directly implies second moment boundedness by a constant.

\subsubsection{Nuisance Penalty}

We have
$$
    \E(\widehat\varphi_1(Z)\mid X, A=a) = \frac{\mu_a(X)-\widehat\mu_a(X)}{\widehat\pi_a(X)/\widehat p(a)}
$$
and
$$
    p_{A\mid X}(a\mid X)=\frac{\pi_a(X)}{\widehat p(a)}.
$$
Hence,
\begin{align*}
    &\left| \Pb\left\{ \E(\widehat\varphi_1(Z)\mid X, A=a)p_{A\mid X}(a\mid X) + \widehat\varphi_{2,a}(Z) - V_a^T\psi \right\} \right| \\
    &=\left| \Pb\left\{ \frac{\mu_a(X)-\widehat\mu_a(X)}{\widehat\pi_a(X)/\widehat p(a)}\frac{\pi_a(X)}{\widehat p(a)} + \widehat\mu_a(X) - \psi_0 - V_a^T\psi \right\} \right| \\
    &=\left| \Pb\left\{ \frac{\pi_a(X)}{\widehat\pi_a(X)}(\mu_a(X)-\widehat\mu_a(X)) + \widehat\mu_a(X) - \mu_a(X) \right\} \right| \\
    &=\left| \Pb\left\{ \left(\frac{\pi_a(X)}{\widehat\pi_a(X)} - 1\right)(\mu_a(X)-\widehat\mu_a(X)) \right\} \right| \\
    &\leq \anorm{\frac{\pi_a}{\widehat\pi_a} - 1}\norm{\mu_a - \widehat\mu_a} = r_{2,a}(n,k).
\end{align*}
This verifies Assumption (iii) of Theorem \ref{thm:master}.

\subsubsection{Restricted Eigenvalue Condition}

Since the eigenvalues of $\int_\mathcal{A} aa^T\widehat p(a)da$ are lower bounded by $C$, we obtain
$$
    v^T\left( \int_\mathcal{A} aa^T\widehat p(a)da \right) v \geq C\norm{v}_2^2,
$$
which verifies Assumption (iv) of Theorem \ref{thm:master}.

\subsubsection{Difference in Measures}

Since $\widehat\Pb = \widetilde\Pb$, Assumption (v) of Theorem \ref{thm:master} is immediately satisfied.

\subsubsection{Integrability}

The $\sigma$-finiteness and integrability in Assumption (vi) of Theorem \ref{thm:master} directly follows from the boundedness assumptions.

\subsubsection{Final Rate}

We verified all Assumptions of Theorem \ref{thm:master}. Therefore, we obtain
\begin{align*}
    &\E\left( \int_\mathcal{A}\left\{ a^T(\widehat\psi - \psi) \right\}^2 \widehat p(a)da \right)\\ &\lesssim s\frac{r_1(k)\log k}{n} + s\E\left( \max_{j\in\{1,\dots,k\}} \left(\int_\mathcal{A} r_{2,a}(n,k)|V_{a,j}| d\widehat\Pb(a) \right)^2\right) \\
    &=s\frac{\log k}{n} + s\E\left( \max_{j\in\{1,\dots,k\}} \left(\int_\mathcal{A} \anorm{\frac{\pi_a}{\widehat\pi_a} - 1}\norm{\mu_a - \widehat\mu_a}|A_j| \widehat p(a) da\right)^2\right) \\
    &\lesssim s\frac{\log k}{n} + s\E\left( \left(\int_\mathcal{A} \anorm{\frac{\pi_a}{\widehat\pi_a} - 1}\norm{\mu_a - \widehat\mu_a} \widehat p(a) da\right)^2\right) \\
    &\lesssim s\frac{\log k}{n} + s(\delta_n\epsilon_n)^2.
\end{align*}
This proves the theorem.

\subsection{Proof of Theorem \ref{thm:master-discrete}}

We start by decomposing the error into an expected in-sample error and an empirical process term:
\begin{equation}\label{eq:decomposition-of-error}
    \E\left( \left\{ V_A^T(\widehat\psi - \psi) \right\}^2 \right) = \E\left( \Pb^2_{n}\left[ \left\{ V_A^T(\widehat\psi - \psi) \right\}^2 \right] \right) + \E\left( (\Pb - \Pb^2_{n})\left[ \left\{ V_A^T(\widehat\psi - \psi) \right\}^2 \right] \right).
\end{equation}

First, we handle the expected in-sample error and want to bound it by the desired order, i.e., we wish to show that
\begin{align*}
    \E\left( \Pb^2_{n}\left[ \left\{ V_A^T(\widehat\psi - \psi) \right\}^2 \right] \right) &\lesssim s\frac{r_1(k)\log k}{n} + s\E\left( \max_{j\in\{1,\dots,k\}} \left(\int_\mathcal{A} r_{2,a}(n,k)|V_{a,j}| d\Pb^1_{n}(a) \right)^2\right).
\end{align*}
by verifying the assumptions of Theorem \ref{thm:master}. After, we show that the empirical process term can be bounded by a term of the same order.\\

To verify the conditions of Theorem \ref{thm:master}, we set
$$
    \widehat\varphi_1(Z)=0,\quad \widehat\varphi_{2,a}(Z)=\widehat\varphi_a(Z),\quad \widehat\Pb=\Pb_n^1, \quad \widetilde\Pb=\Pb_n^2.
$$
In the following subsections, we verify each of the assumptions.

\subsubsection{Boundedness and Nuisance Penalty}

Since $\widehat\varphi_1(Z)=0$ and $\widehat\Pb=\Pb_n^1$, Assumptions (i), (ii), and (iii) of Theorem \ref{thm:master-discrete} directly verify Assumptions (i), (ii), and (iii) of Theorem \ref{thm:master}.

\subsubsection{Restricted Eigenvalue Assumption}

Let $v\in\mathbb{R}^k$ such that $\norm{v}_2=1$ and $\norm{v}_1\leq 2\sqrt{s}$. Then, we obtain
\begin{equation*}
    v^T\Pb_n^2\left\{ V_AV_A^T \right\}v = \frac{1}{n}\norm{\mathbb{A}v}_2^2 = v^T\left(\frac{1}{n}\mathbb{A}^T\mathbb{A}\right)v = v^T\widehat\Sigma v \geq v^T\Sigma v - \frac{C_3}{2} \geq \frac{C_3}{2}
\end{equation*}
with probability at least $1-\alpha$, where we used Assumptions (iv) and (v) to obtain the inequalities. Using a rescaling argument, the sample restricted eigenvalue condition holds with high probability for all vectors $v$ such that $\norm{v}_1 \leq 2\sqrt{s}\norm{v}_2$, i.e.,
\begin{equation}\label{eq:re}
    v^T\Pb_n^2\left\{ V_AV_A^T \right\}v \geq \frac{C_3}{2}\norm{v}_2^2
\end{equation}
with probability at least $1-\alpha$. Since $\widetilde\Pb=\Pb_n^2$, we have with probability at least $1-\alpha$ that
$$
    v^T\int_\mathcal{A} V_aV_a^T d\widetilde\Pb(a) v= v^T\int_\mathcal{A} V_aV_a^T d\Pb_n^2(a) v=v^T\Pb_n^2\left\{ V_AV_A^T \right\}v \geq \frac{C_3}{2}\norm{v}_2^2,
$$
which verifies Assumption (iv) of Theorem \ref{thm:master}.

\subsubsection{Difference in Measures}

To verify Assumption (v) of Theorem \ref{thm:master}, we need to show that
$$
    \E\left[ \max_{j\in\{1,\dots,k\}} (\Pb_n^2 - \Pb_n^1)\left\{V_A^T\psi V_{A,j}\right\}^2 \right] \lesssim \frac{r_1(k)\log k}{n}.
$$
We can decompose the left-hand side into two empirical process terms, i.e.,
\begin{align*}
    &\E\left[ \max_{j\in\{1,\dots,k\}} (\Pb_n^2 - \Pb_n^1)\left\{V_A^T\psi V_{A,j}\right\}^2 \right] \\
    &\leq 2\E\left[ \max_{j\in\{1,\dots,k\}} (\Pb_n^2 - \Pb)\left\{V_A^T\psi V_{A,j}\right\}^2 \right] + 2\E\left[ \max_{j\in\{1,\dots,k\}} (\Pb_n^1 - \Pb)\left\{V_A^T\psi V_{A,j}\right\}^2 \right].
\end{align*}

$V_A^T\psi V_{A,j}$ is almost surely and second-moment bounded by $r_1(k)$ up to constants by Assumptions (i) and (ii). Hence, with an analogous analysis as in Section \ref{sec-t1} (applying Bernstein inequality and using Lemma \ref{lemma-l2}), we obtain the upper bound $\frac{r_1(k) \log k}{n}$.

\subsubsection{Integrability}

Integrability and $\sigma$-finiteness is satisfied since $\widehat\Pb=\Pb_n^1$ is discrete.

\subsubsection{Bounding the Empirical Process Term}

In the previous subsections, we verified all conditions of Theorem \ref{thm:master}. Hence, we obtain
\begin{align*}
    \E\left( \Pb^2_{n}\left[ \left\{ V_A^T(\widehat\psi - \psi) \right\}^2 \right] \right) &\lesssim s\frac{r_1(k)\log k}{n} + s\E\left( \max_{j\in\{1,\dots,k\}} \left(\sum_{a\in\mathcal{A}} r_{2,a}(n,k)|V_{a,j}| \widehat\varpi_a \right)^2\right),
\end{align*}
which means that we bounded the expected in-sample error by the desired rate. However, it remains to bound the empirical process term in the decomposition (\ref{eq:decomposition-of-error}), which is done in the following: By $\mathrm{SC}$, we denote the event that the sample covariance condition of Assumption~(v) holds. Then,
\begin{align*}
    &\E\left( (\Pb-\Pb^2_{n})\left[ \left\{ V_A^T(\widehat\psi - \psi) \right\}^2 \right] \right) \\
    &= \E\left( (\Pb-\Pb^2_{n})\left[ \left\{ V_A^T(\widehat\psi - \psi) \right\}^2 \right] \mathds{1}_{\mathrm{SC}}\right) + \E\left( (\Pb-\Pb^2_{n})\left[ \left\{ V_A^T(\widehat\psi - \psi) \right\}^2 \right] \mathds{1}_{\mathrm{SC}^c}\right).
\end{align*}
We first examine the term where the sample covariance condition holds. We have
\begin{align*}
    (\Pb-\Pb^2_{n})\left[ \left\{ V_A^T(\widehat\psi - \psi) \right\}^2 \right] &= (\Pb-\Pb^2_{n})\left\{ (\widehat\psi - \psi)^TV_AV_A^T(\widehat\psi - \psi) \right\}\\
    &= (\widehat\psi - \psi)^T(\Pb-\Pb^2_{n})\left\{ V_AV_A^T \right\}(\widehat\psi - \psi).
\end{align*}
Recalling the definitions $\Sigma=\Pb(V_AV_A^T)$ and $\widehat\Sigma = \frac{1}{n}\mathbb{A}^T\mathbb{A}$, the previous line equals
$$
    (\widehat\psi - \psi)^T(\Sigma - \widehat\Sigma)(\widehat\psi - \psi).
$$
Due to the sample covariance condition (and rescaling), we obtain
$$
    (\widehat\psi - \psi)^T(\Sigma - \widehat\Sigma)(\widehat\psi - \psi) \leq \frac{C_3}{2}\norm{\widehat\psi - \psi}_2^2.
$$
On the event that the sample covariance condition holds, we already proved the sample restricted eigenvalue condition in (\ref{eq:re}), from which we obtain
$$
    \norm{\widehat\psi - \psi}_2^2 \leq \frac{2}{C_3}\Pb^2_{n}\left[ \left\{ V_A^T(\widehat\psi - \psi) \right\}^2 \right]
$$
where we used that $\norm{\widehat\psi - \psi}_1\leq 2\sqrt{s}\norm{\widehat\psi - \psi}_2$ (which we already verified in \eqref{eq-bound-l1-by-l2}). Therefore,
$$
    (\widehat\psi - \psi)^T(\Sigma - \widehat\Sigma)(\widehat\psi - \psi) \leq \Pb^2_{n}\left[ \left\{ V_A^T(\widehat\psi - \psi) \right\}^2 \right].
$$
The expectation of the right-hand side is the expected in-sample error and can be bounded, as above, by the desired rate. Finally, we need to bound the term where the sample covariance condition is violated. We use that
\begin{equation}
\begin{split}\label{eq:emp-process-term}
    &(\Pb-\Pb^2_{n})\left[ \left\{ V_A^T(\widehat\psi - \psi) \right\}^2 \right] \\
    &= (\Pb-\Pb^2_{n})\left[ \left\{ \sum_{j=1}^k V_{A,j}(\widehat\psi - \psi)_j \right\}^2 \right]\\
    &= \sum_{j,l}(\Pb-\Pb^2_{n})\left\{ V_{A,j}V_{A,l}(\widehat\psi - \psi)_j(\widehat\psi - \psi)_l \right\} \\
    &= \sum_{j,l}(\widehat\psi - \psi)_j(\widehat\psi - \psi)_l(\Pb-\Pb^2_{n})\left\{ V_{A,j}V_{A,l} \right\}\\
    &\leq \norm{\widehat\psi - \psi}_1^2\max_{j,l}\left| (\Pb-\Pb^2_{n})\left\{ V_{A,j}V_{A,l} \right\} \right|\\
    &\leq s^2\max_{j,l}\left| (\Pb-\Pb^2_{n})\left\{ V_{A,j}V_{A,l} \right\} \right|.
\end{split}
\end{equation}
Therefore,
\begin{align*}
    \E\left( (\Pb-\Pb^2_{n})\left[ \left\{ V_A^T(\widehat\psi - \psi) \right\}^2 \right] \mathds{1}_{\mathrm{SC}^c}\right) \leq s^2\E\left[ \max_{j,l}\left| (\Pb-\Pb^2_{n})\left\{ V_{A,j}V_{A,l} \right\} \right| \mathds{1}_{\mathrm{SC}^c} \right].
\end{align*}
We can use Holder's inequality to bound the previous line by
\begin{equation}\label{eq1}
    \leq s^2\E\left[ \max_{j,l}\left( (\Pb-\Pb^2_{n})\left\{ V_{A,j}V_{A,l} \right\} \right)^4 \right]^{\frac{1}{4}} \left(1-\Pb(\mathrm{SC})\right)^{1-\frac{1}{4}}.
\end{equation}
Now, by Assumptions (i) and (ii), $|V_{A,j}V_{A,l}|\lesssim r_1(k)$ and $\Pb(V_{A,j}^2V_{A,l}^2)\lesssim r_1(k)$, which allows us to conclude
$$
    \E\left[ \max_{j,l}\left( (\Pb-\Pb^2_{n})\left\{ V_{A,j}V_{A,l} \right\} \right)^4 \right]^{\frac{1}{4}} \lesssim \sqrt{\frac{r_1(k)\log k}{n}}
$$
by using Bernstein's inequality, Lemma \ref{lemma-orlicz}, and using the fact that $\norm{X}_{L_4}\leq 4!\norm{X}_{\Psi_1}$. By assumption, we have $1-\Pb(\mathrm{SC})\lesssim \frac{\log k}{n}$. Hence, up to constants, we can further upper bound (\ref{eq1}) by
$$
    s^2\sqrt{\frac{r_1(k)\log k}{n}}\left(\frac{\log k}{n}\right)^{\frac{3}{4}} = s\frac{r_1(k)\log k}{n}\cdot \frac{s}{\sqrt{r_1(k)}}\left(\frac{\log k}{n}\right)^{\frac{1}{4}}.
$$
Since $n\gtrsim \frac{s^4}{r_1(k)^2}\log k$ by assumption, we can upper bound the previous line by the desired rate
$$
    \lesssim s\frac{r_1(k)\log k}{n}.
$$
This concludes the proof.

\subsection{Proof of Theorem \ref{thm:master-approximate-sparsity}}\label{proof-master-approximate-sparsity}

The proof follows the same logic as the proof of Theorem \ref{thm:master}. We only highlight the necessary adjustments here.\\

The basic Lasso inequality can be applied in the same way, giving us
\begin{align*}
    &= \E\left[ 2\Pb^2_{n}\left\{ \widehat\varphi_1(Z)V_A^T(\widehat\psi-\psi)+\int_\mathcal{A}\widehat\varphi_{2,a}(Z)V_a^T(\widehat\psi - \psi) d\widehat\Pb(a)\right\} + 2\int_\mathcal{A}(V_a^T\psi)V_a^T(\psi - \widehat\psi)d\widetilde\Pb(a) \right] \\
    &\leq 2\E\left[ \norm{\widehat\psi - \psi}_1 \max_{j\in\{1,\dots,k\}} \left| \Pb^2_{n}\left\{\widehat\varphi_1(Z)V_{A,j} + \int_\mathcal{A} (\widehat\varphi_{2,a}(Z)-V_a^T\psi)V_{a,j}d\widehat\Pb(a)\right\} + (\widetilde\Pb - \widehat\Pb)\left\{V_A^T\psi V_{A,j}\right\} \right| \right]
\end{align*}
Then, the proof of Theorem \ref{thm:master} further upper-bounds $\norm{\widehat\psi - \psi}_1$ using the exact sparsity and restricted eigenvalues to achieve the fast Lasso rate. However, these assumptions are not available to us here, so we simply use the approximate sparsity assumption to bound $\norm{\widehat\psi - \psi}_1\leq 2s$ instead. Consequently,
\begin{align*}
    &= \E\left[ 2\Pb^2_{n}\left\{ \widehat\varphi_1(Z)V_A^T(\widehat\psi-\psi)+\int_\mathcal{A}\widehat\varphi_{2,a}(Z)V_a^T(\widehat\psi - \psi) d\widehat\Pb(a)\right\} + 2\int_\mathcal{A}(V_a^T\psi)V_a^T(\psi - \widehat\psi)d\widetilde\Pb(a) \right] \\
    &\leq 4s\E\left[ \max_{j\in\{1,\dots,k\}} \left| \Pb^2_{n}\left\{\widehat\varphi_1(Z)V_{A,j} + \int_\mathcal{A} (\widehat\varphi_{2,a}(Z)-V_a^T\psi)V_{a,j}d\widehat\Pb(a)\right\} + (\widetilde\Pb - \widehat\Pb)\left\{V_A^T\psi V_{A,j}\right\} \right| \right]
\end{align*}
The term on the right-hand side is analogous to the one obtained in the proof of Theorem \ref{thm:master}, with the only difference that the term does not appear squared. We apply the same decomposition into the terms $T_1$, $T_2$ and $T_3$, i.e., the right-hand side equals
\begin{align*}
    &4s\underbrace{\E\left[ \max_{j\in\{1,\dots,k\}} \left|(\Pb^2_{n}-\Pb)\left\{ \widehat\varphi_1(Z)V_{A,j} + \int_\mathcal{A} (\widehat\varphi_{2,a}(Z)-V_a^T\psi)V_{a,j}d\widehat\Pb(a) \right\}\right| \right]}_{T_1} \\
    &+ 4s\underbrace{\E\left[ \max_{j\in\{1,\dots,k\}} \left| \Pb\left\{ \widehat\varphi_1(Z)V_{A,j} + \int_\mathcal{A} (\widehat\varphi_{2,a}(Z)-V_a^T\psi)V_{a,j}d\widehat\Pb(a)\right\} \right| \right]}_{T_2} \\
    &+ 4s\underbrace{\E\left[\max_{j\in\{1,\dots,k\}}\left| (\widetilde\Pb - \widehat\Pb)\left\{V_A^T\psi V_{A,j}\right\} \right| \right]}_{T_3}.
\end{align*}
Note that, also here, the terms are not of quadratic order, as is the case in the proof of Theorem~\ref{thm:master}. Nonetheless, we can use the same tools to analyze each term. $T_3$ can directly be upper bounded by the desired rate using Assumption (iv). Furthermore, we obtain
\begin{equation*}
    T_2 \leq C_2\E\left\{ \max_{j\in\{1,\dots,k\}} \int_\mathcal{A} r_{2,a}(n,k)|V_{a,j}|d\widehat\Pb(a) \right\}.
\end{equation*}
with completely analogous calculations. For the term $T_1$, we have to formulate Lemma \ref{lemma-l2} in terms of the $L_1$-norm instead of the $L_2$-norm. This can be done by using Lemma \ref{lemma-orlicz} and the fact that the Orlicz norm can be bounded below by $\norm{X}_1\leq \norm{X}_{\psi_1}$. This implies the following lemma.
\begin{lemma}\label{lemma-l1}
    Under the notation and assumptions of Lemma \ref{lemma-orlicz}, we have
    \begin{equation*}
        \E\left(\max_{1\leq i\leq k} |X_i|\right) \leq K\left( a\log(1+k) + \sqrt{b}\sqrt{\log(1+k)} \right).
    \end{equation*}
\end{lemma}
Using Bernstein's inequality in combination with Assumptions (i) and (ii) to verify the assumptions of this lemma, we obtain
\begin{equation*}
    T_1 \lesssim \frac{r_1(k)\log k}{n} + \sqrt{\frac{r_1(k)\log k}{n}}.
\end{equation*}
This proves Theorem \ref{thm:master-approximate-sparsity}.

\subsection{Proof of Theorem \ref{thm:master-approximate-sparsity-discrete}}

We verify the assumptions of Theorem \ref{thm:master-approximate-sparsity} to show the claim.

\subsubsection{Boundedness and Nuisance Penalty}

Since $\widehat\varphi_1(Z)=0$ and $\widehat\Pb=\Pb_n^1$, Assumptions (i), (ii), and (iii) of Theorem \ref{thm:master-approximate-sparsity-discrete} directly verify Assumptions (i), (ii), and (iii) of Theorem \ref{thm:master-approximate-sparsity}.

\subsubsection{Difference in Measures}

To verify Assumption (iv) of Theorem \ref{thm:master-approximate-sparsity}, we need to show that
$$
    \E\left[ \max_{j\in\{1,\dots,k\}} \left|(\Pb_n^2 - \Pb_n^1)\left\{V_A^T\psi V_{A,j}\right\}\right| \right] \lesssim \sqrt{\frac{r_1(k)\log k}{n}}.
$$
We can decompose the left-hand side into two empirical process terms, i.e.,
\begin{align*}
    &\E\left[ \max_{j\in\{1,\dots,k\}} \left|(\Pb_n^2 - \Pb_n^1)\left\{V_A^T\psi V_{A,j}\right\}\right| \right] \\
    &\leq 2\E\left[ \max_{j\in\{1,\dots,k\}} \left|(\Pb_n^2 - \Pb)\left\{V_A^T\psi V_{A,j}\right\}\right| \right] + 2\E\left[ \max_{j\in\{1,\dots,k\}} \left|(\Pb_n^1 - \Pb)\left\{V_A^T\psi V_{A,j}\right\}\right| \right].
\end{align*}

$V_A^T\psi V_{A,j}$ is almost surely and second-moment bounded by $r_1(k)$ up to constants by Assumptions (i) and (ii). Hence, with an analogous analysis as in Section \ref{proof-master-approximate-sparsity} (applying Bernstein inequality and using Lemma \ref{lemma-l1}), we obtain the upper bound $\sqrt{\frac{r_1(k) \log k}{n}}$.

\subsubsection{Integrability}

Integrability and $\sigma$-finiteness is satisfied since $\widehat\Pb=\Pb_n^1$ is discrete.

\subsubsection{Final Rate}

We verified all assumptions of Theorem \ref{thm:master-approximate-sparsity}, so we obtain
\begin{align*}
    \E\left( \Pb^2_{n}\left[ \left\{ V_A^T(\widehat\psi - \psi) \right\}^2 \right] \right) &\lesssim s\sqrt{\frac{r_1(k)\log k}{n}} + s\E\left( \max_{j\in\{1,\dots,k\}} \sum_{a\in\mathcal{A}} r_{2,a}(n,k)|V_{a,j}| \widehat\varpi_a \right).
\end{align*}
This proves Theorem \ref{thm:master-approximate-sparsity-discrete}.

\subsection{Proof of Theorem \ref{thm:single-treatment-approximate-sparsity}}

We want to derive the rate using Theorem \ref{thm:master-approximate-sparsity-discrete} for the choices
$$
    \widehat\varphi_a(Z)=\sqrt{k}\left(\frac{\mathds{1}(A=a)}{\widehat\pi_a(X)}\{Y-\widehat\mu_a(X)\}+\widehat\mu_a(X)-\psi_0\right),\quad V_a=\sqrt{k}(\mathds{1}(a=1),\dots,\mathds{1}(a=k)).
$$
To verify the assumptions, we can reuse the computations in the proof of Theorem \ref{thm:single-treatment}.

Assumption (i) of Theorem \ref{thm:master-approximate-sparsity-discrete} is satisfied with $r_1(k)=k$ by the calculation in Section \ref{subsection:bounded-weighted-pseudo-outcomes}.

Assumption (ii) is satisfied with $r_1(k)=k$ using the calculations in Section \ref{subsection:bounded-weighted-second-moments}.

Assumption (iii) can be verified with $r_{2,a}(n,k)=\sqrt{k}\anorm{\frac{\pi_a}{\widehat\pi_a} - 1}\anorm{\mu_a - \widehat\mu_a}$ according to Section \ref{subsection:nuisance-estimation-error-rate}.

Using Theorem \ref{thm:master-approximate-sparsity-discrete}, we obtain
\begin{equation*}
    \E\left(\sum_{a=1}^k w(a)(\widehat\psi_a - \psi_a)^2\right) \lesssim s\sqrt{\frac{\log k}{kn}} + s\E\left[ \max_{a\in\{1,\dots,k\}} \left\{\widehat\varpi_a \delta_n(a)\epsilon_n(a) \right\} \right].
\end{equation*}

Finally, we simplify the nuisance penalty by bounding $\delta_n(a)\leq\delta_n$ and $\epsilon_n(a)\leq\epsilon_n$, and the fact that
$$
    \E\left(\max_{a\in\{1,\dots,k\}} \widehat\varpi_a\right) \leq \sqrt{\E\left(\max_{a\in\{1,\dots,k\}} \widehat\varpi_a^2\right)} \lesssim \frac{1}{k}
$$
due to the already shown identity $\E\left(\max_{a\in\{1,\dots,k\}} \widehat\varpi_a^2\right)\asymp 1/k^2$ in Appendix \ref{subsection:final-rate}. This proves that
\begin{equation*}
    \E\left(\sum_{a=1}^k w(a)(\widehat\psi_a - \psi_a)^2\right) \lesssim s\sqrt{\frac{\log k}{kn}} + \frac{s}{k}\delta_n\epsilon_n.
\end{equation*}

\subsection{Proof of Theorem \ref{thm:binary-vector-trt-approximate-sparsity}}

We prove the result by verifying the conditions of Theorem \ref{thm:master-approximate-sparsity-discrete}. In the proof of Theorem \ref{thm:binary-vector-trt}, we already verified that Assumptions (i) and (ii) hold for $r_1(k)=1$, and Assumption (iii) holds for
\begin{equation*}
    r_{2,a}(n,k) = \anorm{\frac{\pi_a}{\widehat\pi_a} - 1}\anorm{\mu_a - \widehat\mu_a}.
\end{equation*}
Applying Theorem \ref{thm:master-approximate-sparsity-discrete} now gives
\begin{equation*}
    \E\left( \Pb^2_{n}\left[ \left\{ A^T(\widehat\psi - \psi) \right\}^2 \right] \right) \lesssim s\sqrt{\frac{\log k}{n}} + s\E\left( \sum_{a\in \mathcal{A}} \widehat\varpi_a \anorm{\frac{\pi_a}{\widehat\pi_a} - 1}\anorm{\mu_a - \widehat\mu_a} \right).
\end{equation*}
Using the bounds $\delta_n(a)\leq\delta_n$ and $\epsilon_n(a)\leq\epsilon_n$, and the fact that $\sum_{a\in\mathcal{A}} \widehat\varpi_a=1$, we can conclude the claim of the theorem.

\subsection{Proof of Theorem \ref{thm:continuous-vector-trt-approximate-sparsity}}

We prove the claim by verifying the conditions of Theorem \ref{thm:master-approximate-sparsity}. This can be done completely analogously as in the proof of the error rate for continuous vector treatments under exact sparsity, so we refer to the proof of Theorem \ref{thm:continuous-vector-trt} for details.

\section{Proofs of Lower Bounds}

\subsection{Proof of Proposition \ref{prop:classiclowerbd}}
\label{sec:lowerbd1}

To derive this result, we use  Fano's method. The following lemma gives a version based on $\chi^2$ divergence, adapted from  Theorem 2.6 of \citet{tsybakov2009introduction}.

\begin{lemma} \label{lem:fano}
Let $\{P_0,P_1,...,P_M\} \subseteq \mathcal{P}$ with $M \geq 2$. Assume
$\frac{1}{M} \sum_{j=1}^M \chi^2(P_j, P_0) \leq \alpha M$
for $0 < \alpha < 1/2$, 
and $d(\theta_j,\theta_{j'}) \geq \Delta > 0$ for all $j \neq j'$. Then 
$$ \inf_{\widehat\theta} \sup_{P \in \mathcal{P}} \E_P \left[ \ell\left\{ d \left( \widehat\theta , \theta(P)\right) \right\} \right] \geq \frac{\ell(\Delta/2)}{2} \left( 1-\alpha - \frac{1}{M} \right) $$
for any monotonic non-negative loss function $\ell$.
\end{lemma}

\subsubsection{Construction}

Consider the $2^k$ distributions $P_\omega$ where $Y \mid A=a \sim \text{Bern}(1/2 + \omega_a \epsilon)$, i.e., 
\begin{align*}
p_\omega(z) &= \Big\{ y (1/2 + \omega_a \epsilon) + (1-y) (1/2-\omega_a \epsilon) \Big\} \pi_a   
\end{align*}
where $\omega =(\omega_1,...,\omega_k) \in \{0,1\}^k$, $\pi_a\geq C'/k$ and $\epsilon$ is specified later. Note that $\mu_a(x)=1/2+\omega_a \epsilon$. These are valid distributions as long as $ \epsilon \leq 1/2$. Let $P_j=P_{\omega^{(j)}}$, $j=1,\dots,M$, denote at least $2^{k/8}$ distributions with probability mass functions $p_{\omega^{(j)}}$ for which any $\omega^{(j)} \neq  \omega^{(j')}$ differ on at least $k/8$ indices. Existence is guaranteed by the following Varshamov-Gilbert lemma (\citet{tsybakov2009introduction}, Lemma 2.9).

\begin{lemma}[Varshamov-Gilbert] \label{lem:gilbert}
Assume $k \geq 8$. Then there exists a subset $\{ \omega^{(0)},..., \omega^{(M)}\}$ of $\Omega = \{0,1\}^k$ such that
$M \geq 2^{k/8}$,  $\omega^{(0)}=(0,...,0)$, and $\| \omega^{(j)} - \omega^{(j')} \|_0 \geq k/8$ for all $j \neq j'$. 
\end{lemma}

\subsubsection{Functional Separation}

Note for any $\omega, \omega'$ differing at index $a$ we have $|\psi_a(P_\omega) - \psi_a(P_{\omega'})| =  \epsilon$. Now since any $\omega,\omega'$ must differ on at least $k/8$ indices we have
\begin{equation}
\frac{1}{k}\sum_{a=1}^k\left\{\psi_a(P_\omega) - \psi_a(P_{\omega'})\right\}^2 \geq \frac{1}{8} \epsilon^2 \label{eq:thm1sep}
\end{equation}
Note we also have for any $\omega$ that $\| \psi(P_\omega) \|_2^2  = \sum_j (1/2 + \omega_j \epsilon)^2 \leq k (1/2 + \epsilon)^2 \leq 2k(1/4 + \epsilon^2)$.   

\subsubsection{Distributional Distance}

For a single observation, the $\chi^2$ divergence satisfies
\begin{align*}
\chi^2(P_\omega,P_0) &=\sum_{a,y} \frac{p_\omega(z)^2}{p_0(z)}  -1  \leq \sum_{a,y} 2 \{y(1/2 + \epsilon) + (1-y)(1/2 - \epsilon)\}^2 \pi_a  -1 \\
&= 2 \Big\{ (1/2+\epsilon)^2 +  (1/2 - \epsilon)^2 \Big\} - 1 
= 4 \epsilon^2 .
\end{align*}
So for the product measures $P_0^{\otimes n} = \prod_{i=1}^n P_0$ and $P_\omega^{\otimes n} = \prod_{i=1}^n P_\omega$ we have 
\begin{align*}
\log \left\{ 1 + \chi^2(P_\omega^{\otimes n},P_0^{\otimes n}) \right\} &= 
\sum_{i=1}^n \log \left\{ 1 + \chi^2(P_\omega, P_0) \right\} \leq  n \log (1 +  4\epsilon^2) \leq 4n \epsilon^2
\end{align*}
using tensorization properties of $\chi^2$ and $\log(1+x) \leq x$. 
If we take 
$$\epsilon^2 = \frac{\log(M/4+1)}{4n} \leq \frac{k \log 2}{4n}$$ 
then $\log \left\{ 1 + \chi^2(P_\omega^{\otimes n},P_0^{\otimes n}) \right\} \leq \log(M/4+1)$, which implies
$\chi^2(P_\omega^{\otimes n},P_0^{\otimes n}) \leq M/4$, 
satisyfing the distance condition in Lemma \ref{lem:fano} with $\alpha=1/4$. Note that the condition $\epsilon \leq 1/2$ needed for $p_\omega$ to be a valid probability mass function is satisfied if $k \leq n/\log 2$. When $k > n /\log 2$ we can take $\epsilon$ to be a constant, e.g., $\epsilon=1/4$, and the same arguments below show a constant minimax lower bound.

\subsubsection{Final Result}

Plugging the above choice of $\epsilon$ into the functional separation \eqref{eq:thm1sep} yields 
\begin{align*}
\frac{1}{k}\sum_{a=1}^k \left\{\psi_a(P_\omega) - \psi_a(P_{\omega'})\right\}^2 \geq (1/8) \epsilon^2 = \frac{\log(M/4+1)}{32n}
\end{align*}
and note $\log (M/4 + 1) \geq \log(2^{k/8}/4) = (k/8-2) \log 2 \geq k \log 2 /16$
as long as $k \geq 32$. Lemma \ref{lem:fano} with $\alpha=1/4$ and 
$\Delta = \sqrt{C_1 k/n}$
where  $C_1 = \frac{\log 2}{512}$, and $d(x,y)=\sqrt{\frac{1}{k}\sum_{a=1}^k (x_a-y_a)^2}$ 
then yields
\begin{align*}
 \inf_{\widehat\theta} \sup_{P \in \mathcal{P}} \E_P \left( \sum_{a=1}^k\varpi_a(P)\left\{\widehat\psi_a - \psi_a(P)\right\}^2 \right) \geq \frac{C'C_1}{8} \left( 1-1/4 - \frac{1}{2^{k/8}} \right)  \frac{k}{n}  \geq C_2\frac{k}{n}
\end{align*}
for $C_2 = 11 C'C_1 / 128$
, using the fact that $1/2^{k/8} \leq 1 /16$ when $k \geq 32$, and $\varpi_a(P)\geq \frac{C'}{k}$ uniformly across all $a\in\{1,\dots,k\}$ and $P\in\mathcal{P}$. 
This proves the proposition. 

\subsection{Proof of Theorem \ref{thm:minimax-lb-sparse-regime}}

To derive this result, we use Fano's method. Recall Lemma \ref{lem:fano}, which gives a version of Fano's method based on the $\chi^2$ distance, adapted from Theorem 2.6 of \cite{tsybakov2009introduction}. Consequently, the main task will be to construct distributions $P_j, j=0,1,\dots, M$ that satisfy the conditions of this lemma. In particular, we want to ensure that the difference of the functional $\psi$ under two different distributions $P_j$ and $P_{j'}$ is large, while the distributions are still close enough, i.e., the $\chi^2$ distance of the two distributions can be bounded as required by the lemma.

\subsubsection{Construction}

Consider the distributions $P_\omega$ where $Y\mid A=a\sim\mathrm{Bern}(1/2 + \omega_a\epsilon)$, i.e.,
\begin{equation*}
    p_\omega(z) = \left\{ y(1/2 + \omega_a\epsilon) + (1-y)(1/2 - \omega_a\epsilon) \right\}\pi_a,
\end{equation*}
where $\epsilon$ will be specified later and $\omega=(\omega_1, \dots, \omega_k)\in\Omega$ with
\begin{equation*}
    \Omega = \left\{ \omega=(\omega_1,\dots,\omega_k)\in\{0,1\}^k : \norm{\omega}_0= s \right\}
\end{equation*}
the set of all $s$-sparse binary vectors, and $\pi_a$ satisfying Assumption (i). These are valid distributions as long as $\epsilon\leq 1/2$. Moreover, note that $\mu_a(x)=1/2+\omega_a\epsilon$. Let $P_0$ denote the distribution corresponding to the choice $\omega_a=0$ for all $a$. Let $P_j=P_{\omega^{(j)}}, j=1,\dots, M$ denote at least $\left(1+\frac{k}{2s}\right)^{s/8}$ distributions with densities $p_{\omega^{(j)}}$ for which any $s$-sparse vectors $\omega^{(j)}\neq\omega^{(j')}$ differ on at least $s/2$ indices. The existence is guaranteed by the following sparse version of the Varshamov-Gilbert lemma (see Lemma 4.14 of \cite{rigollet2023high} for a proof).

\begin{lemma}[Sparse Version of Varshamov-Gilbert]
    Let $1\leq s\leq k/8$. Then there exists a subset $\{\omega^{(1)},\dots.\omega^{(M)}\}\subseteq\Omega$ of $s$-sparse binary vectors such that $\norm{\omega^{(j)}-\omega^{(j')}}_0\geq \frac{s}{2}$ for all $j\neq j'$, and $\log(M)\geq \frac{s}{8}\log\left(1+\frac{k}{2s}\right)$.
\end{lemma}

Note that the distributions $P_j, j=1,\dots,M$ respect the model: If $\omega^{(j)}_a=0$ then $\psi_a(P_j)=\psi_a(P_0)=1/2$. We can only have $\psi_a(P_j)\neq 1/2$ whenever $\omega^{(j)}_a=1$, but this can only be the case for at most $s$ of the treatments $a\in\{1,\dots,k\}$ by construction. Hence, at most $s$ entries of $\psi(P_j)$ differ from $1/2$.

\subsubsection{Functional Separation}\label{sec-functional-separation-sparse}

Note for any $\omega,\omega'$ differing at index $a$ we have
\begin{equation*}
    \left| \psi_a(P_\omega) - \psi_a(P_{\omega'}) \right| = |\omega_a\epsilon - \omega'_a\epsilon| = \epsilon.
\end{equation*}
Now since any $\omega,\omega'$ must differ on at least $s/2$ indices, we conclude
\begin{equation}\label{eq-functional-separation}
    \frac{1}{k}\sum_{a=1}^k \left\{\psi_a(P_\omega) - \psi_a(P_{\omega'})\right\}^2 \geq \frac{s}{2k}\epsilon^2.
\end{equation}

\subsubsection{Distributional Distance}

For the density of a single observation, the $\chi^2$ distance satisfies
\begin{align*}
    \chi^2(P_\omega, P_0) &= \sum_{a=1}^k\sum_{y=0}^1 \frac{p_\omega(z)^2}{p_0(z)} - 1 \leq \sum_{a=1}^k\sum_{y=0}^1 2\left\{ y(1/2 + \omega_a\epsilon) + (1-y)(1/2 - \omega_a\epsilon) \right\}^2 \pi_a - 1 \\
    &= 2\sum_{a=1}^k\left\{ (1/2 + \omega_a\epsilon)^2 + (1/2 - \omega_a\epsilon)^2 \right\}\pi_a - 1 = 2\sum_{a=1}^k\left\{ 1/2 + 2\omega_a\epsilon^2 \right\}\pi_a - 1 \\
    &=\sum_{a=1}^k \pi_a + 4\epsilon^2 \sum_{a=1}^k \pi_a\omega_a - 1\leq 4C''\epsilon^2 \frac{s}{k}.
\end{align*}
So for the product measures $P_0^{\otimes n}=\prod_{i=1}^n P_0$ and $P_\omega^{\otimes n}=\prod_{i=1}^n P_\omega$, we have
\begin{align*}
    \log\left\{ 1+\chi^2(P_\omega^{\otimes n}, P_0^{\otimes n}) \right\} &= \log\left[ \prod_{i=1}^n \left\{ 1 + \chi^2(P_\omega, P_0) \right\} \right] = \sum_{i=1}^n \log\{ 1 + \chi^2(P_\omega, P_0) \} \\
    &\leq n\log\left(1 + 4C''\epsilon^2 \frac{s}{k} \right) \leq n4C''\epsilon^2 \frac{s}{k}
\end{align*}
using the tensorization property of the $\chi^2$ distance in the first line, and $\log(1+x)\leq x$ in the last. We proceed with a case distinction now.

\paragraph{Case $k\log(k/s)\leq \frac{C''}{2}n$:}
In this case, we take
\begin{equation*}
    \epsilon^2 = \frac{k\log(M/4 + 1)}{4C''sn},
\end{equation*}
then $\log\left\{ 1+\chi^2(P_\omega^{\otimes n}, P_0^{\otimes n})\right\}\leq \log(M/4+1)$, which implies
\begin{equation*}
    \chi^2\left(P_\omega^{\otimes n}, P_0^{\otimes n}\right) \leq \frac{M}{4},
\end{equation*}
satisfying the distance condition in Lemma \ref{lem:fano} with $\alpha=1/4$. Note that this is a valid choice for $\epsilon$, as
\begin{align*}
    \epsilon &= \sqrt{\frac{k\log\left(\frac{M}{4}+1\right)}{4C''ns}} \leq \sqrt{\frac{k\log\left(\binom{k}{s}/4+1\right)}{4C''ns}} \leq \sqrt{\frac{k\log\left(\left(\frac{ek}{s}\right)^s/4+1\right)}{4C''ns}} \\
    &\leq \sqrt{\frac{ks\log\left(\frac{ek}{s}\right)}{4C''ns}}\leq \sqrt{\frac{k\log\left(\frac{k}{s}\right)}{2C''n}}\leq \frac{1}{2},
\end{align*}
where we used the assumption $k\log(k/s)\leq \frac{C''}{2}n$ as well as $k\geq 8s$.

\paragraph{Case $k\log(k/s) > \frac{C''}{2}n$:} In this case, we choose $\epsilon=1/34$.
Then, we have
\begin{align*}
    \log\left\{ 1+\chi^2(P_\omega^{\otimes n}, P_0^{\otimes n}) \right\} \leq  n4C''\epsilon^2 \frac{s}{k} = \frac{C''}{289}\frac{ns}{k} \leq \frac{s}{144}\log(k/s).
\end{align*}
By the property of the sparse Varshamov-Gilbert construction, we have
\begin{align*}
    \log(M/4+1) &\geq \log(M) - \log(4) \geq \frac{s}{8}\log\left(1+\frac{k}{2s}\right) - \log(4) \geq \left(\frac{s}{8} - 1\right)\log\left(1+\frac{k}{2s}\right) \\
    &\geq \frac{s}{72}\log\left(\frac{k}{2s}\right) = \frac{s}{72}\left(\log\left(\frac{k}{s}\right) - \log(2)\right) \geq \frac{s}{144}\log\left(\frac{k}{s}\right)
\end{align*}
as long as $s\geq 9$ and $k\geq 8s$. Therefore,
\begin{equation*}
    \log\left\{ 1+\chi^2(P_\omega^{\otimes n}, P_0^{\otimes n}) \right\} \leq \log(M/4+1),
\end{equation*}
hence
\begin{equation*}
    \chi^2(P_\omega^{\otimes n}, P_0^{\otimes n}) \leq \frac{M}{4}.
\end{equation*}
Finally, note that this choice of $\epsilon$ leads to a valid distribution $p_\omega$, as $\varepsilon \leq 1/2$.

\subsubsection{Final Minimax Lower Bound}

Since we used different choices of $\epsilon$ for the case $k\log(k/s)\leq C''n/2$ and $k\log(k/s) > C''n/2$, we also make this case distinction here and derive the minimax lower bounds separately.

\paragraph{Case $k\log(k/s)\leq \frac{C''}{2}n$:} Plugging in the above choice of $\epsilon$ into the functional separation (\ref{eq-functional-separation}) yields
\begin{align*}
    \frac{1}{k}\sum_{a=1}^k \left\{\psi_a(P_\omega) - \psi_a(P_{\omega'})\right\}^2 &\geq \frac{s}{2k}\frac{k\log(M/4 + 1)}{4C''sn} = \frac{\log(M/4 + 1)}{8C''n} \geq \frac{s\log(k/s)}{1152C''n} \\
    &\geq C_1 s\frac{\log(k/s)}{n} =: \Delta_1,
\end{align*}
as long as $s\geq 9$ and $k\geq 8s$, where $C_1=\frac{1}{1152C''}$.

\paragraph{Case $k\log(k/s) > \frac{C''}{2}n$:} In this case, plugging in $\epsilon$ into the functional separation gives
\begin{equation*}
    \frac{1}{k}\sum_{a=1}^k \left\{\psi_a(P_\omega) - \psi_a(P_{\omega'})\right\}^2 \geq \frac{s}{2k}\frac{1}{1156} = C_2\frac{s}{k} =:\Delta_2,
\end{equation*}
where $C_2=\frac{1}{2312}$.\\

Now, considering both cases combined, we obtain
\begin{equation*}
    \sqrt{\frac{1}{k}\sum_{a=1}^k \left\{\psi_a(P_\omega) - \psi_a(P_{\omega'})\right\}^2}\geq \sqrt{\min(\Delta_1, \Delta_2)}\geq \sqrt{\min(C_1, C_2)\min\left(s\frac{\log(k/s)}{n}, \frac{s}{k}\right)}=: \Delta.
\end{equation*}
Applying Lemma \ref{lem:fano} with this choice of $\Delta$, $\alpha=1/2$, and $d(x,y)=\sqrt{\frac{1}{k}\sum_{a=1}^k (x_a-y_a)^2}$ gives
\begin{align*}
    &\inf_{\widehat\psi} \sup_{P\in\mathcal{P}} \E\left\{\sum_{a=1}^k \varpi_a(P)(\widehat\psi_a - \psi_a)^2\right\} \\
    &\geq C'\left(1-\frac{1}{4} - \frac{1}{M}\right)\frac{(\Delta/2)^2}{2} \\
    &\geq C'\frac{\min(C_1, C_2)}{8}\left(\frac{3}{4}-\left(1+\frac{k}{2s}\right)^{-s/8}\right)\cdot \min\left(s\frac{\log(k/s)}{n}, \frac{s}{k}\right) \\
    &\geq C'\frac{\min(C_1, C_2)}{8}\left(\frac{3}{4}-5^{-9/8}\right)\min\left(s\frac{\log(k/s)}{n}, \frac{s}{k}\right)\\
    &= C_3\cdot \min\left(s\frac{\log(k/s)}{n}, \frac{s}{k}\right)
\end{align*}
for $C_3=\frac{C'\min(C_1, C_2)}{8}\left(\frac{3}{4}-5^{-9/8}\right)$, where we used $M\geq \left(1+\frac{k}{2s}\right)^{s/8}$, $k/s\geq 8$, $s\geq 9$, and the fact that $\varpi_a(P)\geq \frac{C'}{k}$ uniformly across all $a\in\{1,\dots,k\}$ and $P\in\mathcal{P}$. This proves the theorem.

\subsection{Proof of Theorem \ref{thm:minimax-lb-structure-agnostic}}

Without loss of generality, we can assume $\delta_n, \epsilon_n\leq 1$ (the sequences converge to zero and are therefore bounded, and the constants $C_\pi, C_\mu$ can be adjusted accordingly.) Also, note that the density of a distribution in our model can be written as
\begin{equation*}
    p(z)=f(x)\pi_a(x)\mu_a(x).
\end{equation*}

To prove the theorem, we use Le Cam's method with fuzzy hypotheses. The following lemma presents a version based on the Hellinger distance, adapted from Theorem 2.15 of \citet{tsybakov2009introduction}. \\

\begin{lemma}[\citet{tsybakov2009introduction}] \label{lem:minimax}
Let $P_\lambda$ and $Q_\lambda$ denote distributions in $\mathcal{P}$ indexed by a vector $\lambda=(\lambda_1,\dots,\lambda_k)$, with $n$-fold products denoted by $P_\lambda^n$ and $Q_\lambda^n$, respectively.  
Let $\varpi$ denote a prior distribution over $\lambda$.  
If 
$$ H^2\left(\int P_\lambda^n \ d\varpi(\lambda), \int Q_\lambda^n \ d\varpi(\lambda) \right) \leq \alpha < 2$$ 
and 
$ d(\theta(P_\lambda), \theta(Q_{\lambda'})) \geq  \Delta >0 $
for a semi-distance $d$, functional $\theta: \mathcal{P} \mapsto \R^k$ and for all $\lambda,\lambda'$, 
then
$$ \inf_{\widehat\theta} \sup_{P \in \mathcal{P}} \E_P \left\{ \ell\left( d \left( \widehat\theta , \theta(P)\right) \right) \right\} \geq \ell(\Delta/2) \left( \frac{1 - \sqrt{\alpha(1-\alpha/4)}}{2} \right) $$
for any monotonic non-negative loss function $\ell$. \\
\end{lemma}
Consequently, the main task will be to construct distributions $P_\lambda$ and $Q_\lambda$ that satisfy the conditions of the previous lemma. In particular, we must ensure that the difference in the functionals $\psi$ under these distributions is maximally large while ensuring that the distributions are sufficiently close in Hellinger distance. This can be achieved with a construction similar to the one in \cite{jin2025structureagnostic}.

\subsubsection{Construction}

Let $S=\{a\in\{1,\dots,k\}\mid \psi_a\neq \theta\}$, where $k\notin S$. First, we define the null distribution $P_\lambda$. Let the density of $P_\lambda$ be
\begin{equation*}
    p_\lambda(z)=\begin{cases}
        \widehat\pi_a(x)\widehat\mu_a(x)^y\{ 1 - \widehat\mu_a(x) \}^{1-y} & \text{if } a\in S \\
        \widehat\pi_a(x)\widetilde\mu_a(x)^y\{ 1 - \widetilde\mu_a(x) \}^{1-y} & \text{if } a\notin S \\
    \end{cases}
\end{equation*}
for $x\in[0,1]$, where $\widetilde\mu_a$ is a function satisfying $\norm{\widehat\mu_a - \widetilde\mu_a}\leq C_\mu\epsilon_n$, and $\int\widetilde\mu_a(x)d\Pb(x)=\theta$ for all $a\notin S$. Then, the parameter under the null distribution is
\begin{equation*}
    \psi_a(P_\lambda) = \begin{cases}
        \int \widehat\mu_a(x)dx & \text{if } a\in S \\
        \theta & \text{if } a\notin S. \\
    \end{cases}
\end{equation*}

Next, we define the alternative distribution $Q_\lambda$. This must be done separately for the cases $\delta_n\lesssim\epsilon_n$ and $\epsilon_n\lesssim\delta_n$.

\paragraph{Case $\delta_n\lesssim\epsilon_n$:} Assume that $\delta_n\leq D\epsilon_n$ for all $n\in\mathbb{N}$ for some constant $D>0$. Let $\lambda=(\lambda_1,\dots,\lambda_m)$ with $\lambda_j\in\{1, -1\}$ for all $j$. We define the density of the alternative distribution by setting
\begin{align*}
    \pi_{a,\lambda}(x) &= \widehat\pi_a(x)\left\{ 1-\frac{\alpha}{\widehat\mu_a(x)}\sum_{j=1}^m \lambda_j B_j(x) \right\}, \\
    \mu_{a,\lambda}(x) &= \frac{\widehat\mu_a(x) + \beta\sum_{j=1}^m \lambda_j B_j(x)}{1-\frac{\alpha}{\widehat\mu_a(x)}\sum_{j=1}^m \lambda_j B_j(x)}
\end{align*}
for $a\in S$, where
\begin{align*}
    \alpha &= \frac{\varepsilon^2 C_\pi}{\max\{1,2C_\pi, 2D C_\pi/C_\mu\}}\delta_n, \\
    \beta &= \frac{C_\mu}{4}\epsilon_n,
\end{align*}
and $B_j(x)=\mathds{1}(x\in b_{2j}) - \mathds{1}(x\in b_{2j-1}), j=1,\dots, m$ for $2m$ cubes $b_j\subseteq\mathbb{R}^d$ such that $b_j\cap b_{j'}=\emptyset$, $[0,1]^d=\bigcup_{j=1}^{2m} b_j$, and
\begin{align*}
    \int B_j(x)dx &= 0, \\
    \int B_j(x)^2dx &= 2\mathrm{vol}(b_1)=\frac{1}{m}.
\end{align*}
Note that this also implies $\int B_j(x)^3dx = \int B_j(x)dx = 0$. When $a\notin S$ and $a\neq k$, we set
\begin{align*}
    \pi_{a,\lambda}(x)&=\widehat\pi_a(x),\\
    \mu_{a,\lambda}(x)&=\widetilde\mu_a(x).
\end{align*}
When $a=k$, we define
\begin{align*}
    \pi_{a,\lambda}(x)&=1-\sum_{j=1}^{k-1} \pi_{j,\lambda}(x),\\
    \mu_{a,\lambda}(x)&=\widetilde\mu_a(x).
\end{align*}
The constructed distribution lies in the model because, for $a\in S$,
\begin{align*}
    \anorm{\frac{\pi_{a,\lambda}}{\widehat\pi_a}-1}^2=\int \frac{\alpha^2}{\widehat\mu_a(x)^2}\sum_{j=1}^m B_j(x)^2 dx \leq \frac{\alpha^2}{\varepsilon^2} \sum_{j=1}^m \int B_j(x)^2dx = \frac{\alpha^2}{\varepsilon^2}m\frac{1}{m}=\frac{\alpha^2}{\varepsilon^2}\leq (C_\pi \delta_n)^2
\end{align*}
and
\begin{align*}
    \norm{\widehat\mu_a - \mu_{a,\lambda}}^2 &= \int \left(\widehat\mu_a(x) - \frac{\widehat\mu_a(x)+\beta\sum_{j=1}^m \lambda_j B_j(x)}{1-\frac{\alpha}{\widehat\mu_a(x)}\sum_{j=1}^m \lambda_j B_j(x)}\right)^2 dx \\
    &= \int \left( -\frac{\alpha\sum_{j=1}^m\lambda_j B_j(x)}{1-\frac{\alpha}{\widehat\mu_a(x)}\sum_{j=1}^m \lambda_j B_j(x)} - \beta\frac{\sum_{j=1}^m \lambda_j B_j(x)}{1-\frac{\alpha}{\widehat\mu_a(x)}\sum_{j=1}^m \lambda_j B_j(x)}\right)^2 dx \\
    &= (\alpha+\beta)^2 \int \frac{\sum_{j=1}^m B_j(x)^2}{\left(1-\frac{\alpha}{\widehat\mu_a(x)}\sum_{j=1}^m \lambda_j B_j(x)\right)^2} \\
    &\leq 4(\alpha+\beta)^2\leq 16\beta^2 \leq (C_\mu\epsilon_n)^2,
\end{align*}
where the last inequality also used the fact that
\begin{equation*}
    \frac{\alpha}{\widehat\mu_a(x)}\sum_{j=1}^m \lambda_j B_j(x) \leq \frac{\alpha}{\varepsilon} \leq \frac{1}{2}
\end{equation*}
due to our choice of $\alpha$. Moreover, the propensity scores are almost uniform, as
\begin{equation*}
    \frac{3C_1}{4k}\leq\widehat\pi_a(X)\left( 1-\frac{\alpha}{\varepsilon} \right) \leq \pi_{a,\lambda}(X) \leq \widehat\pi_a(X)\left(1+\frac{\alpha}{\varepsilon}\right)\leq \frac{5C_2}{4k}.
\end{equation*}
For $a\notin S, a\neq k$, the required rate conditions are immediately satisfied, and the propensity scores are nearly uniform by assumption. For $a=k$, $\norm{\widehat\mu_a - \mu_{a,\lambda}}^2\leq(C_\mu\epsilon_n)^2$ and
\begin{align*}
    \anorm{\frac{\pi_{k,\lambda}}{\widehat\pi_k} - 1}^2 &= \anorm{\frac{1-\sum_{t\notin S, t\neq k} \widehat\pi_t(x) -\sum_{t\in S} \widehat\pi_t(x)\left\{1-\frac{\alpha}{\widehat\mu_t(x)}\sum_{j=1}^m \lambda_j B_j(x)\right\}}{1-\sum_{t=1}^{k-1} \widehat\pi_t(x)} - 1}^2 \\
    &=\anorm{\frac{1-\sum_{t=1}^{k-1} \widehat\pi_t(x) +\sum_{t\in S} \widehat\pi_t(x)\frac{\alpha}{\widehat\mu_t(x)}\sum_{j=1}^m \lambda_j B_j(x)}{1-\sum_{t=1}^{k-1} \widehat\pi_t(x)} - 1}^2\\
    &=\anorm{\frac{\sum_{t\in S} \widehat\pi_t(x)\frac{\alpha}{\widehat\mu_t(x)}\sum_{j=1}^m \lambda_j B_j(x)}{\widehat\pi_k(x)}}^2 \\
    &\leq\frac{1}{\varepsilon^2}\int \left(\sum_{j=1}^m \lambda_j B_j(x)\right)^2 \left(\sum_{t\in S} \widehat\pi_t(x)\frac{\alpha}{\widehat\mu_t(x)}\right)^2 dx\\
    &\leq \frac{1}{\varepsilon^2}\int \sum_{j=1}^m B_j(x)^2 \frac{\alpha^2}{\varepsilon^2} dx = \frac{\alpha^2}{\varepsilon^4} \leq (C_\pi\delta_n)^2.
\end{align*}
In addition, we have
\begin{equation*}
    \pi_{k,\lambda}(x)=\widehat\pi_k(x) +\sum_{t\in S} \widehat\pi_t(x)\frac{\alpha}{\widehat\mu_t(x)}\sum_{j=1}^m \lambda_j B_j(x)\geq \varepsilon - \frac{\alpha}{\varepsilon}\geq \frac{\varepsilon}{2}=\varepsilon'
\end{equation*}
for almost every $x$.

\paragraph{Case $\epsilon_n\lesssim\delta_n$:} Suppose that $\epsilon_n\leq D\delta_n$ for some constant $D>0$ and for all $n\in\mathbb{N}$. In this case, we define the alternative, for $a\in S$, as
\begin{align*}
    \mu_{a,\lambda}(x) &= \frac{\widehat\mu_a(x)}{1+\frac{\beta}{\widehat\mu_a(x)}\sum_{j=1}^m \lambda_j B_j(x) - \alpha\beta},\\
    \pi_{a,\lambda}(x) &= \left\{1+\frac{\beta}{\widehat\mu_a(x)}\sum_{j=1}^m \lambda_j B_j(x) - \alpha\beta\right\}\left\{ \widehat\pi_a(x) + \alpha\widehat\pi_a(x)\widehat\mu_a(x)\sum_{j=1}^m \lambda_j B_j(x) \right\}.
\end{align*}
We use analogous definitions for $a\notin S$ as in the previous case $\delta_n\lesssim \epsilon_n$. Now, choose
\begin{align*}
    \alpha &= \frac{\varepsilon C_\pi}{4\max\{1,C_\pi\}}\delta_n, \\
    \beta &= \frac{\varepsilon^2 C_\mu}{4\max\{1,C_\mu, DC_\mu/C_\pi, DC_\mu\}}\epsilon_n.
\end{align*}
Then, for $a\in S$,
\begin{equation*}
    \norm{\widehat\mu_a - \mu_{a,\lambda}}^2 \leq 4\int \left(\beta\sum_{j=1}^m \lambda_j B_j(x) - \alpha\beta\widehat\mu_a(x)\right)^2dx \leq 8(\beta^2 + \alpha^2\beta^2)\leq 16\beta^2\leq (C_\mu\epsilon_n)^2
\end{equation*}
and
\begin{align*}
    \anorm{\frac{\pi_{a,\lambda}}{\widehat\pi_a} - 1}^2 &= \anorm{\left\{1+\frac{\beta}{\widehat\mu_a(x)}\sum_{j=1}^m \lambda_j B_j(x) - \alpha\beta\right\}\left\{ 1 + \alpha\widehat\mu_a(x)\sum_{j=1}^m \lambda_j B_j(x) \right\} - 1}^2 \\
    &\leq \anorm{1+\frac{\beta}{\varepsilon} + \frac{\alpha\beta}{\varepsilon} + \frac{\alpha}{\varepsilon}- 1}^2 \leq \frac{(\beta + \alpha\beta + \alpha)^2}{\varepsilon^2}\leq (C_\pi\delta_n)^2.
\end{align*}
Note that the propensity scores are nearly uniform under this alternative distribution since, for $a\in S$,
\begin{equation*}
    \frac{C_1}{4k}\leq\widehat\pi_a(X)\left(1-\frac{\beta}{\varepsilon}-\alpha\beta\right)\left(1-\alpha\right)\leq \pi_{a,\lambda}(X) \leq \widehat\pi_a(X)\left(1+\frac{\beta}{\varepsilon}+\alpha\beta\right)\left(1+\alpha\right)\leq \frac{9C_2}{4k}.
\end{equation*}
When $a=k$, we obtain
\begin{align*}
    &\anorm{\frac{\pi_{k,\lambda}}{\widehat\pi_k} - 1}^2 \\
    &= \anorm{\frac{1-\sum_{t\notin S, t\neq k} \widehat\pi_t(x) -\sum_{t\in S} \left\{1+\frac{\beta}{\widehat\mu_t(x)}\sum_{j=1}^m \lambda_j B_j(x) - \alpha\beta\right\}\left\{ \widehat\pi_t(x) + \alpha\widehat\pi_t(x)\widehat\mu_t(x)\sum_{j=1}^m \lambda_j B_j(x) \right\}}{1-\sum_{t=1}^{k-1} \widehat\pi_t(x)} - 1}^2 \\
    &= \anorm{\frac{\sum_{t\in S} \left(\left\{\frac{\beta}{\widehat\mu_t(x)}\sum_{j=1}^m \lambda_j B_j(x) - \alpha\beta\right\}\left\{ \widehat\pi_t(x) + \alpha\widehat\pi_t(x)\widehat\mu_t(x)\sum_{j=1}^m \lambda_j B_j(x) \right\} + \alpha\widehat\pi_t(x)\widehat\mu_t(x)\sum_{j=1}^m \lambda_j B_j(x) \right)}{\widehat\pi_k}}^2 \\
    &\leq \frac{1}{\varepsilon^2}\anorm{\sum_{t\in S} \left( \left\{\frac{\beta}{\widehat\mu_t(x)}\sum_{j=1}^m \lambda_j B_j(x) - \alpha\beta\right\}\left\{ \widehat\pi_t(x) + \alpha\widehat\pi_t(x)\widehat\mu_t(x)\sum_{j=1}^m \lambda_j B_j(x) \right\} + \alpha\widehat\pi_t(x)\widehat\mu_t(x)\sum_{j=1}^m \lambda_j B_j(x) \right)}^2 \\
    &\leq \frac{1}{\varepsilon^2}\anorm{\sum_{t\in S} \left\{ \widehat\pi_t(x)\left( \frac{\beta}{\varepsilon} + \alpha^2\beta(1-\varepsilon) \right) + 2\widehat\pi_t(x)\alpha\beta + \widehat\pi_t(x)\alpha(1-\varepsilon)\right\}}^2 \\
    &=\frac{1}{\varepsilon^2}\int \left( \frac{\beta}{\varepsilon} + \alpha^2\beta + 2\alpha\beta + \alpha\right)^2\left( \sum_{t\in S}\widehat\pi_t(x) \right)^2 dx \\
    &\leq \left( \frac{\beta}{\varepsilon^2} + \frac{\beta}{\varepsilon} + 2\frac{\beta}{\varepsilon} + \frac{\alpha}{\varepsilon}\right)^2 \leq (C_\pi\delta_n)^2.
\end{align*}
Additionally,
\begin{align*}
    &\pi_{k,\lambda}(x)\\
    &=\widehat\pi_k(x) - \sum_{t\in S} \left( \left\{\frac{\beta}{\widehat\mu_t(x)}\sum_{j=1}^m \lambda_j B_j(x) - \alpha\beta\right\}\left\{ \widehat\pi_t(x) + \alpha\widehat\pi_t(x)\widehat\mu_t(x)\sum_{j=1}^m \lambda_j B_j(x) \right\} + \alpha\widehat\pi_t(x)\widehat\mu_t(x)\sum_{j=1}^m \lambda_j B_j(x) \right) \\
    &\geq \varepsilon - \left(\frac{\beta}{\varepsilon} + \alpha^2\beta + 2\alpha\beta + \alpha\right)\geq \varepsilon - \frac{3}{4}\varepsilon = \frac{\varepsilon}{4} = \varepsilon'.
\end{align*}

\subsubsection{Functional Separation}\label{sec-functional-separation-sa}

\paragraph{Case $\delta_n\lesssim\epsilon_n$:}
Under $Q_\lambda$, for $a\in S$, the functional is
\begin{align*}
    \int \frac{\widehat\mu_a(x) + \beta\sum_{j=1}^m \lambda_j B_j(x)}{1-\frac{\alpha}{\widehat\mu_a(x)}\sum_{j=1}^m \lambda_j B_j(x)} d\Pb(x) = \int \left\{\widehat\mu_a(x) + \beta\sum_{j=1}^m \lambda_j B_j(x)\right\}\sum_{k=0}^\infty \left\{\frac{\alpha}{\widehat\mu_a(x)}\sum_{j=1}^m \lambda_j B_j(x)\right\}^k dx,
\end{align*}
by using the fact that $\frac{\alpha}{\widehat\mu_a(x)}\sum_{j=1}^m \lambda_j B_j(x)<1$. Considering the first three summands of the series separately, the previous line equals
\begin{align*}
    = &\int \widehat\mu_a(x) + \left(\alpha+\beta\right)\sum_{j=1}^m \lambda_j B_j(x) + \frac{\alpha\beta+\alpha^2}{\widehat\mu_a(x)}\sum_{j=1}^m B_j(x)^2 +\frac{\beta\alpha^2}{\widehat\mu_a(x)^2}\sum_{j=1}^m \lambda_jB_j(x)^3\\
    &+\left\{\widehat\mu_a(x)+\beta\sum_{j=1}^m \lambda_j B_j(x)\right\}\sum_{k=3}^\infty \left\{\frac{\alpha}{\widehat\mu_a(x)}\sum_{j=1}^m \lambda_j B_j(x)\right\}^k dx.
\end{align*}
Therefore, the difference in the functionals under $P_\lambda$ and $Q_\lambda$ is
\begin{align*}
    &\int \frac{\alpha\beta+\alpha^2}{\widehat\mu_a(x)} + \frac{\beta\alpha^2}{\widehat\mu_a(x)^2}\sum_{j=1}^m \lambda_jB_j(x)^3 + \left\{ \widehat\mu_a(x) + \beta\sum_{j=1}^m \lambda_j B_j(x)\right\}\sum_{k=3}^\infty \left\{ \frac{\alpha}{\widehat\mu_a(x)}\sum_{j=1}^m \lambda_jB_j(x) \right\}^k dx \\
    &\geq \int\frac{\alpha\beta+\alpha^2}{\widehat\mu_a(x)} -\frac{\beta\alpha^2}{\widehat\mu_a(x)^2} + \left\{ \widehat\mu_a(x) + \beta\sum_{j=1}^m \lambda_j B_j(x)\right\}\sum_{k=3}^\infty \left\{ \frac{\alpha}{\widehat\mu_a(x)}\sum_{j=1}^m \lambda_jB_j(x) \right\}^k dx.
\end{align*}
Since $\frac{\alpha}{\widehat\mu_a(x)}\sum_{j=1}^m \lambda_jB_j(x)\leq 1/2$, we can use the formula $\sum_{k=3}^\infty t^k=\frac{t^3}{1-t}$. Hence, the previous line equals
\begin{align*}
    &= \int \frac{\alpha\beta+\alpha^2}{\widehat\mu_a(x)} -\frac{\beta\alpha^2}{\widehat\mu_a(x)^2} + \left\{ \widehat\mu_a(x) + \beta\sum_{j=1}^m \lambda_j B_j(x)\right\}\frac{\left\{ \frac{\alpha}{\widehat\mu_a(x)}\sum_{j=1}^m \lambda_jB_j(x) \right\}^3}{1-\frac{\alpha}{\widehat\mu_a(x)}\sum_{j=1}^m \lambda_jB_j(x)} dx \\
    &= \int \frac{\alpha\beta+\alpha^2}{\widehat\mu_a(x)} -\frac{\beta\alpha^2}{\widehat\mu_a(x)^2} + \left\{ \widehat\mu_a(x) + \beta\sum_{j=1}^m \lambda_j B_j(x)\right\}\frac{ \frac{\alpha^3}{\widehat\mu_a(x)^3}\sum_{j=1}^m \lambda_jB_j(x) }{1-\frac{\alpha}{\widehat\mu_a(x)}\sum_{j=1}^m \lambda_jB_j(x)} dx \\
    &\geq \int \frac{\alpha\beta+\alpha^2}{\widehat\mu_a(x)} + -\frac{\beta\alpha^2}{\widehat\mu_a(x)^2} \widehat\mu_a(x)\frac{ -\frac{\alpha^3}{\widehat\mu_a(x)^3} }{1-\frac{\alpha}{\widehat\mu_a(x)}\sum_{j=1}^m \lambda_jB_j(x)} dx + \int \frac{ \frac{\beta\alpha^3}{\widehat\mu_a(x)^3}\sum_{j=1}^m B_j(x)^2 }{1-\frac{\alpha}{\widehat\mu_a(x)}\sum_{j=1}^m \lambda_jB_j(x)} dx \\
    &\geq \int \frac{\alpha\beta+\alpha^2}{\widehat\mu_a(x)} -\frac{\beta\alpha^2}{\widehat\mu_a(x)^2} - \frac{2\alpha^3}{\widehat\mu_a(x)^2} dx\\
    &\geq \int \frac{(\alpha\beta+\alpha^2)\varepsilon - (\alpha\beta + \alpha^2)\alpha - \alpha^3}{\widehat\mu_a(x)^2} dx.
\end{align*}
Now, since $\alpha\leq\varepsilon/2$ and $\alpha^3\leq(\varepsilon/2)\alpha^2$, we can bound the functional separation by
\begin{equation*}
    \geq \int \frac{(\alpha\beta)(\varepsilon/2)}{\widehat\mu_a(x)^2}dx \geq \frac{\varepsilon}{2}\alpha\beta.
\end{equation*}
for $a\in S$.

For $a\notin S$, the function is equal to $\theta$ under both the null and alternative distributions. Hence, the functional separation of $\psi_a$ is zero in that case.

\paragraph{Case $\epsilon_n\lesssim\delta_n$:} In this case, for $a\in S$, the difference of the functional equals
\begin{align*}
    &\int \left(\frac{\widehat\mu_a(x)}{1+\frac{\beta}{\widehat\mu_a(x)}\sum_{j=1}^m \lambda_j B_j(x) - \alpha\beta} - \widehat\mu_a(x)\right)dx \\
    &=\int \widehat\mu_a(x)\sum_{k=0}^\infty \beta^k\left\{\alpha - \frac{1}{\widehat\mu_a(x)}\sum_{j=1}^m \lambda_j B_j(x)\right\}^k - \widehat\mu_a(x)dx\\
    &=\int \widehat\mu_a(x)\beta\left\{\alpha - \frac{1}{\widehat\mu_a(x)}\sum_{j=1}^m \lambda_j B_j(x)\right\} + \widehat\mu_a(x)\sum_{k=2}^\infty \beta^k\left\{\alpha - \frac{1}{\widehat\mu_a(x)}\sum_{j=1}^m \lambda_j B_j(x)\right\}^kdx\\
    &\geq \varepsilon \alpha\beta + \int \widehat\mu_a(x)\frac{\beta^3\left\{\alpha - \frac{1}{\widehat\mu_a(x)}\sum_{j=1}^m \lambda_j B_j(x)\right\}^3}{1-\beta\left\{\alpha - \frac{1}{\widehat\mu_a(x)}\sum_{j=1}^m \lambda_j B_j(x)\right\}}dx \\
    &\geq \varepsilon \alpha\beta - \frac{2\beta^3}{\varepsilon} \geq \varepsilon \alpha\beta - 2\frac{\beta}{\varepsilon}\alpha\beta \geq \varepsilon \alpha\beta - \frac{\varepsilon}{2}\alpha\beta = \frac{\varepsilon}{2}\alpha\beta
\end{align*}
where obtaining the equality in the second line and the inequality in the fourth line relies on
\begin{equation*}
    \beta\left\{\alpha - \frac{1}{\widehat\mu_a(x)}\sum_{j=1}^m \lambda_j B_j(x)\right\}\leq \alpha\beta + \frac{\beta}{\varepsilon}\leq \frac{1}{4}+\frac{1}{4}=\frac{1}{2}.
\end{equation*}

\subsubsection{Distributional Distance}\label{sec-distributional-distance-sa}

Next, we verify that the constructed distributions $P_\lambda$ and $Q_\lambda$ are close. More specifically, our goal is to bound the Hellinger distance by $1/2$. This is done using the following lemma from \citet{robins2009semiparametric}, which allows one to bound the Hellinger distance between $n$-fold products of mixtures using single-observation densities.

\begin{lemma}\label{lemma-bound-hellinger-distance}
    Let $P_\lambda$ and $Q_\lambda$ denote distributions indexed by a vector $\lambda=(\lambda_1,\dots,\lambda_m)$, and let $\mathcal{Z}=\bigcup_{j=1}^m \mathcal{Z}_j$ denote a partition of the sample space. Assume
    \begin{enumerate}
        \item $P_\lambda(\mathcal{Z}_j)=Q_\lambda(\mathcal{Z}_j)=p_j$ for all $\lambda$, and
        \item the conditional distributions $\mathds{1}_{\mathcal{Z}_j}dP_\lambda/p_j$ and $\mathds{1}_{\mathcal{Z}_j}dQ_\lambda/p_j$ (given an observation is in $\mathcal{Z}_j$) do not depend on $\lambda_l$ for $l\neq j$.
    \end{enumerate}
    For a prior distribution $\varpi$ over $\lambda$, let $\bar p=\int p_\lambda d\varpi(\lambda)$ and $\bar q=\int q_\lambda d\varpi(\lambda)$, and define
    \begin{align*}
        \delta_1 &= \max_{j\in \{1,\dots,m\}} \sup_\lambda \int_{\mathcal{Z}_j} \frac{(p_\lambda - \bar p)^2}{p_\lambda p_j}d\nu, \\
        \delta_2 &= \max_{j\in \{1,\dots,m\}} \sup_\lambda \int_{\mathcal{Z}_j} \frac{(q_\lambda - p_\lambda)^2}{p_\lambda p_j}d\nu, \\
        \delta_3 &= \max_{j\in \{1,\dots,m\}} \sup_\lambda \int_{\mathcal{Z}_j} \frac{(\bar q - \bar p)^2}{p_\lambda p_j}d\nu \\
    \end{align*}
    for a dominating measure $\nu$. If $\bar p/p_\lambda\leq b<\infty$ and $np_j\max (1,\delta_1,\delta_2)\leq b$ for all $j$, then
    \begin{equation*}
        H^2\left(\int P_\lambda^{\otimes n}d\omega(\lambda), \int Q_\lambda^{\otimes n}d\omega(\lambda)\right) \leq Cn\left\{ n\left(\max_{j\in \{1,\dots,m\}} p_j\right)\left(\delta_1\delta_2 + \delta_2^2\right) + \delta_3\right\}
    \end{equation*}
    for a constant $C$ that depends on $b$.
\end{lemma}

\paragraph{Case $\delta_n\lesssim\epsilon_n$:}
We want to verify the assumptions of the previous lemma. We first calculate the densities $p_\lambda, q_\lambda$, and $\bar p, \bar q$. We have
\begin{align*}
    p_\lambda(z) = &\mathds{1}(a\in S)\widehat\pi_a(x)\widehat\mu_a(x)^y\{1-\widehat\mu_a(x)\}^{1-y} + \mathds{1}(a\notin S, a\neq k)\widehat\pi_a(x)\widetilde\mu_a(x)^y\{1-\widetilde\mu_a(x)\}^{1-y} \\
    &+ \mathds{1}(a=k)\left\{ 1-\sum_{j=1}^{k-1}\widehat\pi_t(x) \right\}\widetilde\mu_k(x)^y\{1-\widetilde\mu_k(x)\}^{1-y}.
\end{align*}
In particular, $p_\lambda$ does not depend on $\lambda$, hence
\begin{equation*}
    \bar p = \int p_\lambda d\varpi(\lambda) = p_\lambda.
\end{equation*}
Further, we have
\begin{align*}
    q_\lambda(z)&=\mathds{1}(a\neq k)\pi_{a,\lambda}(x)\mu_{a,\lambda}(x)^y\{1-\mu_{a,\lambda}(x)\}^{1-y} + \mathds{1}(a=k)\left\{ 1-\sum_{t=1}^{k-1} \pi_{t,\lambda}(x) \right\}\widetilde\mu_a(x)^y\{1-\widetilde\mu_a(x)\}^{1-y} \\
    &= \mathds{1}(a\in S)\widehat\pi_a(x)\left\{1-\frac{\alpha}{\widehat\mu_a(x)}\sum_{j=1}^m \lambda_j B_j(x)\right\}\left\{\frac{\widehat\mu_a(x) + \beta\sum_{j=1}^m \lambda_j B_j(x)}{1-\frac{\alpha}{\widehat\mu_a(x)}\sum_{j=1}^m \lambda_j B_j(x)}\right\}^y\\
    &\cdot\left\{\frac{1-\widehat\mu_a(x) - \left(\beta+\frac{\alpha}{\widehat\mu_a(x)}\right)\sum_{j=1}^m \lambda_j B_j(x)}{1-\frac{\alpha}{\widehat\mu_a(x)}\sum_{j=1}^m \lambda_j B_j(x)}\right\}^{1-y} + \mathds{1}(a\notin S, a\neq k)\widehat\pi_a(x)\widetilde\mu_a(x)^y\{1-\widetilde\mu_a(x)\}^{1-y}\\
    &+\mathds{1}(a=k)\left(1-\sum_{t\in S}\widehat\pi_t(x)\left\{1-\frac{\alpha}{\widehat\mu_t(x)}\sum_{j=1}^m \lambda_j B_j(x) \right\}-\sum_{t\notin S, t\neq k}\widehat\pi_t(x)\right)\widetilde\mu_a(x)^y\{1-\widetilde\mu_a(x)\}^{1-y} \\
    &=\mathds{1}(a\in S)\widehat\pi_a(x)\left\{\widehat\mu_a(x) + \beta\sum_{j=1}^m \lambda_j B_j(x)\right\}^y\left\{1-\widehat\mu_a(x) - \left(\beta+\frac{\alpha}{\widehat\mu_a(x)}\right)\sum_{j=1}^m \lambda_j B_j(x)\right\}^{1-y}\\
    &+\mathds{1}(a\notin S, a\neq k)\widehat\pi_a(x)\widetilde\mu_a(x)^y\{1-\widetilde\mu_a(x)\}^{1-y}\\
    &+\mathds{1}(a=k)\left(1-\sum_{t\in S}\widehat\pi_t(x)\left\{1-\frac{\alpha}{\widehat\mu_t(x)}\sum_{j=1}^m \lambda_j B_j(x) \right\}-\sum_{t\notin S, t\neq k}\widehat\pi_t(x)\right)\widetilde\mu_a(x)^y\{1-\widetilde\mu_a(x)\}^{1-y}.
\end{align*}
Choose the uniform prior $\varpi(\lambda)=\frac{1}{2^m}\prod_j\mathds{1}(\lambda_j\in\{1,-1\})$. Then
\begin{align*}
    \bar q(z) &= \int q_\lambda(z) d\varpi(\lambda) \\
    &= \mathds{1}(a\in S)\widehat\pi_a(x)\widehat\mu_a(x)^y\{1-\widehat\mu_a(x)\}^{1-y} + \mathds{1}(a\notin S, a\neq k)\widehat\pi_a(x)\widetilde\mu_a(x)^y\{1-\widetilde\mu_a(x)\}^{1-y} \\
    &+ \mathds{1}(a=k)\left\{1-\sum_{t=1}^{k-1}\widehat\pi_t(x)\right\}\widetilde\mu_a(x)^y\{1-\widetilde\mu_a(x)\}^{1-y}.
\end{align*}
by using the fact that $q_\lambda$ is linear in $\sum_j \lambda_jB_j$ and $\int f(z) \sum_j \lambda_jB_j(x) d\varpi(\lambda)=0$ due to the choice of prior.

Now, let the partition of the lemma be given by $\mathcal{Z}_j=\{0,1\}\times\{1,\dots,k\}\times b_{2j-1}\cup b_{2j}$. Then, we have
\begin{equation*}
    p_j = P_\lambda(\mathcal{Z}_j) = Q_\lambda(\mathcal{Z}_j) = 2\mathrm{vol}(b_1) =\frac{1}{m}
\end{equation*}
by assumption on the size of the cubes, i.e., $\int \mathds{1}(x\in b_j(x))dx=\frac{1}{2m}$, and the fact that the marginal distribution of $P_\lambda$ and $Q_\lambda$ with respect to $X$ is the uniform distribution.
Hence, Assumption 1 of Lemma \ref{lemma-bound-hellinger-distance} is satisfied. Assumption 2 is also satisfied since $\sum_j \lambda_j B_j(x)$ depends only on $\lambda_s$ whenever $x\in\mathcal{Z}_s$ because $x\in\mathcal{Z}_s$ implies $B_j(x)=0$ for all $j\neq s$.

We can now characterize the constants $\delta_1,\delta_2$, and $\delta_3$ in Lemma \ref{lemma-bound-hellinger-distance}: Since $p_\lambda = \bar p$, we have $\delta_1=0$. Moreover, $\bar q = \bar p$, so also $\delta_3=0$. It remains to determine $\delta_2$. We have
\begin{align*}
    \delta_2 &= \max_{j\in \{1,\dots,m\}} \sup_\lambda \int_{\mathcal{Z}_j} \frac{(q_\lambda - p_\lambda)^2}{p_\lambda p_j} d\nu = m\max_{j\in \{1,\dots,m\}} \sup_\lambda \int_{\mathcal{Z}_j} \frac{(q_\lambda - p_\lambda)^2}{p_\lambda} d\nu.
\end{align*}
Now, note that
\begin{equation*}
    p_\lambda(z) \geq \mathds{1}(a\in S)\varepsilon \widehat\pi_a(x)+\mathds{1}(a=k)\widehat\pi_k(x)\widetilde\mu_k(x)^y\{1-\widetilde\mu_k(x)\}^{1-y}.
\end{equation*}
Then,
\begin{align*}
    \frac{(q_\lambda(z) - p_\lambda(z))^2}{p_\lambda(z)} &\leq \mathds{1}(a\in S)\frac{\widehat\pi_a(x)}{\varepsilon}\left\{ \beta^2\sum_{j=1}^m B_j(x)^2 \right\}^y\left\{\left(\beta + \frac{\alpha}{\widehat\mu_a(x)}\right)^2\sum_{j=1}^m B_j(x)^2\right\}^{1-y} \\
    &+ \mathds{1}(a=k)\left(\sum_{t\in S} \frac{\widehat\pi_t(x)}{\sqrt{\widehat\pi_k(x)}}\left\{ \frac{\alpha}{\widehat\mu_t(x)}\sum_{j=1}^m \lambda_j B_j(x)\right\}\right)^2\widetilde\mu_k(x)^{y}\{1-\widetilde\mu_k(x)\}^{1-y}\\
    &\leq \mathds{1}(a\in S) \frac{1}{\varepsilon}\beta^{2y}\left(2\beta^2 + 2\frac{\alpha^2}{\widehat\mu_a(x)^2}\right)^{1-y}\widehat\pi_a(x)\sum_{j=1}^m B_j(x)^2 \\
    &+ \mathds{1}(a=k)(1-\varepsilon)\frac{\alpha^2}{\varepsilon^2}\frac{1}{\varepsilon}\sum_{j=1}^m B_j(x)^2 \\
    &\leq \mathds{1}(a\in S) \frac{2}{\varepsilon^3}\left(\beta^2 + \alpha^2\right)\widehat\pi_a(x)\sum_{j=1}^m B_j(x)^2 + \mathds{1}(a=k)\frac{1-\varepsilon}{\varepsilon^3}\alpha^2\sum_{j=1}^m B_j(x)^2 \\
    &\leq \frac{2}{\varepsilon^3}\left( \alpha^2+ \beta^2\right)\left(\sum_{j=1}^m B_j(x)^2\right)\left( \mathds{1}(a\in S)\widehat\pi_a(x) + \mathds{1}(a=k) \right).
\end{align*}
Since the previous line does not depend on $\lambda$, when choosing the dominating measure $\nu$ as the Lebesgue measure $\mu$, we obtain
\begin{align*}
    \delta_2 &\leq m\max_{j\in \{1,\dots,m\}} \int_{\mathcal{Z}_j} \frac{(q_\lambda - p_\lambda)^2}{p_\lambda} d\nu(z) \\
    &\leq m \max_{j\in \{1,\dots,m\}} \int_{b_{2j-1}\cup b_{2j}}\sum_{a=1}^k\sum_{y=0}^1 \frac{2}{\varepsilon^3}\left( \alpha^2+ \beta^2\right)\sum_{r=1}^m B_r(x)^2 \left( \mathds{1}(a\in S)\widehat\pi_a(x) + \mathds{1}(a=k) \right)dx \\
    &\leq m \max_{j\in \{1,\dots,m\}} \int_{b_{2j-1}\cup b_{2j}} \frac{8}{\varepsilon^3}\left( \alpha^2+ \beta^2\right)\sum_{r=1}^m B_r(x)^2 dx\\
    &=\max_{j\in \{1,\dots,m\}} \frac{8}{\varepsilon^3}(\alpha^2+\beta^2)\cdot m\cdot\mathrm{vol}(b_{2j-1}\cup b_{2j}) \\
    &= \frac{8}{\varepsilon^3}(\alpha^2+\beta^2)\cdot m\cdot\frac{1}{m}\\
    &=\frac{8}{\varepsilon^3}(\alpha^2+\beta^2).
\end{align*}
Lastly, choose
\begin{equation*}
    m = \max\left\{n, 2Cn^2\left(\frac{8}{\varepsilon^3}\right)^2(\alpha^2+\beta^2)^2\right\},
\end{equation*}
where $C$ is the constant given by Lemma \ref{lemma-bound-hellinger-distance} if it was applied with $b=\max\left\{1, \frac{8}{\varepsilon^3}(C_\mu^2+C_\pi^2)\right\}$.
Now, we have $\bar p/p_\lambda = 1$ and
\begin{equation*}
    np_j\max(1,\delta_1,\delta_2) \leq \frac{n}{m}\max\left\{1, \frac{8}{\varepsilon^3}(C_\mu^2+C_\pi^2)\right\}\leq \max\left\{1, \frac{8}{\varepsilon^3}(C_\mu^2+C_\pi^2)\right\}.
\end{equation*}
Therefore, we can apply Lemma \ref{lemma-bound-hellinger-distance} with $b$ defined as above and obtain
\begin{equation*}
    H^2\left(\int P_\lambda^{\otimes n}d\omega(\lambda), \int Q_\lambda^{\otimes n}d\omega(\lambda)\right) \leq C\frac{n^2}{m}\frac{64}{\varepsilon^6}(\alpha^2+\beta^2)^2\leq \frac{1}{2}.
\end{equation*}

\paragraph{Case $\epsilon_n\lesssim\delta_n$:} First compute that
\begin{align*}
    &q_\lambda(z) \\
    &= \mathds{1}(a\in S)\left\{ 1 + \frac{\beta}{\widehat\mu_a(x)}\sum_{j=1}^m \lambda_jB_j(x) - \alpha\beta \right\}\left\{  \widehat\pi_a(x) + \alpha\widehat\pi_a(x)\widehat\mu_a(x)\sum_{j=1}^m \lambda_j B_j(x) \right\}\\
    &\left( \frac{\widehat\mu_a(x) }{1+\frac{\beta}{\widehat\mu_a(x)}\sum_{j=1}^m \lambda_jB_j(x) - \alpha\beta} \right)^y\left( 1-\frac{\widehat\mu_a(x) }{1+\frac{\beta}{\widehat\mu_a(x)}\sum_{j=1}^m \lambda_jB_j(x) - \alpha\beta} \right)^{1-y} \\
    &+ \mathds{1}(a\notin S, a\neq k)\widehat\pi_a(x)\widetilde\mu_a(x)^y\left\{ 1-\widetilde\mu_a(x) \right\}^{1-y} \\
    &+\mathds{1}(a=k)\left( 1-\sum_{t\in S}\left\{1+\frac{\beta}{\widehat\mu_t(x)}\sum_{j=1}^m \lambda_j B_j(x) - \alpha\beta\right\}\left\{ \widehat\pi_t(x) + \alpha\widehat\pi_t(x)\widehat\mu_t(x)\sum_{j=1}^m \lambda_j B_j(x) \right\} - \sum_{t\notin S, t\neq k}\widehat\pi_t(x) \right)\\
    &\widetilde\mu_t(x)^y\left\{ 1-\widetilde\mu_t(x) \right\}^{1-y} \\
    &=\mathds{1}(a\in S)\left\{  \widehat\pi_a(x) + \alpha\widehat\pi_a(x)\widehat\mu_a(x)\sum_{j=1}^m \lambda_j B_j(x) \right\}\widehat\mu_a(x)^y\left\{ 1 + \frac{\beta}{\widehat\mu_a(x)}\sum_{j=1}^m \lambda_jB_j(x) - \alpha\beta - \widehat\mu_a(x)\right\}^{1-y}\\
    &+\mathds{1}(a\notin S, a\neq k)\widehat\pi_a(x)\widetilde\mu_a(x)^y\left\{ 1-\widetilde\mu_a(x) \right\}^{1-y}\\
    &+\mathds{1}(a=k)\Bigg( 1-\sum_{t=1}^{k-1}\widehat\pi_t(x) - \sum_{t\in S}\left\{\frac{\beta}{\widehat\mu_t(x)}\sum_{j=1}^m \lambda_j B_j(x) - \alpha\beta\right\}\left\{ \widehat\pi_t(x) + \alpha\widehat\pi_t(x)\widehat\mu_t(x)\sum_{j=1}^m \lambda_j B_j(x) \right\} \\
    &- \sum_{t\in S} \alpha\widehat\pi_t(x)\widehat\mu_t(x)\sum_{j=1}^m \lambda_j B_j(x) \Bigg)\widetilde\mu_a(x)^y\left\{ 1-\widetilde\mu_a(x) \right\}^{1-y}.
\end{align*}
Therefore,
\begin{align*}
    \bar q_\lambda(z) &= \mathds{1}(a\in S)\widehat\pi_a(x)\widehat\mu_a(x)^y\left\{ 1-\widehat\mu_a(x) \right\}^{1-y} \\
    &+ \mathds{1}(a\notin S, a\neq k)\widehat\pi_a(x)\widetilde\mu_a(x)^y\left\{ 1-\widetilde\mu_a(x) \right\}^{1-y}\\
    &+ \mathds{1}(a=k)\left( 1-\sum_{t=1}^{k-1} \widehat\pi_t(x) \right)\widetilde\mu_a(x)^y\left\{ 1-\widetilde\mu_a(x) \right\}^{1-y}.
\end{align*}
Hence, we again obtain $\delta_1=0$ and $\delta_3=0$. It remains to determine $\delta_2$. Using again that
\begin{equation*}
    p_\lambda(z) \geq \mathds{1}(a\in S)\varepsilon \widehat\pi_a(x)+\mathds{1}(a=k)\widehat\pi_k(x)\widetilde\mu_k(x)^y\{1-\widetilde\mu_k(x)\}^{1-y},
\end{equation*}
we obtain
\begin{align*}
    &\frac{(q_\lambda(z) - p_\lambda(z))^2}{p_\lambda(z)} \\
    &\leq \mathds{1}(a\in S)\frac{\widehat\pi_a(x)}{\varepsilon}\Bigg[ \left\{ 1+\alpha\widehat\mu_a(x)\sum_{j=1}^m\lambda_j B_j(x) \right\}\widehat\mu_a(x)^y\left\{ 1+\frac{\beta}{\widehat\mu_a(x)}\sum_{j=1}^m \lambda_j B_j(x) - \alpha\beta - \widehat\mu_a(x) \right\}^{1-y} \\
    &- \widehat\mu_a(x)^y\left\{ 1-\widehat\mu_a(x) \right\}^{1-y} \Bigg]^2 \\
    &+\mathds{1}(a=k)\widetilde\mu_a(x)^{y}\left\{ 1-\widetilde\mu_a(x) \right\}^{1-y}\frac{1}{\widehat\pi_k(x)}\Bigg( - \sum_{t\in S}\left\{\frac{\beta}{\widehat\mu_t(x)}\sum_{j=1}^m \lambda_j B_j(x) - \alpha\beta\right\}\\
    &\left\{ \widehat\pi_t(x) + \alpha\widehat\pi_t(x)\widehat\mu_t(x)\sum_{j=1}^m \lambda_j B_j(x) \right\} - \sum_{t\in S} \alpha\widehat\pi_t(x)\widehat\mu_t(x)\sum_{j=1}^m \lambda_j B_j(x) \Bigg)^2\\
    &\leq \mathds{1}(a\in S) \frac{\widehat\pi_a(x)}{\varepsilon}\left[ \left(\sum_{j=1}^m \lambda_j B_j(x)\right)\left\{ \frac{\beta}{\widehat\mu_a(x)} + \alpha\widehat\mu_a(x) - \alpha^2\beta\widehat\mu_a(x) - \alpha\widehat\mu_a(x)^2 \right\} \right]^2 \\
    &+\mathds{1}(a=k)(1-\varepsilon)\frac{1}{\varepsilon}\left(-\left(\sum_{j=1}^m \lambda_j B_j(x)\right)\sum_{t\in S}\left(\frac{\beta\widehat\pi_t(x)}{\widehat\mu_t(x)} - \alpha^2\beta \widehat\pi_t(x)\widehat\mu_t(x)+\alpha\widehat\pi_t(x)\widehat\mu_t(x)\right) \right)^2\\
    &\leq \mathds{1}(a\in S)\frac{1}{\varepsilon}\left( \frac{\beta}{\varepsilon} + \alpha \right)^2\widehat\pi_a(x)\sum_{j=1}^m B_j(x)^2+\mathds{1}(a=k) (1-\varepsilon)\frac{1}{\varepsilon}\left( \frac{\beta}{\varepsilon} + \alpha \right)^2\sum_{j=1}^m B_j(x)^2 \\
    &\leq \frac{2}{\varepsilon^3}(\alpha^2 + \beta^2)\left(\sum_{j=1}^m B_j(x)^2\right) \left( \mathds{1}(a\in S)\widehat\pi_a(x) + \mathds{1}(a=k) \right).
\end{align*}
Consequently, $\delta_2$ is the same as for the case $\delta_n\lesssim\epsilon_n$, which means that we can choose the same $m$ to bound the Hellinger distance.

\subsubsection{Final Minimax Lower Bound}

From Section \ref{sec-functional-separation-sa}, we know that
\begin{align*}
    &\sqrt{\sum_{a=1}^k \frac{1}{k}\left\{\psi_a(P_\lambda) - \psi_a(Q_{\lambda'})\right\}^2} \\
    &\geq \sqrt{\frac{s}{k}}\frac{\varepsilon}{2}\begin{cases}
        \left( \frac{\varepsilon^2 C_\pi}{\max\{1,2C_\pi, 2D C_\pi/C_\mu\}}\right)\left(
        \frac{C_\mu}{4}\right)\delta_n\epsilon_n & \text{for }\delta_n\lesssim \epsilon_n\\
        \left(\frac{\varepsilon C_\pi}{4\max\{1,C_\pi\}}\right)
        \left(\frac{\varepsilon^2 C_\mu}{4\max\{1,C_\mu, DC_\mu/C_\pi, DC_\mu\}}\right)\delta_n\epsilon_n & \text{for }\epsilon_n\lesssim \delta_n
    \end{cases}\\
    &\geq \sqrt{C'}\sqrt{\frac{s}{k}}\delta_n\epsilon_n=:\Delta
\end{align*}
for
\begin{align*}
    C'=\frac{\varepsilon^2}{4}\min\Bigg\{&\left( \frac{\varepsilon^2 C_\pi}{\max\{1,2C_\pi, 2D C_\pi/C_\mu\}}\right)^2\left(
        \frac{C_\mu}{4}\right)^2;\\
        &\left(\frac{\varepsilon C_\pi}{4\max\{1,C_\pi\}}\right)^2
        \left(\frac{\varepsilon^2 C_\mu}{4\max\{1,C_\mu, DC_\mu/C_\pi, DC_\mu\}}\right)^2\Bigg\},
\end{align*}
and from Section \ref{sec-distributional-distance-sa} that
\begin{equation*}
    H^2\left(\int P_\lambda^{\otimes n}d\omega(\lambda), \int Q_\lambda^{\otimes n}d\omega(\lambda)\right) \leq \frac{1}{2}=:\alpha.
\end{equation*}
We can now apply Lemma \ref{lem:minimax} with $d(x,y)=
\sqrt{\frac{1}{k}\sum_{a=1}^k (x_a-y_a)^2}$ for this $\Delta$ and $\alpha$ and any $P\in\mathcal{P}$, which gives
\begin{align*}
    \inf_{\widehat\psi} \sup_{P\in\mathcal{P}} \E_P\left\{\sum_{a=1}^k\varpi_a(P)\left(\widehat\psi_a - \psi_a(P)\right)^2\right\} &\geq \min\{C_1',\varepsilon'\}\left(\frac{\Delta}{2}\right)^2\left( \frac{1 - \sqrt{\alpha(1-\alpha/4)}}{2} \right)\\
    &\geq \min\{C_1',\varepsilon'\}C'\frac{4-\sqrt{7}}{8}\cdot \frac{s}{k}\left(\delta_n\epsilon_n\right)^2
\end{align*}
by using the fact that $\varpi_a(P)\geq \min\{C_1',\varepsilon'\}/k$ uniformly across all $a\in\{1,\dots,k\}$ and $P\in\mathcal{P}$.
This proves Theorem \ref{thm:minimax-lb-structure-agnostic}.

\end{document}